\def\isarxivversion{1} 

\ifdefined\isarxivversion
\documentclass[11pt]{article}
\else
\documentclass[anon,12pt]{colt2021}
\fi

\ifdefined\isarxivversion
\usepackage{amsfonts}
\fi

\usepackage{algorithm}
\usepackage{array}
\usepackage{multirow}
\usepackage{color}
\usepackage[english]{babel}
\usepackage{graphicx}
\usepackage{wrapfig}%,epsfig}
\usepackage{epstopdf}
\usepackage{url}
\usepackage{graphicx}
\usepackage{color}
\usepackage{epstopdf}
\usepackage{algpseudocode}
\usepackage{scrextend}
\usepackage[T1]{fontenc}
\usepackage{bbm}
\usepackage{comment}
\usepackage{tikz}
\usepackage{tikz}
\usepackage{hyperref}
\hypersetup{colorlinks=true,citecolor=blue,linkcolor=red}

\usetikzlibrary{arrows}

\graphicspath{{./figs/}}

\ifdefined\isarxivversion
\usepackage[margin=1in]{geometry}
\else 

\fi

\definecolor{b2}{RGB}{51,153,255}
\definecolor{mygreen}{RGB}{80,180,0}
\definecolor{trz}{HTML}{4169e1}
\definecolor{shunhua}{RGB}{0,180,80}
\definecolor{mycrimson}{RGB}{165,28,48}

\newcommand{\Zhao}[1]{{\color{b2}[Zhao: #1]}}

\newcommand{\Baihe}[1]{{\color{mygreen}[Baihe: #1]}}
\newcommand{\Runzhou}[1]{{\color{trz}[Runzhou: #1]}}
\newcommand{\Shunhua}[1]{{\color{red}[Shunhua: #1]}}

\newcommand{\Ruizhe}[1]{{\color{mycrimson}[Ruizhe: #1]}}
\renewcommand{\Zhao}[1]{{}}
\renewcommand{\Baihe}[1]{{}}
\renewcommand{\Ruizhe}[1]{{}}
\renewcommand{\Shunhua}[1]{{}}
\renewcommand{\Runzhou}[1]{{}}

\ifdefined\isarxivversion
\usepackage{amsmath}
\usepackage{amsthm}
\usepackage{subfig}

\newtheorem{theorem}{Theorem}[section]
\newtheorem{lemma}[theorem]{Lemma}
\newtheorem{definition}[theorem]{Definition}

\newtheorem{corollary}[theorem]{Corollary}

\newtheorem{remark}[theorem]{Remark}

\else
\fi

\newtheorem{assumption}[theorem]{Assumption}

\newtheorem{fact}[theorem]{Fact}

\newcommand{\wt}{\widetilde}

\newcommand{\N}{\mathcal{N}}
\newcommand{\R}{\mathbb{R}}

\newcommand{\LHS}{\mathrm{LHS}}

\renewcommand{\d}{\mathrm{d}}

\newcommand{\tr}{\mathrm{tr}}
\newcommand{\new}{\mathrm{new}}

\newcommand{\midd}{\mathrm{mid}}
\newcommand{\diag}{\mathrm{diag}}

\newcommand{\g}{{\mathsf{g}}}
\newcommand{\A}{{\mathsf{A}}}
\newcommand{\Tmat}{{\cal T}_\mathrm{mat}}
\newcommand{\vect}{\mathrm{vec}}

\newcommand{\supp}{\mathrm{supp}}

\makeatletter
\newcommand*{\RN}[1]{\expandafter\@slowromancap\romannumeral #1@}
\makeatother

\begin{document}

\ifdefined\isarxivversion
\title{Solving SDP Faster: A Robust IPM Framework and Efficient Implementation}

\else

\fi

\ifdefined\isarxivversion
\author{
	Baihe Huang\thanks{\texttt{baihehuang@pku.edu.cn}. Peking University.}
	\and
	Shunhua Jiang\thanks{\texttt{sj30050@columbia.edu}. Columbia University.}
	\and
	Zhao Song\thanks{\texttt{zsong@adobe.com}. Adobe Research.}
	\and
	Runzhou Tao\thanks{\texttt{runzhou.tao@columbia.edu}. Columbia University.}
	\and
	Ruizhe Zhang\thanks{\texttt{ruizhe@utexas.edu}. The University of Texas at Austin.}
}

\date{}
\else

\fi

\ifdefined\isarxivversion
 \begin{titlepage}
     \maketitle
     \begin{abstract}
         
This paper introduces a new robust interior point method analysis for semidefinite programming (SDP). This new robust analysis can be combined with either logarithmic barrier or hybrid barrier.

Under this new framework, we can improve the running time of semidefinite programming (SDP) with variable size $n \times n$ and $m$ constraints up to $\epsilon$ accuracy.

We show that for the case $m = \Omega(n^2)$, we can solve SDPs in $m^{\omega}$ time. This suggests solving SDP is nearly as fast as solving the linear system with equal number of variables and constraints. This is the first result that tall dense SDP can be solved in the nearly-optimal running time, and it also improves the state-of-the-art SDP solver [Jiang, Kathuria, Lee, Padmanabhan and Song, FOCS 2020].

In addition to our new IPM analysis, we also propose a number of techniques that might be of further interest, such as, maintaining the inverse of a Kronecker product using lazy updates, a general amortization scheme for positive semidefinite matrices.

     \end{abstract}
     \thispagestyle{empty}
 \end{titlepage}
 \else
 
 \maketitle
\begin{abstract}

This paper introduces a new robust interior point method analysis for semidefinite programming (SDP). This new robust analysis can be combined with either logarithmic barrier or hybrid barrier.

Under this new framework, we can improve the running time of semidefinite programming (SDP) with variable size $n \times n$ and $m$ constraints up to $\epsilon$ accuracy.

We show that for the case $m = \Omega(n^2)$, we can solve SDPs in $m^{\omega}$ time. This suggests solving SDP is nearly as fast as solving the linear system with equal number of variables and constraints. This is the first result that tall dense SDP can be solved in the nearly-optimal running time, and it also improves the state-of-the-art SDP solver [Jiang, Kathuria, Lee, Padmanabhan and Song, FOCS 2020].

In addition to our new IPM analysis, we also propose a number of techniques that might be of further interest, such as, maintaining the inverse of a Kronecker product using lazy updates, a general amortization scheme for positive semidefinite matrices.

\end{abstract}

 \fi
 
 \iffalse
\newpage

\pagenumbering{roman}

{\hypersetup{linkcolor=black}
\tableofcontents
}
\newpage

\pagenumbering{arabic}
\setcounter{page}{1}
\fi

\section{Introduction}
Semidefinite programming (SDP) optimizes a linear objective function over the intersection of the positive semidefinite (PSD) cone with an affine space. SDP is of great interest both in theory and in practice. Many problems in operations research, machine learning, and theoretical computer science can be modeled or approximated as semidefinite programming problems. In machine learning, SDP has applications in adversarial machine learning \cite{rsl18}, learning structured distribution \cite{clm20}, sparse PCA \cite{aw08,dejl07}, robust learning \cite{dkklms16,dhl19,jlt20}. %, sparse matrix factorization \cite{csz20}. 
In theoretical computer science, SDP has been used in approximation algorithms for max-cut \cite{gw94}, coloring $3$-colorable graphs \cite{kms94}, and sparsest cut \cite{arv09}, quantum complexity theory \cite{jjuw11}, robust learning and estimation \cite{cg18,cdg19,cdgw19}, graph sparsification \cite{ls17}, algorithmic discrepancy and rounding \cite{bdg16,bg17,b19}, sum of squares optimization \cite{bs16,fkp19}, terminal embeddings \cite{cn21}, and matrix discrepancy \cite{hrs21}.

SDP is formally defined as follows: 
\begin{definition}[Semidefinite programming] \label{def:sdp_primal}
Given symmetric\footnote{We can without loss of generality assume that $C,A_1,\cdots,A_m$ are symmetric. Given any $A \in \R^{n \times n}$, we have $\sum_{i, j} A_{ij} X_{ij} = \sum_{i, j} A_{ij} X_{ji} = \sum_{i, j} (A^\top)_{ij} X_{ij}$ since $X$ is symmetric, so we can replace $A$ with $(A+ A^\top)/2$.} matrices $C, A_1, \cdots, A_m \in \mathbb{R}^{n \times n}$ and a vector $b \in \R^m$, the goal is to solve the following optimization problem:
\begin{align}
\label{eq:sdp_primal}
\max_{X \in \R^{n \times n}} ~ \langle C, X \rangle %\notag \\
\textup{ subject to } ~ \langle A_i, X \rangle = b_i, ~~\forall i \in [m], %\\
~ X \succeq 0,% \notag
\end{align}
where $\langle A,B\rangle := \sum_{i,j} A_{i,j} B_{i,j}$ is the matrix inner product. 
\end{definition}

The input size of an SDP instance is $m n^2$, since there are $m$ constraint matrices each of size $n \times n$. The well-known linear programming (LP) is a simpler case than SDP, where $X \succeq 0$ and $C, A_1, \cdots, A_m$ are restricted to be $n \times n$ diagonal matrices. The input size of an LP instance is thus $m n$.

Over the last many decades, there are three different lines of high accuracy SDP solvers (with logarithmic accuracy dependence in the running time). The first line of work is using the cutting plane method, such as \cite{s77,yn76,k80,kte88,nn89,v89,bv02,km03,lsw15,jlsw20}. This line of work uses $m$ iterations, and each iteration uses some SDP-based oracle call. The second line of work is using interior point method (IPM) and log barrier function such as \cite{nn92,jklps20}. The third line of work is using interior point method and hybrid barrier function such as \cite{nn94,a00}. 

Recently, a line of work uses robust analysis and dynamic maintenance to speedup the running time of linear programming \cite{cls19,b20,blss20,jswz21,b21}. One major reason made solving SDP much more harder than solving linear programming is: in LP the slack variable is a vector(can be viewed as a diagonal matrix), and in SDP {\it the slack variable is a positive definite matrix}. Due to that reason, the gradient/Hessian computation requires some complicated and heavy calculations based on the Kronecker product of matrices, while LP only needs the basic matrix-matrix product \cite{v89_lp,cls19,jswz21}. Therefore, handling the errors in each iteration and maintaining the slack matrices are way more harder in SDP. Thus, we want to ask the following question:
\begin{center}
    {\it Can we efficiently solve SDP \underline{without} computing exact gradient, Hessian, and Newton steps?}
\end{center}

In this work, we will answer the above question by introducing new framework for both IPM analysis and variable maintenance. For IPM analysis, we build a robust IPM framework for arbitrary barrier functions that supports errors in computing gradient, Hessian, and Newton steps. For variable maintenance, we provide a general amortization method that gives improved guarantees on reducing the computational complexity by lazily updating the Hessian matrices.

For solving SDP using IPM with log barrier, the current best  algorithm (due to Jiang, Kathuria, Lee, Padmanabhan and Song \cite{jklps20}) runs in $O(\sqrt{n} (mn^2 + m^{\omega} + n^{\omega}))$ time. Since the input size of SDP is $m n^2$, ideally we would want an SDP algorithm that runs in $O(mn^2 + m^{\omega} + n^{\omega})$ time, which is roughly the running time to solve linear systems\footnote{We note that a recent breakthrough result by Peng and Vempala \cite{pv21} showed that a sparse linear system can be solved  faster than matrix multiplication. However, their algorithm essentially rely on the sparsity of the problems. And it is still widely believed that general linear system requires matrix multiplication time.}. The current best algorithms are still at least a $\sqrt{n}$ factor away from the optimal.

Inspired by the result \cite{cls19} which solves LP in the current matrix multiplication time, a natural and fundamental question for SDP is 
\begin{center}
{\it
    Can we solve SDP in the current matrix multiplication time?
}
\end{center}
More formally, for the above formulation of SDP (Definition~\ref{def:sdp_primal}), is that possible to solve it in $mn^2+m^{\omega} + n^{\omega}$ time? In this work, we give a positive answer to this question by using our new techniques. For the tall dense SDP where $m = \Omega(n^2)$, our algorithm runs in $m^{\omega} + m^{2+1/4}$ time, which matches the current matrix multiplication time. The tall dense SDP finds many applications and is one of the two predominant cases in \cite{jklps20}\footnote{See Table~1.2 in \cite{jklps20} and Section~\ref{sec:related_work}.}. 
This is the first result that shows SDP can be solved as fast as solving linear systems.%, and we improve it in Table~\ref{table:runtime}).

Finally, we also show that our techniques and framework are quite versatile and can be used to directly speedup the SDP solver via the hybrid barrier \cite{nn89,a00}.

%\vspace{-3mm}
\paragraph{Our results.}
%\label{subsec:results}
We present the simplified version of our main result in the following theorem. The formal version can be found in Theorem~\ref{thm:general_SDP_formal}.
%\vspace{-1mm}
\begin{theorem}[Main result, informal version of Theorem~\ref{thm:general_SDP_formal}]\label{thm:tall_SDP_informal}
%Let $m = \Omega(n^2)$. 
For $\epsilon$-accuracy, there is a classical algorithm that solves a general SDP instance with variable size $n \times n$ and $m$ constraints in time\footnote{We use $O^\ast(\cdot)$ to hide $n^{o(1)}$ and $\log^{O(1)}(mn/\epsilon)$ factors, and $\Tilde{O}(\cdot)$ to hide $\log^{O(1)}(mn/\epsilon)$ factors.} %\footnote{We use $O^*$ to hide $n^{o(1)}$ and $\log^{O(1)}(n/\epsilon)$ factors and $\wt{O}$ to hide $\log^{O(1)}(n/\epsilon)$ factors, where $\epsilon$ is the accuracy parameter. The running time actually has a $\sqrt{n} \cdot m^2$ term. For simplicity we write it as $m^{2+1/4}$ since we are considering the case $m \geq n^2$.} 
$
    O^\ast( ( \sqrt{n} ( m^{2} + n^4 )+ m^{\omega} + n^{2\omega} ) \log(1/\epsilon)),
$
where $\omega$ is the exponent of matrix multiplication.

In particular, for $m = \Omega(n^2)$, our algorithm takes matrix multiplication time $m^{\omega}$ for current $\omega \approx 2.373$.
\end{theorem}
%\vspace{-1mm}
\begin{remark}
For any $m \geq n^{2 - 0.5/ \omega} \approx n^{1.79}$ with current $\omega \approx 2.37286$ \cite{w12,l14,aw21}, our algorithm runs faster than \cite{jklps20}. %See Remark~\ref{remark:running_time}. 
\end{remark}

Theorem~\ref{thm:tall_SDP_informal} and \cite{jklps20} are focusing on the log barrier method for solving SDP. However, the area of speeding up the hybrid barrier-based SDP solver is quite blank. We also improve the state-of-the-art implementation of the hybrid barrier-based SDP solver \cite{nn89,a00} in all parameter regimes. See Section~\ref{sec:intro_hybrid} and Theorem~\ref{thm:vol_intro} for more details.

%\vspace{-2mm}

\paragraph{Roadmap.} In Section~\ref{sec:previous_technique}, we review the previous approaches for solving SDP and discuss their bottlenecks. In Section~\ref{sec:robust_central_path}, we introduce our robust framework for IPM. In Section~\ref{sec:our_technique}, we show our main techniques and sketch the proof of our main result (Theorem~\ref{thm:tall_SDP_informal}). In Section~\ref{sec:intro_hybrid}, we overview the approach of applying our robust framework to speedup the hybrid barrier-based SDP solver. Related works are provided in Section~\ref{sec:related_work}. We define our notations and include several useful tools in Section~\ref{sec:preli}. In Section~\ref{sec:alg_and_result}, we give the formal version of our algorithm and the main theorem, where the proof is given Section~\ref{sec:correct} and \ref{sec:time}. Our general robust IPM framework is displayed in Section~\ref{sec:correct_previous}. Our fast implementation of the hybrid barrier-based SDP solver can be found in Section~\ref{sec:hybrid_barrier}.

%\vspace{-2mm}
%\vspace{-2mm}
\section{An Overview of Previous Techniques}\label{sec:previous_technique}
%\vspace{-2mm}

Under strong duality, the primal formulation of the SDP in Eq.~\eqref{eq:sdp_primal} is equivalent to the following dual formulation: 
\begin{definition}[Dual problem]\label{def:sdp_dual}
Given symmetric matrices $C,A_1,\dots,A_m\in \R^{n \times n}$ and $b_i\in \R$ for all $i \in [m]$, the goal is to solve the following convex optimization problem:
\begin{align}\label{eq:sdp_dual}
    \min_{y \in \R^m} ~ b^\top y ~~~ \mathrm{~subject~to} ~ S = \sum_{i=1}^m y_i A_i - C, 
    ~~~~ S \succeq 0. 
\end{align}
\end{definition}
%\vspace{-5mm}

%Using interior point methods to solve SDP is fast both in theory and practice.
Interior point methods (IPM) solve the above problem by (approximately) following a central path in the feasible region $\{ y \in \R^m : S = \sum_{i=1}^m y_i A_i - C \succeq 0 \}$. As a rich subclass of IPM, barrier methods \cite{nn92,a00} define a point on the central path as the solution to the following optimization problem parametrized by $\eta > 0: ~ \min_{y \in \R^m } f_{\eta} (y)$ where 
\begin{align}\label{eq:augmented_obj}
    f_{\eta}(y) : = \eta \cdot \langle b , y \rangle + \phi(y)
\end{align}
is the augmented objective function and $\phi : \R^m \rightarrow \R$ is a barrier function\footnote{The choice of the barrier function leads to different numbers of iterations. Nesterov and Nemirovski \cite{nn92} utilize the log barrier function $\phi_{\log}(y) = - \log \det ( S ) $ which guarantees convergence in $\wt{O}(\sqrt{n})$ iterations. Anstreicher \cite{a00} uses the Hybrid barrier $\phi_{\mathrm{hybrid}}(y) = 225(n/m)^{1/2} \cdot \left(\phi_{\mathrm{vol}}(y) + \phi_{\log}(y) \cdot (m-1)/(n-1)\right)$ where $\phi_{\mathrm{vol}}(y)$ is the volumetric barrier $\phi_{\mathrm{vol}}(y) = \frac{1}{2}\log\det (\nabla^2 \phi_{\log}(y))$. Hybrid barrier guarantees convergence in $\wt{O}((mn)^{1/4})$ iterations.} that restricts $y$ to the feasible region since $\phi(y)$ increases to infinity when $y$ approaches the boundary of the feasible region. Barrier methods usually start with an initial feasible $y$ for a small $\eta$, and increase $\eta$ in each iteration until $y$ is close to the optimal solution of the SDP. In short-step barrier methods with log barrier, $\eta^\new = (1+1/\sqrt{n})\eta$. It takes a Newton step $ -H(y)^{-1} g(y,\eta)$ in each iteration to keep $y$ in the proximity of the central path. Here $g(y,\eta)$ and $H(y) $ are the gradient and the Hessian of $f_{\eta}(y)$.

%\vspace{-2mm}
\subsection*{Techniques and bottlenecks of existing algorithms}

Fast solvers of SDP include the cutting plane method and interior point method. 
The fastest known algorithms for SDP based on the cutting plane method \cite{lsw15,jlsw20} have $m$ iterations and run in $O^*(m(mn^2+n^\omega+m^2))$ time. 
The fastest known algorithm for SDP based on the interior point method \cite{jklps20} has $\sqrt{n}$ iterations and runs in $O^*(\sqrt{n}(mn^2+n^\omega+m^\omega))$ time.
In most applications of SDP where $m \geq n$, interior point method of \cite{jklps20} runs faster. In the following we briefly discuss the techniques and bottlenecks of interior point methods.
%\vspace{-2mm}
\paragraph{Central path.}
Interior point method updates the dual variable $y$ by Newton step $-H(y)^{-1}g(y,\eta)$ to keep it in the proximity of central path.
This proximity is measured by the potential function $\|H(y)^{-1}g(y,\eta)\|_{H(y)}$.\footnote{For symmetric PSD matrix $A$, let $\| x \|_A=\sqrt{x^\top A x}$ denote matrix norm of $x$.}
In classical interior point literature, this potential function is well controlled by taking exact Newton step (see e.g. \cite{r01}).
%Due to computational considerations, 
\cite{jklps20} relaxes this guarantee and allows PSD approximation to the Hessian matrix $H(y)$. However, their convergence also relies on exact computation of slack matrix $S$ and gradient $g$. This leads to a $mn^{2.5}$ term in their running time. %In fact, In Section~\ref{sec:robust_central_path}, we analyze a more robust the potential perturbations of slack matrix $S$ and gradient $g$.

%\vspace{-2mm}
\paragraph{Amortization techniques.}
\cite{jklps20} keeps a PSD approximation $\wt{H}$ of the Hessian $H$ and updates $\wt{H}$ by a low rank matrix in each iteration. The running time of this low rank update is then controlled by a delicate amortization technique. This technique also appears in linear programming \cite{cls19} and empirical risk minimization \cite{lsz19}. \cite{jklps20} brings this technique to SDP, and costs $n^{0.5}m^\omega$ time in computing the inverse of Hessian matrix. When $m$ becomes larger, this term dominates the complexity and becomes undesirable.
%\vspace{-3mm}

%\vspace{-1mm}
\section{The Robust SDP Framework}\label{sec:robust_central_path}

\begin{algorithm}[!ht]
\caption{The general robust barrier method framework for SDP.
}
\label{alg:robust_ipm:intro}
%\small
\begin{algorithmic}[1]
\Procedure{\textsc{GeneralRobustSDP}}{$\mathsf{A} \in \R^{m \times n^2}$, $b \in \R^m$, $C \in \R^{n \times n}$} 
%\State \Comment{Initialization} 
\State Choose $\eta$ and $T$
%\State $\eta \leftarrow \frac{1}{n+2}$, ~~$T \leftarrow \frac{40}{\epsilon_N} \sqrt{n} \log (\frac{n}{\epsilon})$
\State Find initial feasible dual vector $y \in \R^m$ %according to Lemma~\ref{lem:initialization}
\Comment{{\bf Condition 0} in Lemma~\ref{lem:invariant_newton:intro}}
\For {$t = 1 \to T$} {\bf do} \Comment {Iterations of approximate barrier method}
	\State $\eta^{\new} \leftarrow \eta \cdot (1 + \frac{\epsilon_N}{20 \sqrt{\theta}})$
	\State $\wt{S} \leftarrow \textsc{ApproxSlack} ()$ %\Comment{{\bf Condition 1} in Lemma~\ref{lem:invariant_newton:intro}}
	\State $\wt{H} \leftarrow \textsc{ApproxHessian}()$ \Comment{{\bf Condition 1} in Lemma~\ref{lem:invariant_newton:intro}}
	\State $\wt{g} \leftarrow \textsc{ApproxGradient}( )$ \Comment{{\bf Condition 2} in Lemma~\ref{lem:invariant_newton:intro}}
	\State $\wt{\delta}_y \leftarrow \textsc{ApproxDelta}()$ \Comment{{\bf Condition 3} in Lemma~\ref{lem:invariant_newton:intro}}
	\State $y^{\new} \leftarrow y + \delta_y$
	\State $y \leftarrow y^{\new}$ \Comment{Update variables}
\EndFor
\EndProcedure
\end{algorithmic}
\end{algorithm}

In section~\ref{sec:robust_central_path}, we introduce our robust SDP framework. This framework works for general barrier functions and finds applications in both Algorithm~\ref{alg:ss_ipm_woodbury} and Algorithm~\ref{alg:sdp_vol}-\ref{alg:sdp_vol_cont}. We consider self-concordant barrier function $\phi$ with complexity $\theta$ (Definition~\ref{def:self_concordant_complexity})\footnote{The barrier function being ``self-concordant'' is a key assumption in the interior-point method \cite{n88,n88a,nn89,r01}. It is also useful in many optimization tasks \cite{hil14,nar16,llv20}.}. 
For the regularized objective $f_\eta$ in Eq.~\eqref{eq:augmented_obj}, we define the gradient $\g: \R^m \times \R \rightarrow \R^m$ as 
\begin{align*}
    \g(y,\eta) = \eta \cdot b - \nabla \phi(y) .
\end{align*}

Interior point method takes Newton step $(\nabla^2 \phi(y))^{-1}g(y,\eta)$ and guarantees the variables in the proximity of the central path by bounding the potential function $\Phi(z,y,\eta) = \|\g(y,\eta)\|_{(\nabla^2 \phi(z))^{-1}}$.
In practical implementations, there are perturbations in the Newton step due to errors in slack matrix $S$, gradient $\g(y,\eta)$, Hessian matrix $\nabla^2 \phi(y)$ and Newton step $(\nabla^2 \phi(y))^{-1} \cdot \g(y,\eta)$. Many fast algorithms maintain approximations to these quantities to reduce the running time. 
We propose a more general robust framework (compared with \cite{r01,jklps20}) which captures all these errors. 
We show that as long as these errors are bounded by constants in the local norm\footnote{See condition 0-3 in Lemma~\ref{lem:invariant_newton:intro} for details.}, the potential function stays bounded, which guarantees the closeness to central path. 
Therefore this analysis is currently the most robust possible. The main component of our robust analysis is the following one step error control.

\begin{lemma}[One step error control of the robust framework, informal version of Lemma~\ref{lem:invariant_newton}]\label{lem:invariant_newton:intro}
Let the potential function of IPM defined by \begin{align*}
\Psi(z,y,\eta):= \| \g(y,\eta)\|_{(\nabla^2 \phi(z))^{-1}}.
\end{align*}
Given any parameters $ \alpha_S \in [1, 1+10^{-4}]$, $c_H \in [10^{-1},1]$, $ \epsilon_g , \epsilon_{\delta} \in [0, 10^{-4}]$, and $\epsilon_N  \in (0, 10^{-1})$, $\eta >0$. Suppose that there is
\begin{itemize}
    \item {\bf Condition 0.} a feasible dual solution $y \in \R^m$ satisfies $ \Phi(y , y , \eta) \leq \epsilon_N$, 
    %\item {\bf Condition 1.} a symmetric matrix $\wt{S} \in \mathbb{S}^{n \times n}_{>0}$ satisfies $\alpha_S^{-1} \cdot S(y) \preceq \wt{S} \preceq \alpha_S \cdot S(y)$,
    \item {\bf Condition 1.} a symmetric matrix $\wt{H} \in \mathbb{S}^{n \times n}_{>0}$ satisfies $ c_H \cdot \nabla^2 \phi(y)  \preceq \wt{H} \preceq  \nabla^2 \phi(y)$,
    %\item {\bf Condition 2.} a symmetric matrix $\wt{H} \in \mathbb{S}^{n \times n}_{>0}$ satisfies $ \alpha_H^{-1} \cdot H( \wt{S} )  \preceq \wt{H} \preceq \alpha_H \cdot H( \wt{S} )$,
    \item {\bf Condition 2.} a vector $\wt{g} \in \R^m$ satisfies $\|\wt{g} - \g(y,\eta^{\new} )\|_{(\nabla^2 \phi(y))^{-1}} \leq \epsilon_g \cdot \|\g(y,\eta^{\new})\|_{(\nabla^2 \phi(y))^{-1}}$,
    \item {\bf Condition 3.} a vector $\wt{\delta}_y \in \R^m$ satisfies $\|\wt{\delta}_y - \wt{H}^{-1} \wt{g}\|_{\nabla^2 \phi(y)} \leq \epsilon_{\delta} \cdot \| \wt{H}^{-1} \wt{g} \|_{\nabla^2 \phi(y)}$.
\end{itemize}
Then $\eta^{\new} = \eta (1 + \frac{\epsilon_N}{20 \sqrt{\theta}} ) $ and $y^{\new} = y - \wt{\delta}_y$ satisfy
\begin{align*} 
\Psi(y^{\new}, y^{\new}, \eta^{\new}) \leq \epsilon_N. 
\end{align*}
\end{lemma}

This result suggests that as long as we find an initial dual variable $y$ in the proximity of central path, i.e. $\Phi(y,y,\eta) \leq \epsilon_N$, Lemma~\ref{lem:invariant_newton:intro} will guarantee that the invariant $\Phi(y,y,\eta) \leq \epsilon_N$ holds throughout Algorithm~\ref{alg:robust_ipm:intro}, even when there exist errors in the slack matrices, Hessian, gradient and Newton steps. As shown in Section~\ref{sec:correct_previous}, the duality gap is upper bounded by $\theta \cdot \Phi(y,y,\eta) / \eta$. In at most $O(\sqrt{\theta} \cdot \log(\theta/\epsilon))$ iterations, $\eta$ will become greater than $\theta \cdot \Phi(y,y,\eta) / \epsilon$. Therefore Algorithm~\ref{alg:robust_ipm:intro} finds $\epsilon$-optimal solution within $O(\sqrt{\theta} \cdot \log(\theta/\epsilon))$ iterations. 

We note that \cite{a00} and \cite{r01} only consider Condition 0 and requires the $c_H = 1, \epsilon_g = \epsilon_\delta = 0$ in Condition 1, 2, and 3. \cite{jklps20} considered Condition 0 and Condition 1 in Lemma~\ref{lem:invariant_newton:intro} and requires the $\epsilon_g = \epsilon_\delta = 0$ in Condition 2 and 3. Moreover, the Condition 1 in \cite{jklps20} requires $c_H$ to be very close to 1, and we relax this condition to support any constant in $[10^{-1},1]$. In addition, our framework also relaxes the computation of gradient and Newton direction to allow some approximations, which makes it possible to apply more algorithmic techniques in the interior-point method. More details are provided in Section~\ref{sec:our_roubst_ns}.

\section{Our Techniques}\label{sec:our_technique}
In this section, we introduce our main techniques, and provide a self-contained proof sketch of our main result Theorem~\ref{thm:tall_SDP_informal}. We tackle the two bottlenecks of $m^{\omega}$ cost per iteration in \cite{jklps20} by proposing two different techniques:

{\bf Bottleneck 1:} Instead of inverting the Hessian matrix from scratch in each iteration, we make use of the already-computed Hessian inverse of the previous iteration. We prove that using low-rank updates, the change to the inverse of Hessian matrices (computed using Kronecker product) is low-rank, and thus we can use Woodbury identity to efficiently update the Hessian inverse.  In Section~\ref{sec:tech_slack_low_rank} we introduce the low-rank update to the Hessian, and in Section~\ref{sec:tech_inverse_hessian_woodbury} we describe how to compute the Hessian inverse efficiently using Woodbury identity and fast matrix rectangular multiplication.

 {\bf Bottleneck 2:} We propose a better amortization scheme for PSD matrices that improves upon the previous $m^{\omega}$ amortized cost.  We give a proof sketch of our amortized analysis in Section~\ref{sec:tech_amortize}.

%\paragraph{Preliminary.} 
%We use $\otimes$ to denote matrix Kronecker product: for matrices $A \in \R^{m \times n}$ and $B \in \R^{p \times q}$, $A \otimes B \in \R^{pm \times qn}$ is defined to be $(A\otimes B)_{p(i-1) + s, q(j-1) + t} = A_{i,j} \cdot B_{s,t}$ for any $i \in [m]$, $j \in [n]$, $s \in [p]$, $[t] \in [q]$.
\iffalse
i.e., 
\begin{align*}
     A \otimes B  = 
     \begin{bmatrix} 
     A_{1,1} \cdot B & A_{1,2} \cdot B & \dots & A_{1,n} \cdot B\\ 
     A_{2,1} \cdot B & A_{2,2} \cdot B & \cdots & A_{2,n} \cdot B \\
     \vdots & \vdots & \ddots &\vdots \\
     A_{m,1} \cdot B & A_{m,2} \cdot B & \dots & A_{m,n} \cdot B
     \end{bmatrix} 
     \in \R^{pm \times qn}.
 \end{align*}
 \fi
 
 \iffalse
To overcome the bottleneck of inverting the Hessian matrix, we leverage the Woodbury's identity that shows a low-rank update of a matrix leads to a low-rank update of its inverse:
\begin{fact}[Woodbury matrix identity, \cite{w49,w50}]\label{fac:woodbury_simple}
For matrices $M\in \R^{n\times n}$, $U\in \R^{n\times k}$, $C\in \R^{k\times k}$, $V\in \R^{k\times n}$,
    $
        ( M + U C V )^{-1} = M^{-1} - M^{-1} U (C^{-1} + V M^{-1} U )^{-1} V M^{-1} .
    $
\end{fact}
\fi

%\paragraph{Algorithm.} A simplified version of our algorithm is given in Algorithm~\ref{alg:informal}.

 \begin{algorithm}[!ht]\caption{Informal version of Alg.~\ref{alg:ss_ipm_woodbury}. An implementation of \textsc{GeneralRobustSDP}}\label{alg:informal}
  %\small
 \begin{algorithmic}[1]
 \Procedure{\textsc{SolveSDP}}{ $\mathsf{A} \in \R^{m \times n^2}$, $b \in \R^m$, $C \in \R^{n \times n}$} \For {$t = 1 \to T$} \Comment{$T = \wt{O}(\sqrt{n})$}
 	\State $\eta^{\new} \leftarrow \eta \cdot ( 1 + {1}/{\sqrt{n}})$
 	\State $g_{\eta^{\new}}(y)_j \leftarrow  \eta^{\new} \cdot b_j -   \tr [ S^{-1} \cdot A_j ] $,~ $\forall j \in m$ \Comment{Gradient computation}
 	\State $\delta_y \leftarrow - \wt{H}^{-1} \cdot g_{\eta^{\new}}(y)$ \Comment{Compute Newton step} 	
 	\State $y^{\new} \leftarrow y + \delta_y$ \Comment{Update dual variables}
 	\State {$S^{\new} \leftarrow \sum_{i \in [m]} (y^{\new})_i A_i - C$} \Comment{Compute slack matrix}
 	\State Compute $V_1, V_2 \in \R^{n \times r_t}$ such that $\wt{S}^{\new} = \wt{S} + V_1 \cdot V_2^{\top}$    \Comment{Step 1 of Sec.~\ref{sec:tech_slack_low_rank}} 	\State Compute $V_3, V_4 \in \R^{n \times r_t}$ such that $(\wt{S}^{\new})^{-1} = (\wt{S})^{-1} + V_3 \cdot V_4^{\top}$   \Comment{ Step 2 of Sec.~\ref{sec:tech_slack_low_rank}}
 	\State Compute $\mathsf{A} Y_1, \mathsf{A} Y_2 \in \R^{m \times n r_t}$ such that $\wt{H}^{\new} = \wt{H} + (\mathsf{A} Y_1) \cdot (\mathsf{A} Y_2)^{\top}$   \Comment{ Step 3 of Sec.~\ref{sec:tech_slack_low_rank}}
 	\State $(\wt{H}^{\new})^{-1} \leftarrow \wt{H}^{-1} + \text{low-rank update}$ %, use fast rectangular matrix multiplication
  \Comment{Sec.~\ref{sec:tech_inverse_hessian_woodbury}} %Efficiently update Hessian inverse,
 	\State $y \leftarrow y^{\new}$,  $S \leftarrow S^{\new}$, $\wt{S} \leftarrow \wt{S}^{\new}, \wt{H}^{-1}\leftarrow (\wt{H}^{\new})^{-1}$ \Comment{Update variables}
 \EndFor
 \EndProcedure
 \end{algorithmic}
 \end{algorithm}

\subsection{Low rank update of Hessian}\label{sec:tech_slack_low_rank}
Low-rank approximation of Kronecker product itself is an interesting problem and has been studied in \cite{swz19}. In this section, we describe how the low-rank update of the slack matrix leads to a low-rank update of the Hessian matrix that involves Kronecker product.

The Hessian matrix is defined as %(equivalent to Eq.~\eqref{eq:hessian_prev})
$
    H = \mathsf{A} \cdot ({S}^{-1} \otimes {S}^{-1} ) \cdot \mathsf{A}^{\top}.
$ 
%where $\mathsf{A} \in \R^{m \times n^2}$ stacks the vectorization of $A_1, \cdots, A_m \in \R^{n \times n}$. 
We take the following three steps to construct the low-rank update of $H$.

\paragraph{Step 1: low-rank update of the slack matrix.}
%follow \cite{jklps20} to
We use an approximate slack matrix that yields a low-rank update. In the $t$-th iteration of Algorithm~\ref{alg:informal}, we use $\wt{S}$ to denote the current approximate slack matrix, and ${S}^\new$ to denote the new exact slack matrix. We will use $\wt{S}$ and ${S}^\new$ to find the new approximate slack matrix $\wt{S}^{\new}$.

Define $Z = (S^{\new})^{-1/2} \wt{S} (S^\new)^{-1/2} -I$ which captures the changes of the slack matrix. We compute the spectral decomposition: $Z = U \cdot \diag(\lambda) \cdot U^\top$. 
%It is observed in \cite{jklps20} 
We show that 
\begin{align*}
    \sum_{i = 1}^n \lambda_i^2 = \|S^{-1/2}S^\new S^{-1/2} -I\|_F = O(1),
\end{align*}
which implies that only a few eigenvalues of $Z$ are significant, say e.g. $\lambda_1,\dots,\lambda_{r_t}$. We only keep these eigenvalues and set the rest to be zero. In this way we get a low-rank approximation of $Z$: $\wt{Z} = U \cdot \diag(\wt{\lambda}) \cdot U^\top$ where $\wt{\lambda} = [\lambda_1,\cdots, \lambda_{r_t},0,\dots,0]^{\top}$. Now we can use $\wt{Z}$ to update the approximate slack matrix by a low-rank matrix:
$$
    \wt{S}^\new = \wt{S} + (S^{\new})^{1/2} \cdot \wt{Z} \cdot (S^{\new})^{1/2} = \wt{S} + V_1 \cdot V_2^{\top},
$$
where %$V_1 := (S^{\new})^{1/2} \cdot U \cdot \diag(\wt{\lambda})$ and $V_2 := (S^{\new})^{1/2} \cdot U$ 
$V_1$ and $V_2$ both have size $n \times r_t$.
Since $\wt{Z}$ is a good approximation of $Z$, $\wt{S}^\new$ is a PSD approximation of ${S}^\new$, which guarantees that $y$ still lies in the proximity of the central path.

\paragraph{Step 2: low-rank update of inverse of slack.}
Using Woodbury identity, we can show that 
$$
    (\wt{S}^{\new})^{-1} = (\wt{S} + V_1 \cdot V_2^{\top})^{-1} = \wt{S}^{-1} + V_3 V_4^\top,
$$
where $V_3 = - \wt{S}^{-1} V_1 ( I + V_2^\top \wt{S}^{-1} V_1 )^{-1} $ and $V_4 = \wt{S}^{-1} V_2$ both have size $n \times r_t$. Thus, this means $(\wt{S}^{\new})^{-1} - \wt{S}^{-1}$ has a rank $r_t$ decomposition.

\paragraph{Step 3: low-rank update of Hessian.}
Using the linearity and the mixed product property (Part~2 of Fact~\ref{fac:kron_basic}) of Kronecker product, we can find a low-rank update to $(\wt{S}^{\new})^{-1} \otimes (\wt{S}^{\new})^{-1}$. More precisely, we can rewrite $(\wt{S}^{\new})^{-1} \otimes (\wt{S}^{\new})^{-1} $ as follows:
$$  
(\wt{S}^{\new})^{-1} \otimes (\wt{S}^{\new})^{-1} 
    = (\wt{S}^{-1} + V_3 V_4^{\top}) \otimes (\wt{S}^{-1} + V_3 V_4^{\top}) 
    = ~ \wt{S}^{-1} \otimes \wt{S}^{-1} + {\cal S}_{\mathrm{diff}}.
$$
The term ${\cal S}_{\mathrm{diff}}$ is the difference that we want to compute, we can show
\begin{align*}
    {\cal S}_{\mathrm{diff}} 
    = & ~ \wt{S}^{-1} \otimes (V_3 V_4^{\top}) + (V_3 V_4^{\top}) \otimes \wt{S}^{-1} + (V_3 V_4^{\top}) \otimes (V_3 V_4^{\top}) \\
    = & ~ (\wt{S}^{-1/2} \otimes V_3) \cdot (\wt{S}^{-1/2} \otimes V_4^{\top}) + (V_3 \otimes \wt{S}^{-1/2}) \cdot (V_4^{\top} \otimes \wt{S}^{-1/2}) + (V_3 \otimes V_3) \cdot (V_4^{\top} \otimes V_4^{\top}) \\
    = & ~ Y_1 \cdot Y_2^{\top}
\end{align*}
where $Y_1$ and $Y_2$ both have size $n^2 \times n r_t$. In this way we get a low-rank update to the Hessian:
\begin{align*}
    \wt{H}^{\new} = \mathsf{A} \cdot ((\wt{S}^{\new})^{-1} \otimes (\wt{S}^{\new})^{-1}) \cdot \mathsf{A}^{\top} = \wt{H} + (\mathsf{A} Y_1) \cdot (\mathsf{A} Y_2)^{\top}.
\end{align*}

\subsection{Computing Hessian inverse efficiently}\label{sec:tech_inverse_hessian_woodbury}

In this section we show how to compute the Hessian inverse efficiently.

Using Woodbury identity again, we have a low rank update to $\wt{H}^{-1}$: 
\begin{align*}
    (\wt{H}^{\new})^{-1} 
    = ~ \big( \wt{H} + (\mathsf{A} Y_1) \cdot (\mathsf{A} Y_2)^{\top} \big)^{-1} 
    = ~ \wt{H}^{-1} - \wt{H}^{-1} \cdot \mathsf{A} Y_1 \cdot (I + Y_2^{\top} \mathsf{A}^{\top} \cdot \mathsf{A} Y_1)^{-1} \cdot Y_2^{\top} \mathsf{A}^{\top} \cdot \wt{H}^{-1} 
\end{align*}
The second term in the above equation has rank $nr$. Thus $(\wt{H}^{\new})^{-1} - \wt{H}^{-1}$ has a rank $n r$ decomposition. 
%To compute $(\wt{H}^{\new})^{-1}$ we stack the $m$ vectorizations $\vect[A_1], \vect[A_2], \cdots, \vect[A_m] \in \R^{n^2}$ to form $\mathsf{A} \in \R^{m \times n^2}$. 
To compute $(\wt{H}^{\new})^{-1}$ in each iteration, we first compute $\mathsf{A} Y_1, \mathsf{A} Y_2 \in \R^{m \times n r_t}$ and multiply it with $\wt{H}^{-1} \in \R^{m \times m}$ to get $\wt{H}^{-1} \cdot \mathsf{A} Y_1, \wt{H}^{-1} \cdot \mathsf{A} Y_2 \in \R^{m \times n r_t}$. Then we compute $I + (Y_2^{\top} \mathsf{A}^{\top}) \cdot (\mathsf{A} Y_1) \in \R^{n r_t\times n r_t}$ and find its inverse $(I + Y_2^{\top} \mathsf{A}^{\top} \cdot \mathsf{A} Y_1)^{-1} \in \R^{n r_t\times n r_t}$. Finally, we multiply $\wt{H}^{-1} \cdot \mathsf{A} Y_1, \wt{H}^{-1} \cdot \mathsf{A} Y_2 \in \R^{m \times n r_t}$ and $(I + Y_2^{\top} \mathsf{A}^{\top} \cdot \mathsf{A} Y_1)^{-1} \in \R^{n r_t\times n r_t}$ together to obtain $(\wt{H})^{-1} \mathsf{A} Y_1 \cdot (I + Y_2^{\top} \mathsf{A}^{\top} \cdot \mathsf{A} Y_1)^{-1} \cdot Y_2^{\top} \mathsf{A}^{\top} (\wt{H})^{-1} \in \R^{m \times m}$, as desired.
Using fast matrix multiplication in each aforementioned step, the total computation cost is bounded by 
\begin{align}\label{eq:tech_cost_per_iter}
O(\Tmat(m, n^2, n r_t) + \Tmat(m, m, n r_t) + (n r_t)^{\omega}).
\end{align}

\subsection{General amortization method}\label{sec:tech_amortize}
As mentioned in the previous sections, our algorithm relies on the maintenance of the slack matrix and the inverse of the Hessian matrix via low-rank updates. In each iteration, the time to update $\wt{S}$ and $\wt{H}$ to $\wt{S}_{\new}$ and $\wt{H}_\new$ is proportional to the magnitude of low-rank change in $\wt{S}$, namely $r_t = \mathrm{rank}(\wt{S}_\new - \wt{S})$. 
To deal with $r_t$, we propose a general amortization method which extends the analysis of several previous work \cite{cls19,lsz19,jklps20}. %,sy21
%\footnote{
%Theorem~6.1 in \cite{jklps20} is the special case of Theorem~\ref{thm:amortize_informal} where $g_i = 1 / \sqrt{i}$
%Our result doesn't apply the multiple level amortization technique in \cite{jswz21}. It might be possible to get some extra improvements once that technique gets applied.
%} 
We first prove a tool to characterize intrinsic properties of the low-rank updates, which may be of independent interest.

\begin{theorem}[Informal version of Theorem~\ref{thm:general_amortize_guarantee}]\label{thm:amortize_informal}
Given a sequence of approximate slack matrices $\wt{S}^{(1)}, \wt{S}^{(2)}, \dots, \wt{S}^{(T)} \in \mathbb{R}^{n \times n}$ generated by Algorithm~\ref{alg:ss_ipm_woodbury}, let $r_t = \mathrm{rank}(\wt{S}^{(t+1)} - \wt{S}^{(t)})$ denotes the rank of update on $\wt{S}^{(t)}$. Then for any non-increasing vector $g \in \R_+^n$, we have 
$$
    \sum_{t=1}^T r_t \cdot g_{r_t} \leq \wt{O}(T \cdot \|g\|_2).
$$
\end{theorem}
%\begin{remark}
%Theorem~6.1 in \cite{jklps20} is the special case of Theorem~\ref{thm:amortize_informal} where $g_i = 1 / \sqrt{i}$.
%\end{remark}
Next, we show a proof sketch of Theorem~\ref{thm:amortize_informal}.
\begin{proof}%[Proof sketch]
For any matrix $Z$, let $|\lambda(Z)|_{[i]}$ denotes its $i$-th largest absolute eigenvalue. We use the following potential function
$
    \Phi_g(Z) := \sum_{i=1}^n g_i \cdot |\lambda(Z)|_{[i]}
$. Further, for convenient, we define $\Phi_g(S_1,S_2) := \Phi_g( S_1^{-1/2} S_2 S_1^{-1/2} - I )$. 
Our proof consists of the following two parts (Lemma~\ref{lem:S_move} and Lemma~\ref{lem:wt_S_move}):
\begin{itemize}
    \item {The change of the exact slack matrix increases the potential by a small amount}, specifically  $\Phi_g( S^{\new} , \wt{S} ) - \Phi_g( S, \wt{S} ) \leq \| g \|_2$.
    %\begin{align*}
    %\Phi_g((S^{\new})^{-1/2} \cdot \wt{S} \cdot (S^{\new})^{-1/2} - I) - \Phi_g(S^{-1/2} \cdot \wt{S} \cdot S^{-1/2} - I) \leq \|g\|_2.
    %\end{align*}
    \item {The change of the approximate slack matrix decreases the potential proportionally to the update rank}, specifically
    $\Phi_g( S^{\new} , \wt{S}^{\new} ) - \Phi_g( S^{\new} , \wt{S} ) \leq - r_t \cdot g_{r_t}$.
    %\begin{align*}
    %\Phi_g((S^{\new})^{-1/2} \cdot \wt{S}^{\new} \cdot (S^{\new})^{-1/2} - I) - \Phi_g((S^{\new})^{-1/2} \cdot \wt{S} \cdot (S^{\new})^{-1/2} - I) \leq - r_t \cdot g_{r_t}.
    %\end{align*}
\end{itemize}
In each iteration, the change of potential is composed of the changes of the exact and the approximate slack matrices: 
\begin{align*}
    \Phi_g( S^{\new} , \wt{S}^{\new} ) - \Phi_g( S , \wt{S} ) =  \Phi_g( S^{\new} , \wt{S} ) -  \Phi_g( S , \wt{S} ) + \Phi_g( S^{\new} , \wt{S}^{\new} ) - \Phi_g( S^{\new} , \wt{S} ).
\end{align*}

Note that $\Phi_g(S, \wt{S}) = 0$ holds in the beginning of our algorithm and $\Phi_g(S, \wt{S} ) \geq 0$ holds throughout the algorithm, combining the observations above we have $  T\cdot \|g\|_2 -\sum_{t=1}^T r_t \cdot g_{r_t} \geq 0$ as desired.
\end{proof}

\paragraph{Amortized analysis.}

Next we show how to use Theorem~\ref{thm:amortize_informal} to prove that our algorithm has an amortized cost of $m^{\omega - 1/4} + m^2$ cost per iteration when $m = \Omega(n^2)$. Note that in this case there are $\sqrt{n} = m^{1/4}$ iterations.

When $m = \Omega(n^2)$, the dominating term in our cost per iteration (see Eq.~\eqref{eq:tech_cost_per_iter}) is $\Tmat(m, m, n r_t)$. We use fast rectangular matrix multiplication to upper bound this term by
$$
    \Tmat(m, m, n r_t) \leq m^2 + m^{2 - \frac{\alpha(\omega-2)}{1-\alpha}} \cdot n^{\frac{\omega-2}{1-\alpha}} \cdot r_t^{\frac{\omega-2}{1-\alpha}}.
$$ 
We define a non-increasing sequence $g \in \R^n$ as $g_i = i^{\frac{\omega-2}{1-\alpha} - 1}$. This $g$ is tailored for the above equation, and its $\ell_2$ norm is bounded by $\|g\|_2 \leq n^{\frac{(\omega-2)}{1-\alpha}-1/2}$. Then using Theorem~\ref{thm:amortize_informal} we have 
\begin{align*}
    \sum_{t=1}^T r_t^{\frac{\omega-2}{1-\alpha}} =\sum_{t=1}^T r_t\cdot  r_t^{\frac{\omega-2}{1-\alpha}-1}= \sum_{t=1}^T r_t \cdot g_{r_t} \leq T \cdot n^{\frac{(\omega-2)}{1-\alpha}-1/2}. 
\end{align*}
Combining this and the previous equation, and since we assume $m = \Omega(n^2)$, we have 
\begin{align*}
    \sum_{t=1}^T \Tmat(m, m, n r_t) \leq T \cdot (m^2 + m^{2 - \frac{\alpha(\omega-2)}{1-\alpha}} \cdot n^{\frac{2(\omega-2)}{1-\alpha}-1/2}) = T \cdot (m^2 + m^{\omega - 1/4}).
\end{align*}
Since $T = \wt{O}(m^{1/4})$, we proved the desired computational complexity in Theorem~\ref{thm:tall_SDP_informal}.

%\vspace{-1mm}
\section{Solving SDP With Hybrid Barrier}\label{sec:intro_hybrid}
%\vspace{-1mm}

Volumetric barrier was first proposed by Vaidya \cite{v89} for the polyhedral, and was generalized to the spectrahedra $\{y\in \R^m: y_1A_1+\cdots +y_mA_m\succeq 0\}$ by Nesterov and Nemirovski \cite{nn94}. They showed that the volumetric barrier $\phi_{\mathrm{vol}}$ can make the interior point method converge in $\sqrt{m}n^{1/4}$ iterations, while the log barrier $\phi_{\log}$ need $\sqrt{n}$ iterations. By combining the volumetric barrier and the log barrier, they also showed that the hybrid barrier achieves $(mn)^{1/4}$ iterations. Anstreicher \cite{a00} gave a much simplified proof of this result.

We show that the hybrid barrier also fits into our robust IPM framework. And we can apply our newly developed low-rank update and amortization techniques in the log barrier case to efficiently implement the SDP solver based on hybrid barrier. The informal version of our result is stated in below. 
\begin{theorem}[Informal version of Theorem~\ref{thm:vol_formal}]\label{thm:vol_intro}
There is an SDP algorithm based on hybrid barrier which takes $(mn)^{1/4}\log(1/\epsilon)$ iterations with cost-per-iteration $O^\ast\left(m^{2} n^{\omega} +  m^4\right)$.

In particular, our algorithm improves \cite{a00} in nearly all parameter regimes.
For example, if $m =n^2$, our new algorithm takes $n^{8.75}$ time while \cite{a00} takes $n^{10.75}$ time. If $m = n$, our new algorithm takes $n^{\omega+2.5}$ time, while \cite{a00} takes $n^{6.5}$ time. 
\end{theorem}

The hybrid barrier function is as follows:
\begin{align*}
    \phi(y):=225\sqrt{\frac{n}{m}} \cdot \left(\phi_{\mathrm{vol}}(y) + \frac{m-1}{n-1}\cdot \phi_{\log}(y)\right),
\end{align*}
where $\phi_{\mathrm{vol}}(y) = \frac{1}{2}\log \det(\nabla^2\phi_{\log}(y))$. According to our general IPM framework (Algorithm~\ref{alg:robust_ipm:intro}), we need to efficiently compute the gradient and Hessian of $\phi(y)$. Recall from \cite{a00} that the gradient of the volumetric barrier is:
\begin{align*}
    (\nabla \phi_{\mathrm{vol}}(y))_i = -\tr[H(S)^{-1} \cdot \A  (S^{-1}A_i S^{-1}\otimes S^{-1}) \mathsf{A}^\top]~~~\forall i\in [m].
\end{align*}
And the Hessian can be written as $\nabla^2\phi_{\mathrm{vol}}(y) = 2Q(S) + R(S) - 2T(S)$, where for any $i,j\in [m]$,
\begin{align}\label{eq:hessian_vol_intro}
    Q(S)_{i,j} = &~ \tr[H(S)^{-1} \A (S^{-1}A_iS^{-1}A_jS^{-1} \otimes_S S^{-1}) \A^\top],\notag\\
    R(S)_{i,j} = &~ \tr[H(S)^{-1} \A (S^{-1}A_iS^{-1} \otimes_S S^{-1}A_jS^{-1}) \A^\top],\\
    T(S)_{i,j}=&~ \tr[H(S)^{-1} \A (S^{-1}A_iS^{-1} \otimes_S S^{-1}) \A^\top H(S)^{-1}\A (S^{-1}A_jS^{-1} \otimes_S S^{-1}) \A^\top ].\notag
\end{align}
Here, $\otimes_S$ is the symmetric Kronecker product\footnote{$X\otimes_S Y := \frac{1}{2}(X\otimes Y + Y\otimes X)$.}. 

A straight-forward implementation of the hybrid barrier-based SDP algorithm can first compute the matrices $S^{-1}A_i$ and $S^{-1}A_iS^{-1}A_j$ for all $i \in \{1,2,\cdots,m\}$ for all $j\in \{1,2,\cdots,m\}$ in time $O(m^2n^\omega)$. The gradient $\nabla \phi(y)$ and the Hessian of $\phi_{\log}(y)$ can be computed by taking traces of these matrices. To compute $\nabla \phi_{\mathrm{vol}}(y), Q(S), R(S), T(S)$, we observe that each entry of these matrices can be written as the inner-product between $H(S)^{-1}$ and some matrices formed in terms of $\tr[S^{-1}A_iS^{-1}A_jS^{-1}A_k]$ and $\tr[S^{-1}A_iS^{-1}A_jS^{-1}A_kS^{-1}A_l]$ for $i,j,k,l\in [m]$. Hence, we can spend $O(m^4n^2)$-time computing these traces and then get $\nabla \phi_{\mathrm{vol}}(y), Q(S), R(S), T(S)$ in $O(m^{\omega+2})$-time. After obtaining the gradient and Hessian of the hybrid barrier function, we finish the implementation of IPM SDP solver by computing the Newton direction $\delta_y = - (\nabla^2 \phi(y))^{-1} (\eta b - \nabla \phi(y))$. (More details are given in Section~\ref{sec:straightforward_vol}).

To speedup the straight forward implementation, we observe two bottleneck steps in each iteration:
\begin{enumerate}
    \item Computing the traces $\tr[S^{-1}A_iS^{-1}A_jS^{-1}A_kS^{-1}A_l]$ for $i,j,k,l\in [m]$.
    \item Computing the matrices $Q(S), R(S), T(S)$.
\end{enumerate}

To handle the first issue, we use the low-rank update and amortization techniques introduced in the previous section to approximate the change of the slack matrix $S$ by a low-rank matrix. One challenge for the volumetric barrier is that its Hessian (Eq.~\eqref{eq:hessian_vol_intro}) is much more complicated than the log barrier's Hessian $H(S)$. For $H(S)$, if we replace $S$ with its approximation $\wt{S}$, then $H(\wt{S})$ will be a PSD approximation of $H(S)$. However, this may not hold for the volumetric barrier's Hessian if we simply replace all the $S$ in $\nabla^2\phi(y)$ by its approximation $\wt{S}$. We can resolve this challenge by carefully choosing the approximation place: if we approximate the second $S$ in the trace, i.e., $\tr[S^{-1}A_i\wt{S}^{-1}A_jS^{-1}A_kS^{-1}A_l]$, then the resulting matrix will be a PSD approximation of $\nabla^2\phi(y)$. In other words, the Condition 1 in our robust IPM framework (Lemma~\ref{lem:invariant_newton:intro}) is satisfied. Notice that in each iteration, we only need to maintain the change of  $\tr[S^{-1}A_i\wt{S}^{-1}A_jS^{-1}A_kS^{-1}A_l]$, which by the low-rank guarantee, can be written as
\begin{align*}
    \tr[A_lS^{-1}A_i\cdot V_3V_4^\top\cdot  A_jS^{-1}A_kS^{-1}],
\end{align*}
where $V_3,V_4\in \R^{n\times r_t}$. Then, we can first compute the matrices 
\begin{align*}
    \big\{A_l S^{-1}A_iV_3\in \R^{n\times r_t}\big\}_{i,l\in [m]}~~\text{and}~~\big\{V_4^\top A_jS^{-1}A_kS^{-1}\in \R^{r_t\times n}\big\}_{j,k\in [m]}.
\end{align*}
It takes $m^2\cdot {\cal T}_{\mathrm{mat}}(n,n,r_t)$-time. And we can compute all the traces $\tr[S^{-1}A_i\wt{S}^{-1}A_jS^{-1}A_kS^{-1}A_l]$ simultaneously in ${\cal T}_{\mathrm{mat}}(m^2, nr_t, m^2)$ by batching them together and using fast matrix multiplication on a $m^2$-by-$nr_t$ matrix and a $nr_t$-by-$m^2$ matrix. A similar amortized analysis in the log barrier case can also be applied here to get the amortized cost-per-iteration for the low-rank update. One difference is that the potential function $\Phi_g(Z)$ (defined in Section~\ref{sec:tech_amortize}) changes more drastically in the hybrid barrier case. And we can only get $\sum_{t=1}^T r_t \cdot g_{r_t} \leq O(T \cdot (n/m)^{1/4}\cdot \|g\|_2 \cdot \log n)$. 

For the second issue, we note that computing the $T(S)$ matrix is the most time-consuming step, which need $m^{\omega+2}$-time. In \cite{a00}, it is proved that $\frac{1}{3} Q(S) \preceq \nabla^2\phi_{\mathrm{vol}} (y) \preceq Q(S)$.
With this PSD approximation, our robust IPM framework enables us to use $Q(S)$ as a ``proxy Hessian'' of the volumetric barrier. That is, in each iteration, we only compute $Q(S)$ and ignore $R(S)$ and $T(S)$. And computing $Q(S)$ only takes $O(m^4)$-time, which improves the $m^{\omega+2}$ term in the straight forward implementation.

Combining them together, we obtain the running time in Theorem~\ref{thm:vol_intro}. More details are provided in Section~\ref{sec:hybrid_barrier}.

%\vspace{-1mm}
\paragraph{Lee-Sidford barrier for SDP?}
In LP, the hybrid barrier was improved by Lee and Sidford \cite{ls19} to achieve $O^*(\sqrt{\min\{m,n\}})$ iterations. For SDP, we hope to design a barrier function with $O^*(\sqrt{m})$ iterations. However, the Lee-Sidford barrier function does not have a direct correspondence in SDP due to the following reasons. First, \cite{ls19} defined the barrier function in the dual space of LP which is a polyhedron, while for SDP, the dual space is a spectrahedron. Thus, the geometric intuition of the Lee-Sidford barrier (John's ellipsoid) may not be helpful to design the corresponding barrier for SDP. Second, efficient implementation of Lee-Sidford barrier involves a primal-dual central path method \cite{blss20}. However, the cost of following primal-dual central path in SDP is prohibitive since this involves solving Lyapunov equations in $\R^{n \times n}$. Third, the Lewis weights play an important role in the Lee-Sidford barrier. Notice that in LP, the volumetric barrier can be considered as reweighing the constraints in the log barrier based on the leverage score, and the Lee-Sidford barrier uses Lewis weights for reweighing to improve the volumetric barrier. However, in SDP, we have observed that the leverage score vector becomes the leverage score matrix. Thus, we may need some matrix version of Lewis weights to define the Lee-Sidford barrier for SDP. Section~\ref{sec:leverage_matrix} studies several properties of the leverage score matrix and give an algorithm to efficiently maintain this matrix in each iteration of the IPM, which might be the first step towards improving the SDP hybrid barrier.

%\vspace{-2mm}
\section{Related Work}\label{sec:related_work}
%\vspace{-2mm}

\paragraph{Other SDP solvers.}
The interior point method is a second-order algorithm. Second-order algorithms usually have logarithmic dependence on the error parameter $1/\epsilon$. First-order algorithms do not need to use second-order information, but they usually have polynomial dependence on $1/\epsilon$. There is a long list of work focusing on first-order algorithms \cite{ak07,jy11,alo16,gh16,al17,cdst19,lp19,ytfuc19,jllpt20}. Solving SDPs has also attracted attention in the parallel setting \cite{jy11,jy12,alo16,jllpt20}.

%\vspace{-3mm}
\paragraph{Cutting plane method.}

Cutting plane method is a class of optimization algorithms that iteratively queries a separation oracle to cut the feasible set that contains the optimal solution. There has been a long line of work to obtain fast cutting plane methods \cite{s77,yn76,k80,kte88,nn89,v89,av95,bv02,lsw15,jlsw20}.

%\vspace{-2mm}
\paragraph{Low-rank approximation}
Low-rank approximation is a well-studied topic in numerical linear algebra \cite{s06,cw13,bwz16,swz17,swz19}. Many different settings of that problem have been studied. In this paper, we are dealing with Kronecker product type low rank approximation.

%\vspace{-3mm}
\paragraph{Applications of SDP.}

As described by \cite{jklps20}, $m = \Omega(n^2)$ is an essential case of using SDP to solve many practical combinatorial optimization problems. Here we provide a list of examples, e.g., the sparsest cut \cite{arv09}, the $c$-balanced graph separation problem \cite{fhl08} and the minimum uncut \cite{acmm05} can be solved by SDP with $m=\Omega(n^3)$. The optimal experiment design \cite{vbw98}, Haplotype frequencies estimation \cite{hh06} and embedding of finite metric spaces into $\ell^2$ \cite{llr95} need to solve SDPs with $m=\Omega(n^2)$.

\ifdefined\isarxivversion

\section*{Acknowledgments}
The authors would like to thank Haotian Jiang and Yin Tat Lee for many helpful discussions and insightful comments on manuscripts. The authors would like to thank Ainesh Baksh, Sitan Chen and Jerry Li for useful discussions on Sum of Squares. The authors would like to thank Elaine Shi and Kai-min Chung for useful discussions. The authors would like to thank Binghui Peng and Hengjie Zhang for useful discussions at the early stage of this project.

Baihe Huang is supported by the Elite Undergraduate Training Program of School of Mathematical Sciences at Peking University, and this work is done while interning at Princeton University and Institute for Advanced Study (advised by Zhao Song).

Shunhua Jiang is supported by NSF CAREER award CCF-1844887.
\else
\fi

\ifdefined\isarxivversion

\else
\fi

\newpage
{\hypersetup{linkcolor=black}
\tableofcontents
}
\newpage
\section{Preliminary}\label{sec:preli}

\subsection{Notations}

\paragraph{Basic matrix notations.}
For a square matrix $X$, we use $\tr[X]$ to denote the trace of $X$.

We use $\| \cdot \|_2$ and $\|\cdot\|_F$ to denote the spectral norm and Frobenious norm of a matrix. 
Let us use $\|\cdot \|_1$ to represent the Schatten-1 norm of a matrix, i.e., $\|A\|_1 = \tr[(A^*A)^{1/2}]$. 

We say a symmetric matrix $A \in \R^{n \times n}$ is positive semi-definite (PSD, denoted as $A \succeq 0$) if for any vector $x \in \R^n$, $x^{\top} A x \geq 0$. 
We say a symmetric matrix $A \in \R^{n \times n}$ is positive definite (PD, denoted as $A \succ 0$) if for any vector $x \in \R^n$, $x^{\top} A x > 0$. 

We define $\mathbb{S}^{n \times n}_{\succ 0}$ to be the set of all $n$-by-$n$ symmetric positive definite matrices. 

Let us define $\mathbb{S}^{n \times n}_{\succeq 0}$ to be the set of all $n$-by-$n$ symmetric positive semi-definite matrices.

For a matrix $A \in \R^{m \times n}$, we use $\lambda(A) \in \R^n$ to denote the eigenvalues of $A$.

For any vector $v \in \R^n$, we use $v_{[i]}$ to denote the $i$-th largest entry of $v$.

For a matrix $A \in \R^{m \times n}$, and subsets $S_1\subseteq [m], S_2\subseteq [n]$, we define $A_{S_1,S_2} \in \R^{|S_1|\times |S_2|}$ to be the submatrix of $A$ that only has rows in $S_1$ and columns in $S_2$. We also define $A_{S_1,:} \in \R^{|S_1| \times n}$ to be the submatrix of $A$ that only has rows in $S_1$, and $A_{:,S_2} \in \R^{m \times |S_2|}$ to be the submatrix of $A$ that only has columns in $S_2$.

For two symmetric matrices $A,B\in \R^{n\times n}$, we say $A\preceq B$ (or equivalently, $B\succeq A$), if $B-A$ is a PSD matrix.

\begin{fact}[Spectral norm implies Loewner order]\label{fac:spec_psd}
Let $A,B\in \R^{n\times n}$ be two symmetric PSD matrices. Then, for any $\epsilon\in (0,1)$,
\begin{align*}
    \left\|A^{-1/2} B A^{-1/2} - I\right\|_2 \leq \epsilon
\end{align*}
implies
\begin{align*}
    (1-\epsilon)A \preceq B \preceq (1+\epsilon)A.
\end{align*}
\end{fact}

\begin{fact}[Trace property of matrix Loewner order]\label{fac:trace_loewner}
Given symmetric PSD matrices $A,B \in \R^n$. Suppose $(1+\epsilon)^{-1} \cdot A \preceq \wt{A} \preceq (1+\epsilon) \cdot A$, then
\begin{align*}
    (1+\epsilon)^{-1} \cdot \tr[AB] \leq \tr[\wt{A}B] \leq (1+\epsilon) \cdot \tr[AB].
\end{align*}
\end{fact}
\begin{proof}
Consider the spectral decomposition of $B$: $B = \sum_{i=1}^n \lambda_i v_iv_i^\top$ where $\lambda_i \geq 0$. Then
\begin{align*}
    \tr[\wt{A}B] = &~ \tr[\wt{A} \cdot (\sum_{i=1}^n \lambda_i v_iv_i^\top)]\\
    =&~ \sum_{i=1}^n \lambda_i v_i^\top \wt{A} v_i\\
    \leq &~ (1+\epsilon)\cdot (\sum_{i=1}^n \lambda_i v_i^\top {A} v_i)\\
    = &~ \tr[AB].
\end{align*}
Similarly, $\tr[\wt{A}B] \geq (1+\epsilon)^{-1} \cdot \tr[AB]$.
\end{proof}

\paragraph{Matrix related operations}
For two matrices $A,B \in \R^{m \times n}$, we define the matrix inner product $\langle A, B \rangle := \tr[A^{\top} B]$.

We use $\vect[]$ to denote matrix vectorization: for a matrix $A \in \R^{m \times n}$, $\vect[A] \in \R^{mn}$ is defined to be $\vect[A]_{(j-1)\cdot n + i} = A_{i,j}$ for any $i \in [m]$ and $j \in [n]$, i.e., \begin{align*}
   \vect[A]  = \begin{bmatrix} A_{:,1} \\ \vdots\\ A_{:,n}\end{bmatrix} \in \R^{mn}.
\end{align*}

We use $\otimes$ to denote matrix Kronecker product: for matrices $A \in \R^{m \times n}$ and $B \in \R^{p \times q}$, $A \otimes B \in \R^{pm \times qn}$ is defined to be $(A\otimes B)_{p(i-1) + s, q(j-1) + t} = A_{i,j} \cdot B_{s,t}$ for any $i \in [m]$, $j \in [n]$, $s \in [p]$, $[t] \in [q]$, i.e., 
\begin{align*}
     A \otimes B  = 
     \begin{bmatrix} 
     A_{1,1} \cdot B & A_{1,2} \cdot B & \dots & A_{1,n} \cdot B\\ 
     A_{2,1} \cdot B & A_{2,2} \cdot B & \cdots & A_{2,n} \cdot B \\
     \vdots & \vdots & \ddots &\vdots \\
     A_{m,1} \cdot B & A_{m,2} \cdot B & \dots & A_{m,n} \cdot B
     \end{bmatrix} 
     \in \R^{pm \times qn}.
 \end{align*}

\begin{definition}[Stacking matrices]\label{def:stacking_matrices}
Let $A_1, A_2, \cdots, A_m \in \R^{n \times n}$ be $m$ symmetric matrices.
We use $\mathsf{A} \in \R^{m \times n^2}$ to denote the matrix that is constructed by stacking the $m$ vectorizations $\vect[A_1], \cdots, \vect[A_m] \in \R^{n^2}$ as rows of $\mathsf{A}$, i.e.,
\begin{align*}
    \mathsf{A} := \begin{bmatrix}
    \vect[A_1]^\top\\
    \vdots\\
    \vect[A_m]^\top
    \end{bmatrix}\in \R^{m\times n^2}.
\end{align*}
\end{definition}

\begin{fact}\label{fac:sigma_z_sigma_z_norm}
For any $\epsilon_1,\epsilon_2 \in (0,1/10)$.
Let $D \in \R^{ n \times n}$ be a diagonal matrix with non-negative entries and such that $\|D^2 - I \|_F \leq \epsilon_1$. Let $X \in \R^{n \times n}$ be a matrix that has bounded norm, e.g., $\|X\|_2 \leq \epsilon_2$. Then 
\begin{align*}
    \|D X D - X\|_F \leq 3 \cdot \epsilon_1 \cdot \epsilon_2.
\end{align*}
\end{fact}

\begin{proof}
Denote $D = \diag (\sigma_1,\dots,\sigma_n)$. We have
\begin{align*}
    \|D  X D - X\|_F \leq &~ \|(D - I) X (D - I) + (D - I) X + X( D - I) \|_F\\
    \leq &~ \|(D - I) X (D - I)\|_F + 2 \cdot \| (D - I) X \|_F\\
    \leq &~ 3 \cdot \| (D - I) X \|_F\\
    \leq &~ 3 \cdot \left(\tr[(D - I)^2 X^2]\right)^{1/2}\\
    \leq &~3 \cdot  \left( \epsilon_2^2 \cdot \tr[(D - I)^2]\right)^{1/2}\\
    = &~ 3 \cdot \epsilon_2 \cdot \left(\sum_{i=1}^n (\sigma_i - 1)^2 \right)^{1/2}\\
    \leq &~ 3 \cdot  \epsilon_2 \cdot \left(\sum_{i=1}^n (\sigma_i^2 - 1)^2 \right)^{1/2}\\
    \leq &~ 3 \cdot  \epsilon_2 \cdot \epsilon_1
\end{align*}
where the second step uses triangle inequality, the third step uses $-I \preceq D - I \preceq I$, the fifth step uses $\|X\|_2 \leq \epsilon_2$, the penultimate step uses $\sigma_i \geq 0$, and the last step uses $\|D^2 - I \|_F \leq \epsilon_1$.
\end{proof}

\subsection{Tools: Woodbury identity}
We state a common fact on matrix inverse update in \cite{w49,w50}.
\begin{fact}[Woodbury matrix identity]\label{fac:woodbury}
Given two integers $n$ and $k$. Let $n \geq k$.
For square matrix $A\in \R^{n\times n}$, tall matrix $B\in \R^{n\times k}$, square matrix $C\in \R^{k\times k}$, fat matrix $D\in \R^{k\times n}$,
    \begin{align*}
        ( A + B C D )^{-1} = A^{-1} - A^{-1} B (C^{-1} + D A^{-1} B )^{-1} D A^{-1} .
    \end{align*}
\end{fact}

\subsection{Tools: Properties of matrix operations}
\begin{fact}[Matrix inner product]\label{fac:matrix_inner_product}
For two matrices $A, B \in \R^{m \times n}$, we have $\langle A, B \rangle = \tr[A^{\top} B] = \vect[A]^{\top} \vect[B]$.
\end{fact}

\begin{fact}[Basic properties of Kronecker product]\label{fac:kron_basic}
The Kronecker product $\otimes$ satisfies the following properties.
\begin{enumerate}
    \item \label{part:kron_transpose} For matrices $A \in \R^{a \times n}$ and $B \in \R^{b \times m}$, we have $(A \otimes B)^{\top} = A^{\top} \otimes B^{\top} \in \R^{nm \times ab}$.
    \item \label{part:kron_mixed_product} For matrices $A \in \R^{a \times n}$, $B \in \R^{b \times m}$, $C \in \R^{n \times c}$, $D \in \R^{m \times d}$, we have $(A \otimes B) \cdot (C \otimes D) = (AC \otimes BD)\in \R^{ab\times cd}$.
\end{enumerate}
\end{fact}

\begin{fact}[Spectral properties of Kronecker product]\label{fac:kron_spectral}
The Kronecker product satisfies the following spectral properties.
\begin{enumerate}
    \item For matrices $A,B$, if $A$ and $B$ are PSD matrices, then $A\otimes B$ is also PSD.
    \item For two PSD matrices $A$ and $B\in \R^{n\times n}$, if $A\preceq B$, then $A\otimes A\preceq B\otimes B$.
\end{enumerate}
\end{fact}

The following result is often used in SDP-related calculations.
\begin{fact}[Kronecker product and vector multiplication]\label{fac:kron_vec}
Given $A \in \R^{m \times n},B\in \R^{n \times k}, C \in \R^{k \times l}$, $D \in \R^{l \times m}$, we have
\begin{enumerate}
    \item $\vect[ABC] = (C^\top \otimes A) \cdot \vect[B]$. 
    
    Note that $ABC \in \R^{m \times l}$, $C^\top \otimes A \in \R^{ml \times nk}$, and $\vect[B] \in \R^{nk}$.
    \item $\tr[ABCD] = \vect[D]^{\top} \cdot (C^\top \otimes A) \cdot \vect[B]$.
    
    Note that $ABCD \in \R^{m \times m}$, $\vect[D] \in \R^{ml}$, $C^\top \otimes A \in \R^{ml \times nk}$, and $\vect[B] \in \R^{nk}$.
\end{enumerate}
\end{fact}

We state a standard fact for Kronecker product.
\begin{fact}[Positive Semidefinite property of Kronecker product]\label{lem:approx_Hessian}
Let $m,n$ denote two positive integers.
Given a matrix $\mathsf{A} \in \R^{m \times n^2}$, let $S, \wt{S} \in \R^{n \times n}$ be two PSD matrices. 
Define
\begin{align*}
    H := \mathsf{A} \cdot (S^{-1} \otimes S^{-1}) \cdot \mathsf{A}^{\top} \in \R^{m \times m}, \text{~~~and~~~} \wt{H}:= \mathsf{A} \cdot (\wt{S}^{-1} \otimes \wt{S}^{-1}) \cdot \mathsf{A}^{\top} \in \R^{m \times m}.
\end{align*}
Then, for any accuracy parameter $\alpha \geq 1$, if $\wt{S}$ is an $\alpha$-PSD approximation of $S$, i.e., $\alpha^{-1} S \preceq \wt{S} \preceq \alpha S$, then 
\begin{align*}
    \alpha^{-2} H \preceq \wt{H} \preceq \alpha^2 H.
\end{align*}
\end{fact}

\begin{proof}
Given any vector $v \in \R^m$, we can write $v^\top H v$ and $v^\top \wt{H} v$ in the following way:
\begin{align}
    v^\top H v = & ~ \sum_{i = 1}^m\sum_{j=1}^m v_i v_jH_{i,j} = \sum_{i = 1}^m\sum_{j=1}^m v_i v_j \tr[S^{-1}A_iS^{-1}A_j] \notag\\
    = & ~ \tr \Big[ S^{-1/2} \Big( \sum_{i\in [m]} v_i A_i \Big) S^{-1} \Big( \sum_{i\in [m]} v_i A_i \Big) S^{-1/2} \Big]\notag\\
    = &~ \left\|\vect\Big[ S^{-1/2}\Big(\sum_{i\in[m]} v_i A_i \Big)
    S^{-1/2}\Big]\right\|_2^2\notag\\
    = &~\left\|\Big(S^{-1/2}
    \otimes S^{-1/2}\Big) \vect\Big[\sum_{i=1}^m v_i A_i\Big] \right\|_2^2,
    \label{eq:Hessian_norm}
\end{align}
where the last line follows from Fact~\ref{fac:kron_vec}. 
Similarly,
\begin{align}
    v^\top \wt{H} v = & ~ \left\|\Big(\widetilde{S}^{-1/2}
    \otimes \widetilde{S}^{-1/2}\Big) \vect\Big[\sum_{i=1}^m v_i A_i\Big] \right\|_2^2.
    \label{eq:approx_Hessian_norm}
\end{align}
Since the right hand side of Eq.~\eqref{eq:Hessian_norm} and Eq.~\eqref{eq:approx_Hessian_norm} are non-negative for any $v\in \R^m$, both $H$ and $\wt{H}$ are PSD matrices. 

Since $\alpha^{-1} S \preceq \wt{S} \preceq \alpha S$, we have $$\alpha^{-1} S^{-1} \preceq \wt{S}^{-1} \preceq \alpha S^{-1}.$$
By Fact~\ref{fac:kron_spectral}, it further implies that
\begin{align*}
    \alpha^{-2}S^{-1}\otimes S^{-1}\preceq \widetilde{S}^{-1}\otimes \widetilde{S}^{-1} \preceq \alpha^2S^{-1}\otimes S^{-1}.
\end{align*}
Let $b:=\vect[\sum_{i=1}^m v_i A_i]^\top$. We have
\begin{align}\label{eq:trace_upper_bound}
    \left\|\Big(\widetilde{S}^{-1/2}
    \otimes \widetilde{S}^{-1/2}\Big) b \right\|_2^2 = &~b^\top (\widetilde{S}^{-1/2}
    \otimes \widetilde{S}^{-1/2}) \cdot (\widetilde{S}^{-1/2}
    \otimes \widetilde{S}^{-1/2})b\notag\\
    = &~ b^\top (\widetilde{S}^{-1} \cdot \widetilde{S}^{-1})b\notag\\
    \leq &~ \alpha^2 \cdot b^\top (S^{-1}\otimes S^{-1})b\notag\\
    = &~ \alpha^2\cdot \left\|\Big(S^{-1/2}
    \otimes S^{-1/2}\Big) b \right\|_2^2.
\end{align}
And
\begin{align}
    \left\|\Big(\widetilde{S}^{-1/2}
    \otimes \widetilde{S}^{-1/2}\Big) b \right\|_2^2 \geq ~\alpha^{-2} \cdot \left\|\Big(S^{-1/2}
    \otimes S^{-1/2}\Big) b \right\|_2^2.
    \label{eq:trace_lower_bound}
\end{align}
Combining Eqs.~\eqref{eq:Hessian_norm}-\eqref{eq:trace_lower_bound}, we come to 
\begin{align*}
    \alpha^{-2}\cdot v^\top H v \leq v^\top \wt{H} v \leq \alpha^2\cdot v^\top H v.
\end{align*}
Since $v$ can be arbitrarily chosen from $\R^m$, we complete the proof.
\end{proof}

We state another fact for Kronecker product in below:
\begin{fact}[Kronecker product with equivalence for matrix norm]\label{fac:transform_norm_to_matrix_norm}
    Given a constraint matrix $\mathsf{A} \in \R^{m \times n^2}$ and vector $b \in \R^m$. Let $\eta >0$ denote a parameter. Let $g(y,\eta) \in \R^m$  be defined as
    \begin{align*}
        g(y,\eta)_i = \eta b_i - \tr[ S^{-1} A_i ]~~~\forall i\in [m].
    \end{align*}
    Let $X \in \R^{n \times n}$ denote a matrix that 
    \begin{align*}
        \langle X, A_i\rangle = \eta b_i~~~\forall i\in [m].
    \end{align*}
    Let $H:= \mathsf{A} (S^{-1} \otimes S^{-1}) \mathsf{A}^\top $.
    If matrix $S$ is a PSD matrix, then we have
    \begin{align*}
        g(y,\eta)^\top H^{-1} g(y,\eta) = v^\top \mathsf{B}^\top ( \mathsf{B} \mathsf{B}^\top )^{-1} \mathsf{B} v,
    \end{align*}
    where $v :=\vect[S^{1/2} X S^{1/2} -I]\in \R^{n^2}$ and $\mathsf{B} \in \R^{m \times n^2}$ is a matrix that $i$-th row is $B_i = \vect[S^{-1/2} A_i S^{-1/2}] \in \R^{n^2}$
\end{fact}
\begin{proof}
We start with re-writing $g(y,\eta) \in \R^m$ as follows: for each $i \in [m]$
\begin{align*}
    g(y,\eta)_i = & ~ b_i \eta - \tr[S^{-1} A_i] \\
    = & ~ \tr[X A_i] - \tr[S^{-1} A_i] \\
    = & ~ \tr[ (X - S^{-1}) A_i ] \\
    = & ~ \tr[ S^{1/2} (X - S^{-1}) S^{1/2} \cdot S^{-1/2} A_i S^{-1/2} ] \\
    = & ~ \tr[ (S^{1/2} X S^{1/2} - I) \cdot B_i ].
\end{align*}
Thus, using the definition of $v$, we have
\begin{align*}
    g(y,\eta) = \mathsf{B} v.
\end{align*}
Our next step is to rewrite $H$ as follows: for each $i,j \in [m]\times [m]$
\begin{align*}
    H_{i,j} = & ~ \tr[ A_i S^{-1} A_j S^{-1} ] \\
    = & ~ \tr[ S^{-1/2} A_i S^{-1/2} \cdot S^{-1/2} A_j S^{-1/2} ] \\
    = & ~ \tr[ B_i \cdot B_j ] 
\end{align*}
which implies that $H = \mathsf{B} \mathsf{B}^\top$.

Thus, combine all the above computations, we have
\begin{align*}
    g(y,\eta)^\top H^{-1} g(y,\eta) = v^\top \mathsf{B}^\top (\mathsf{B} \mathsf{B}^\top)^{-1} \mathsf{B} v.
\end{align*}
Therefore, we complete the proof.
\end{proof}

\subsection{Tools: Fast matrix multiplication}\label{sec:fast_matrix_multiplication}
We use $\Tmat(a,b,c)$ to denote the time of multiplying an $a \times b$ matrix with another $b \times c$ matrix. Fast  matrix multiplication \cite{c82,w12,l14,gu18,cglz20,aw21} is a fundamental tool in theoretical computer science.

For $k \in \R_+$, we define $\omega(k) \in \R_+$ to be the value such that $\forall n \in \N_+$, $\Tmat(n,n,n^k) = O(n^{\omega(k)})$.

For convenience we define three special values of $\omega(k)$. We define $\omega$ to be the fast matrix multiplication exponent, i.e., $\omega := \omega(1)$. We define $\alpha \in \R_+$ to be the dual exponent of %quadratic-time 
matrix multiplication, i.e., $\omega(\alpha) = 2$. We define $\beta := \omega(2)$.

The following fact can be found in Lemma 3.6 of \cite{jklps20}, also see \cite{bcs97}.
\begin{fact}[Convexity of $\omega(k)$]\label{fact:omega_k_convex}
The function $\omega(k)$ is convex.
\end{fact}

The following fact can be found in Lemma~A.5 of \cite{cls19}.%Upper bound of $\Tmat(n,n,r)$
\begin{fact}[Fast rectangular matrix multiplication]\label{fact:approx_Tmat_r_leq_n}
For any two integers $r \leq n$, the time of multiplying an $n \times n$ matrix with another $n \times r$
\begin{align*}
    \Tmat(n,n,r) \leq n^2 + r^{\frac{\omega - 2}{1 - \alpha}} \cdot n^{2 - \frac{\alpha (\omega-2)}{(1 - \alpha)}}.
\end{align*}
\end{fact}

The following fact can be found in Lemma~A.4 of \cite{cls19}.
\begin{fact}[Relation of $\omega$ and $\alpha$]\label{fact:bound_omega_alpha}
$\frac{\omega-2}{1-\alpha} - 1 \leq 0$; that is, $\omega + \alpha \leq 3$.
\end{fact}

\iffalse
The following fact can be found in Lemma~3.9-3.10 of \cite{jklps20}.
\begin{fact}[Upper bound of $\Tmat (m,n^2,m)$]\label{fact:bound_tmat_m_n2_m}
Let $m$ and $n$ denote two positive integers. If $n \geq m$, we have
\begin{align*}
    \Tmat (m,n^2,m) = O(mn^{\omega +o(1)}).
\end{align*}
\end{fact}
\fi

\section{Our Algorithm and Result}\label{sec:alg_and_result}
%\subsection{Main result}

We state our main result of Algorithm~\ref{alg:ss_ipm_woodbury} as follows:

\begin{theorem}[Main result for Algorithm~\ref{alg:ss_ipm_woodbury}]\label{thm:general_SDP_formal}
Given symmetric matrices $C, A_1,\cdots,A_m \in \R^{n \times n}$, and a vector $b \in \R^m$. Define matrix $\mathsf{A} \in \R^{m \times n^2}$ by stacking the $m$ vectors $\vect[A_1], \cdots, \vect[A_m] \in \R^{n^2}$ as rows. Consider the following SDP instance:
\begin{align*}
     \max_{X \in \R^{n \times n}} & ~ \langle C,X \rangle \\
    \mathrm{~s.t.~} &~ \langle A_i, X \rangle = b_i, ~ \forall i \in [m], \\
    &~ X \succeq 0,
\end{align*}
There is a SDP algorithm (Algorithm~\ref{alg:ss_ipm_woodbury}) that runs 
in time
\begin{align*}
    O^\ast\Big( \left( \sqrt{n} ( m^{2} + n^4 ) + m^{\omega}  + n^{2 \omega}\right) \cdot \log (1/\epsilon) \Big) .
\end{align*}
and outputs a PSD matrix $X \in \mathbb{R}^{n \times n}$ that satisfies 
\begin{align}\label{eq:approx_optimality_formal}
\langle C, X \rangle \geq \langle C, X^* \rangle - \epsilon \cdot \| C \|_{2} \cdot R \quad \text{and} \quad \sum_{i = 1}^m\left| \langle A_i, X \rangle - b_i \right| \leq 4 n \epsilon \cdot \Big( R  \sum_{i = 1}^m \| A_i \|_1 + \| b \|_1 \Big) ,
\end{align}
where $X^*$ is an optimal solution of the SDP instance, and $\| A_i \|_1$ is the Schatten $1$-norm of matrix $A_i$.
\end{theorem}
\begin{proof}
The correctness (Eq.~\eqref{eq:approx_optimality_formal}) follows from Theorem~\ref{thm:correctness_ss_ipm_woodbury}. The running time follows from Theorem~\ref{thm:running_time}.
\end{proof}

\begin{remark}\label{remark:running_time}
For current  matrix multiplication time $\omega \approx 2.373$ (\cite{l14}), the running time of our algorithm can be written as 
\begin{align*}
    O(\max\{n^{2\omega},m^{\omega}\} \cdot \log (1/\epsilon)).
\end{align*}
Therefore, when $m \geq n^{2-0.5/\omega} \approx n^{1.79}$, we have $\max\{n^{2\omega},m^{\omega}\} \leq \sqrt{n}\cdot m^\omega$ and thus our algorithm is better than \cite{jklps20}. For the regimes when $m$ is smaller, we can apply Algorithm~\ref{alg:sdp_vol}-\ref{alg:sdp_vol_cont} in Section~\ref{sec:hybrid_barrier} or Algorithm~1 in \cite{jklps20}. 
\end{remark}

\begin{corollary}[Tall SDPs]
When $m = \Omega(n^2)$, we can solve SDP in $O(m^{\omega} \cdot \log(1/\epsilon))$ time for current $\omega \approx 2.373$.
\end{corollary}
\begin{proof}
When $m = \Omega(n^2)$, the running time of Theorem~\ref{thm:general_SDP_formal} is
\begin{align*}
    O\Big( (m^{2} \cdot \sqrt{n} + m^{\omega} + n^{4.5} + n^{2 \omega}) \cdot \log (1/\epsilon) \Big) = O((m^2 \cdot \sqrt{n} + m^{\omega}) \cdot \log (1/\epsilon)).
\end{align*}
For current $\omega \approx 2.373$ , $m^2 \cdot \sqrt{n} = m^{2.25} < m^{\omega}$.
\end{proof}
Note that the running time of \cite{jklps20} is $O(\sqrt{n} \cdot (n^{\omega} + m^{\omega} + m n^2)) = O(m^{\omega + 0.25})$ when $m = \Omega(n^2)$.

\begin{algorithm}[!ht]
\caption{ Our SDP solver with log barrier.}
\label{alg:ss_ipm_woodbury}
\begin{algorithmic}[1]
\item {\bf procedure} {\textsc{SolveSDP}}{$(m, n, C, \{ A_i \}_{i=1}^m$, $\mathsf{A} \in \R^{m \times n^2}$, $b \in \R^m)$} \\
\hspace{4mm} \hspace{\fill} $\triangleright$ {Initialization} \\
\hspace{4mm} Construct $\mathsf{A} \in \R^{m \times n^2}$ by stacking $m$ vectors $\vect[A_1], \vect[A_2], \cdots, \vect[A_m] \in \R^{n^2}$ \label{line:mathsf_A}\\
\hspace{4mm} $\eta \leftarrow \frac{1}{n+2}$, ~~$T \leftarrow \frac{40}{\epsilon_N} \sqrt{n} \log (\frac{n}{\epsilon})$\\
\hspace{4mm} Find initial feasible dual vector $y \in \R^m$ according to Lemma~\ref{lem:initialization}\\
\hspace{4mm} $S \leftarrow \sum_{i \in [m]} y_i \cdot A_i - C$, ~~ $\wt{S} \leftarrow S$ \hspace{\fill} $\triangleright$ {$S, \wt{S} \in \R^{n \times n}$} \label{line:S_init}\\
\hspace{4mm} $G \leftarrow (\mathsf{A} \cdot (\wt{S}^{-1} \otimes \wt{S}^{-1}) \cdot \mathsf{A}^{\top})^{-1}$ \hspace{\fill} $\triangleright$ {$G \in \R^{m \times m}$} \label{line:G_init}\\
\hspace{4mm} \hspace{\fill} $\triangleright$ {Maintain $G = \wt{H}^{-1}$ where $\wt{H} := \mathsf{A} \cdot (\wt{S}^{-1} \otimes \wt{S}^{-1}) \cdot \mathsf{A}^{\top}$}\\
\hspace{4mm}{\bf for} {$t = 1 \to T$} {\bf do} \hspace{\fill} $\triangleright$ {Iterations of approximate barrier method}\\
	\hspace{8mm} $\eta^{\new} \leftarrow \eta \cdot (1 + \frac{\epsilon_N}{20 \sqrt{n}})$\\
	\hspace{8mm} {\bf for} {$j = 1, \cdots, m$} {\bf do} \\
		\hspace{12mm} $g_{\eta^{\new}}(y)_j \leftarrow  b_j \cdot\eta^{\new}  -   \tr [ S^{-1} \cdot A_j ] $ \hspace{\fill} $\triangleright$ {Gradient computation, $g_{\eta^{\new}}(y) \in \R^m$} \label{line:gradient_sdp_woodbury}\\
	\hspace{8mm} {\bf end for}\\
	\hspace{8mm} $\delta_y \leftarrow - G \cdot g_{\eta^{\new}}(y)$ \hspace{\fill} $\triangleright$ {Update on $y \in \R^m$}\label{line:delta_y_sdp_woodbury}\\
	%\hspace{8mm} Let $\wt{\delta}_y$ denote the changes generated by Quantum machine \Zhao{I add this}\\
	\hspace{8mm} $y^{\new} \leftarrow y + \delta_y$\\
	\hspace{8mm} $S^{\new} \leftarrow \sum_{i \in [m]} (y^{\new})_i \cdot A_i - C$ \label{line:S_new_sdp_woodbury}\\
	\hspace{8mm} $V_1, V_2 \leftarrow \textsc{LowRankSlackUpdate}(S^{\new}, \wt{S})$ \hspace{\fill} $\triangleright$ {$V_1,V_2 \in \R^{n \times r_t}$. Algorithm~\ref{alg:low_rank_slack_update}.} \label{line:low_rank_slack_update} \\
	\hspace{8mm} $\wt{S}^{\new} \leftarrow \wt{S} + V_1 V_2^{\top}$ \hspace{\fill} $\triangleright$ {Approximate slack computation}\label{line:slack_update} \\
	\hspace{8mm} $V_3 \leftarrow - \wt{S}^{-1} V_1 (I + V_2^{\top} \wt{S}^{-1} V_1)^{-1}$ \Comment{$V_3 \in \R^{n \times r_t}$} \label{line:V_3}\\
	\hspace{8mm} $V_4 \leftarrow \wt{S}^{-1} V_2$ \hspace{\fill} $\triangleright$ {$V_4 \in \R^{n \times r_t}$} \label{line:V_4}\\
	\hspace{8mm} $Y_1 \leftarrow [(\wt{S}^{-1/2} \otimes V_3), (V_3 \otimes \wt{S}^{-1/2}), (V_3 \otimes V_3^{\top})]$ \hspace{\fill} $\triangleright$ {$Y_1 \in \R^{n^2 \times (2n r_t + r_t^2)}$} \label{line:Y_1}\\
	\hspace{8mm} $Y_2 \leftarrow [(\wt{S}^{-1/2} \otimes V_4), (V_4 \otimes \wt{S}^{-1/2}), (V_4 \otimes V_4^{\top})]$ \hspace{\fill} $\triangleright$ {$Y_2 \in \R^{n^2 \times (2n r_t + r_t^2)}$} \label{line:Y_2}\\
	\hspace{8mm} $G^{\new} \leftarrow G - G \cdot \mathsf{A} Y_1 \cdot (I + Y_2^{\top} \mathsf{A}^{\top} \mathsf{A} Y_1)^{-1} \cdot Y_2^{\top} \mathsf{A}^{\top} \cdot G$ \Comment{$G^{\new} \in \R^{m \times m}$} \label{line:hessian_sdp_woodbury}\\
	\hspace{8mm} \hspace{\fill} $\triangleright$ {Hessian inverse computation using Woodbury identity}\\
	\hspace{8mm} $y \leftarrow y^{\new}$ \\
	\hspace{8mm} $S \leftarrow S^{\new}$ \\
	\hspace{8mm} $\wt{S} \leftarrow \wt{S}^{\new}$ \\
	\hspace{8mm} $G \leftarrow G^{\new}$ \hspace{\fill} $\triangleright$ {Update variables}\\
\hspace{4mm} {\bf end for}\\
\hspace{4mm} \Return an approximate solution to the original problem \Comment{ Lemma~\ref{lem:initialization}}\\
{\bf end procedure}
\end{algorithmic}
\end{algorithm}

\begin{algorithm}[!ht]
\caption{Low Rank Slack Update 
}\label{alg:low_rank_slack_update}
\begin{algorithmic}[1]
\item {\bf procedure} {\textsc{LowRankSlackUpdate}}{$(S^{\new}, \wt{S})$} \\
\hspace{4mm} \hspace{\fill} $\triangleright$ {$S^{\new}, \wt{S} \in \mathbb{S}^{n \times n}_{\geq 0}$ are positive definite matrices}\\
\hspace{4mm} $\epsilon_S \leftarrow 10^{-5}$ \hspace{\fill} $\triangleright$ {Spectral approximation constant}\\
\hspace{4mm} $Z^{\midd} \leftarrow (S^{\new})^{-1/2} \cdot \wt{S} \cdot (S^{\new})^{-1/2} - I_n$ \label{line:Z_mid}\\
\hspace{4mm} Compute spectral decomposition $Z^{\midd} = U \cdot \diag(\lambda) \cdot U^\top$ \label{line:Z_mid_decomposition}\\
\hspace{4mm} \hspace{\fill} $\triangleright$ {$\lambda = [\lambda_1,\cdots, \lambda_n]^{\top} \in \R^n$ are the eigenvalues of $Z^{\midd}$, and $U \in \R^{n \times n}$ is orthogonal}\\
\hspace{4mm} Let $\pi: [n] \rightarrow [n]$ be a sorting permutation such that $|\lambda_{\pi(i)}| \geq |\lambda_{\pi(i+1)}|$ \label{line:pi}\\
\hspace{4mm} {\bf if} {$|\lambda_{\pi(1)}| \leq \epsilon_S$} {\bf then}\\
	\hspace{8mm} $\wt{S}^{\new} \leftarrow \wt{S}$ \\
\hspace{4mm} {\bf else}\\
	\hspace{8mm} $r \leftarrow 1$\\
	\hspace{8mm} {\bf while} {$r \leq n/2$ and $(|\lambda_{\pi(2r)}| > \epsilon_S$ or $|\lambda_{\pi(2r)}| > (1 - 1/\log n) |\lambda_{\pi(r)}|)$} {\bf do} \label{line:while_slack_update}\\
		\hspace{12mm} $r \leftarrow r + 1$\\
	\hspace{8mm} {\bf end while}\\

    \hspace{8mm} $(\lambda^{\new})_{\pi(i)} \leftarrow
		\begin{cases}
		 0, & \text{~if~} i = 1, 2, \cdots, 2r; \\
		 \lambda_{\pi(i)}, & \text{~otherwise.}
		\end{cases}$\label{line:eigenvalue_update}\\
	\hspace{8mm} $L \leftarrow \supp(\lambda^{\new} - \lambda)$ \hspace{\fill} $\triangleright$ {$|L| = 2r$}\\
	\hspace{8mm} $V_1 \leftarrow ((S^{\new})^{1/2} \cdot U \cdot \diag (\lambda^{\new} - \lambda) )_{:,L}$ \hspace{\fill} $\triangleright$ {$V_1 \in \R^{n \times 2r}$}\label{line:V_1}\\
	\hspace{8mm} $V_2 \leftarrow ((S^{\new})^{1/2} \cdot U )_{:,L}$ \hspace{\fill} $\triangleright$ {$V_2 \in \R^{n \times 2r}$}\label{line:V_2}\\
	\hspace{8mm} \hspace{\fill} $\triangleright$ {$V_1 \cdot V_2^{\top} = (S^{\new})^{1/2} \cdot U \cdot \diag (\lambda^{\new} - \lambda) \cdot U^\top \cdot (S^{\new})^{1/2}$}\\
\hspace{4mm} {\bf end if}\\
\hspace{4mm} \Return $\wt{S}^{\new}$\\
{\bf end procedure}
\end{algorithmic}
\end{algorithm}

\section{Correctness}\label{sec:correct}

In this section we prove the correctness of our  SDP solver Algorithm~\ref{alg:ss_ipm_woodbury}. 
In Section~\ref{sec:approx_slack} we prove that $\wt{S} \in \R^{n \times n}$ updated by Algorithm~\ref{alg:low_rank_slack_update} is a PSD approximation to the true slack matrix $S \in \R^{n \times n}$. In Section~\ref{sec:approx_hessian} we prove that the algorithm maintains $G = \wt{H}^{-1} \in \R^{m \times m}$, and $\wt{H}$ is a PSD approximation to the true Hessian matrix $H \in \R^{m \times m}$.

\begin{theorem}[Correctness of Algorithm~\ref{alg:ss_ipm_woodbury}]\label{thm:correctness_ss_ipm_woodbury}
Consider an SDP instance as in Definition~\ref{def:sdp_primal} with no redundant constraints. Let us assume that the feasible region is bounded, i.e., any feasible solution $X \in \R^{n \times n}_{\succeq 0}$ satisfies $\| X \|_{2} \leq R$. Then for any error parameter $0 < \epsilon \leq 0.01$ and Newton step size $\epsilon_N$ satisfying  $\sqrt{\epsilon} < \epsilon_N \leq 0.1$, Algorithm~\ref{alg:ss_ipm_woodbury} outputs,  in $T = 40 \epsilon_N^{-1} \sqrt{n} \log(n/\epsilon)$ iterations,  a PSD matrix $X \in \mathbb{R}^{n \times n}_{\succeq 0}$ that satisfies 
\begin{align}\label{eq:approx_optimality}
\begin{array}{l}
\langle C, X \rangle \geq \langle C, X^* \rangle - \epsilon \cdot \| C \|_2 \cdot R,~~ \text{and}\\ \sum_{i = 1 }^m \left| \langle A_i, X \rangle - b_i \right| \leq 4 n \epsilon \cdot \Big( R  \sum_{i = 1}^m \| A_i \|_1 + \| b \|_1 \Big) ,
\end{array}
\end{align}
where $X^*$ is any optimal solution to the SDP instance, and $\| A_i \|_1$ is the Schatten $1$-norm of matrix $A_i$. 

Furthermore, in each iteration of Algorithm~\ref{alg:ss_ipm_woodbury}, the following invariant holds for $\alpha_H = 1+10^{-5}$: 
\begin{align}
\| S^{-1/2} S^{\new} S^{-1/2} - I \|_F \leq \alpha_H \cdot \epsilon_N. 
\end{align}
\end{theorem}

\begin{proof}
Combining Lemma~\ref{lem:close_form_slack} and Lemma~\ref{lem:approx_slack} we have that $\alpha_S^{-1} S \preceq \wt{S} \preceq \alpha_S S$ for $\alpha_S = 1+10^{-5}$. 
Therefore condition 1' in Lemma~\ref{lem:invariant_newton} is satisfied by Fact~\ref{lem:approx_Hessian} and condition 2 \& 3 holds trivially. %Moreover,  condition 3 and 4 in Lemma~\ref{lem:invariant_newton} trivially holds. 
Notice $\theta_{\phi_{\mathrm{log}}} = n$. 
Then directly applying Theorem~\ref{thm:approx_central_path} completes the proof.

\end{proof}

\subsection{Approximate slack maintenance}\label{sec:approx_slack}
The following lemma gives a closed-form formula for the updated $\wt{S}^{\new}$ in each iteration.
\begin{lemma}[Closed-form formula of slack update]\label{lem:close_form_slack}
In each iteration of Algorithm~\ref{alg:ss_ipm_woodbury}, the update of the slack variable $\wt{S}^{\new}$ (on Line~\ref{line:slack_update}) satisfies 
\begin{align*}
    \wt{S}^{\new} = \wt{S} + (S^{\new})^{1/2} \cdot U \cdot \diag (\lambda^{\new} - \lambda) \cdot U^\top \cdot (S^{\new})^{1/2},
\end{align*}
where $\widetilde{S}$ is the slack variable in previous iteration, and $U, \lambda, \lambda^{\new}$ are defined in Algorithm~\ref{alg:low_rank_slack_update}. 

Moreover, it implies that $\wt{S}^{\new}$ is a symmetric matrix in each iteration.
\end{lemma}
\begin{proof}
From Line~\ref{line:V_1} and \ref{line:V_2} of Algorithm~\ref{alg:low_rank_slack_update} we have $V_1 \cdot V_2^{\top} = (S^{\new})^{1/2} \cdot U \cdot \diag (\lambda^{\new} - \lambda) \cdot U^\top \cdot (S^{\new})^{1/2}$. Therefore
\begin{align*}
    \wt{S}^{\new} = \wt{S} + V_1 \cdot V_2^{\top} = \wt{S} + (S^{\new})^{1/2} \cdot U \cdot \diag (\lambda^{\new} - \lambda) \cdot U^\top \cdot (S^{\new})^{1/2}.
\end{align*}

In the first iteration, we have $\widetilde{S}=S$ which is a symmetric matrix. 

By the definition of $V_1, V_2$, we know that $V_1\cdot V_2^\top$ is symmetric. Hence, $S^{\new}$ is also symmetric in each iteration.
\end{proof}

The following lemma proves that we always have $\wt{S} \approx S$ throughout the algorithm.
\begin{lemma}[Approximate Slack]\label{lem:approx_slack}
In each iteration of Algorithm~\ref{alg:ss_ipm_woodbury}, the approximate slack variable $\wt{S}$ satisfies that $\alpha_S^{-1} S \preceq \wt{S} \preceq \alpha_S S$, where $\alpha_S = 1+10^{-5}$.
\end{lemma}
\begin{proof}
Notice that
\begin{align*}
    \wt{S}^{\new} &= \wt{S} + (S^{\new})^{1/2} \cdot U \cdot \diag (\lambda^{\new} - \lambda) \cdot U^\top \cdot (S^{\new})^{1/2}\\
    &= \left({S}^{\new} + (S^{\new})^{1/2} Z^{\midd} (S^{\new})^{1/2}\right) + (S^{\new})^{1/2} \cdot U \cdot \diag (\lambda^{\new} - \lambda) \cdot U^\top \cdot (S^{\new})^{1/2}\\
    &= {S}^{\new} + (S^{\new})^{1/2} \cdot U \cdot \diag (\lambda^{\new}) \cdot U^\top \cdot (S^{\new})^{1/2},
\end{align*}
where the first step comes from Lemma~\ref{lem:close_form_slack}, the second step comes from definition $Z^{\midd} = (S^{\new})^{-1/2} \cdot \wt{S} \cdot (S^{\new})^{-1/2} - I$ (Line~\ref{line:Z_mid} of Algorithm~\ref{alg:low_rank_slack_update}), and the final step comes from $Z^{\midd} = U \cdot \diag(\lambda_1,\cdots, \lambda_n) \cdot U^\top$ (Line~\ref{line:Z_mid_decomposition} of Algorithm~\ref{alg:low_rank_slack_update}). 

By Line~\ref{line:eigenvalue_update} of Algorithm~\ref{alg:low_rank_slack_update} we have $(\lambda^{\new})_i \leq \epsilon_S$ for all $i \in [n]$, so 
\begin{align*}
    \left\| (S^{\new})^{-1/2} \cdot \wt{S}^{\new} \cdot (S^{\new})^{-1/2} - I\right\|_2 = \left\|U \cdot \diag (\lambda^{\new}) \cdot U^\top\right\|_2 \leq \epsilon_S.
\end{align*}
This implies that for $\alpha_S = 1+\epsilon_S$, by Fact~\ref{fac:spec_psd}, in each iteration of Algorithm~\ref{alg:ss_ipm_woodbury} the slack variable $\wt{S}$ satisfies $\alpha_S^{-1} S \preceq \wt{S} \preceq \alpha_S S$.
\end{proof}

\subsection{Approximate Hessian inverse maintenance}\label{sec:approx_hessian}
The following lemma shows that the maintained matrix $G$ equals to the inverse of approximate Hessian.
\begin{lemma}[Close-form formula for Hessian inverse]\label{lem:close_form_hessian_inverse}
In each iteration of Algorithm~\ref{alg:ss_ipm_woodbury}, we have $G = \wt{H}^{-1} \in \R^{m \times m}$, where $\wt{H} := \mathsf{A} \cdot (\wt{S}^{-1} \otimes \wt{S}^{-1}) \cdot \mathsf{A}^{\top} \in \R^{m \times m}$.
\end{lemma}
\begin{proof}
We prove this lemma by induction.

In the beginning of the algorithm, the initialization of $G$ (Line~\ref{line:G_init} of Algorithm~\ref{alg:ss_ipm_woodbury}) satisfies the formula $G = \wt{H}^{-1}$. 

Assume the induction hypothesis that $G = \wt{H}^{-1}$ in the beginning of each iteration, next we will prove that $G^{\new} = (\wt{H}^{\new})^{-1}$. Note that $G^{\new}$ is updated on Line~\ref{line:hessian_sdp_woodbury} of Algorithm~\ref{alg:ss_ipm_woodbury}. And $\wt{H}^{\new} := \mathsf{A} \cdot ((\wt{S}^{\new})^{-1} \otimes (\wt{S}^{\new})^{-1}) \cdot \mathsf{A}^{\top} \in \R^{m \times m}$, where $\wt{S}^{\new}$ is updated on Line~\ref{line:slack_update} of Algorithm~\ref{alg:ss_ipm_woodbury}.

We first compute $(\wt{S}^{\new})^{-1} \in \R^{n \times n}$:
\begin{align}\label{eq:wt_S_new_inverse}
    (\wt{S}^{\new})^{-1} = &~ (\wt{S} + V_1 V_2^{\top})^{-1} \notag \\
    = &~ \wt{S}^{-1} - \wt{S}^{-1} V_1 \cdot (I + V_2^{\top} \wt{S}^{-1} V_1)^{-1} \cdot V_2^{\top} \wt{S}^{-1} \notag \\
    = &~ \wt{S}^{-1} + V_3 \cdot V_4^{\top},
\end{align}
where the reason of the first step is $\wt{S}^{\new} = \wt{S} + V_1 V_2^{\top}$ (Line~\ref{line:slack_update} of Algorithm~\ref{alg:ss_ipm_woodbury}), the second step follows from Woodbury identity (Fact~\ref{fac:woodbury}), and the third step follows from $V_3 = - \wt{S}^{-1} V_1 (I + V_2^{\top} \wt{S}^{-1} V_1)^{-1} \in \R^{n \times r_t}$ and $V_4 = \wt{S}^{-1} V_2 \in \R^{n \times r_t}$ (Line~\ref{line:V_3} and \ref{line:V_4} of Algorithm~\ref{alg:ss_ipm_woodbury}).

We then compute a close-form formula of $(\wt{S}^{\new})^{-1} \otimes (\wt{S}^{\new})^{-1} \in \R^{n^2 \times n^2}$:
\begin{align}\label{eq:wt_S_new_kron}
    &~ (\wt{S}^{\new})^{-1} \otimes (\wt{S}^{\new})^{-1} \notag \\
    = &~ (\wt{S}^{-1} + V_3 V_4^{\top}) \otimes (\wt{S}^{-1} + V_3 V_4^{\top}) \notag \\
    = &~ \wt{S}^{-1} \otimes \wt{S}^{-1} + \wt{S}^{-1} \otimes (V_3 V_4^{\top}) + (V_3 V_4^{\top}) \otimes \wt{S}^{-1} + (V_3 V_4^{\top}) \otimes (V_3 V_4^{\top}) \notag \\
    = &~ \wt{S}^{-1} \otimes \wt{S}^{-1} + (\wt{S}^{-1/2} \otimes V_3) \cdot (\wt{S}^{-1/2} \otimes V_4^{\top}) + (V_3 \otimes \wt{S}^{-1/2}) \cdot (V_4^{\top} \otimes \wt{S}^{-1/2}) \notag \\
    &~ + (V_3 \otimes V_3) \cdot (V_4^{\top} \otimes V_4^{\top}) \notag \\
    = &~ \wt{S}^{-1} \otimes \wt{S}^{-1} + Y_1 Y_2^{\top},
\end{align}
where the first step follows from Eq.~\eqref{eq:wt_S_new_inverse}, the second step follows from linearity of Kronecker product, the third step follows from mixed product property of Kronecker product (Part~\ref{part:kron_mixed_product} of Fact~\ref{fac:kron_basic}), the fourth step follows from $Y_1 = [(\wt{S}^{-1/2} \otimes V_3), (V_3 \otimes \wt{S}^{-1/2}), (V_3 \otimes V_3^{\top})] \in \R^{n^2 \times (2n r_t + r_t^2)}$ and $Y_2 = [(\wt{S}^{-1/2} \otimes V_4), (V_4 \otimes \wt{S}^{-1/2}), (V_4 \otimes V_4^{\top})] \in \R^{n^2 \times (2n r_t + r_t^2)}$ (Line~\ref{line:Y_1} and \ref{line:Y_2} of Algorithm~\ref{alg:ss_ipm_woodbury}), and the transpose of Kronecker product (Part~\ref{part:kron_transpose} of Fact~\ref{fac:kron_basic}).

Thus we can compute $(\wt{H}^{\new})^{-1} \in \R^{m \times m}$ as follows:
\begin{align}\label{eq:H_new_woodbury}
    (\wt{H}^{\new})^{-1} = &~ \big( \mathsf{A} \cdot \big((\wt{S}^{\new})^{-1} \otimes (\wt{S}^{\new})^{-1} \big) \cdot \mathsf{A}^{\top} \big)^{-1} \notag \\
    = &~ \Big( \mathsf{A} \cdot (\wt{S}^{-1} \otimes \wt{S}^{-1}) \cdot \mathsf{A}^{\top} + \mathsf{A} \cdot Y_1 Y_2^{\top} \cdot \mathsf{A}^{\top} \Big)^{-1} \notag \\
    = &~ G - G \cdot \mathsf{A} Y_1 \cdot (I + Y_2^{\top} \mathsf{A}^{\top} \cdot \mathsf{A} Y_1)^{-1} \cdot Y_2^{\top} \mathsf{A}^{\top} \cdot G \notag \\
    = &~ G^{\new},
\end{align}
where the first step follows from the definition of $\wt{H}^{\new}$, the second step follows from Eq.~\eqref{eq:wt_S_new_kron}, the third step follows from Woodbury identity (Fact~\ref{fac:woodbury}) and the induction hypothesis that $G = \wt{H}^{-1} = (\mathsf{A} \cdot (\wt{S}^{-1} \otimes \wt{S}^{-1}) \cdot \mathsf{A}^{\top})^{-1} \in \R^{m \times m}$, the fourth step follows from the definition of $G^{\new}$ on Line~\ref{line:hessian_sdp_woodbury} of Algorithm~\ref{alg:ss_ipm_woodbury}.

The proof is then completed.
\end{proof}

\section{Time Analysis}\label{sec:time}
In this section we analyze the running time of Algorithm~\ref{alg:ss_ipm_woodbury}. We first present the main theorem.
\begin{theorem}[Running time of Algorithm~\ref{alg:ss_ipm_woodbury}]\label{thm:running_time}
Algorithm~\ref{alg:ss_ipm_woodbury} runs 
in time
\begin{align*}
    O^\ast\left( (m^{2} \cdot \sqrt{n} + m^{\omega} + n^{4.5} + n^{2 \omega}) \cdot \log (1/\epsilon) \right) .
\end{align*}
Specifically, when $m \geq n^2$ the total running time is
\begin{align*}
    O^\ast\left( ( m^{\omega} + m^{2} \cdot \sqrt{n}) \cdot \log (1/\epsilon) \right).
\end{align*}
\end{theorem}
\begin{proof}
The running time of Algorithm~\ref{alg:ss_ipm_woodbury} consists of two parts:

{\bf 1. Initialization.} $O(\Tmat(m, n^2, n^2) + m^{\omega}) \leq O^\ast(n^{2 \omega} + m^{\omega})$ time from Lemma~\ref{lem:initialization_time}.

{\bf 2. Cost of $T$ iterations.}
\begin{align*}
    &~ \sum_{t=1}^T \Big( \Tmat(m, n^2, n r_t) + \Tmat(m, m, n r_t) + (n r_t)^{\omega} + m^2 + m n^2 \Big) \\
    \leq &~ T \cdot O(m^2 + n^4 + n^{2 \omega - 1/2} + m^{2 - \frac{\alpha(\omega-2)}{1-\alpha}} \cdot n^{\frac{2(\omega - 2)}{1 - \alpha}-1/2}) \\
    = &~ O^\ast\Big( (m^{2} \cdot \sqrt{n} + n^{4.5} + n^{2 \omega} + m^{2 - \frac{\alpha(\omega-2)}{1-\alpha}} \cdot n^{\frac{2(\omega - 2)}{1 - \alpha}} ) \cdot \log (1/\epsilon) \Big) \\
    \leq &~ O^\ast\Big( (m^{2} \cdot \sqrt{n} + n^{4.5} + n^{2 \omega} + m^{\omega} ) \cdot \log (1/\epsilon) \Big)
\end{align*}
where the first step follows from Lemma~\ref{lem:cost_per_iter}, the second step follows from Lemma~\ref{lem:amortize_hessian} and $m n^2 \leq O(m^2 + n^4)$, the third step follows from $T = O(\sqrt{n} \log (1/\epsilon))$, the fourth step follows from the inequalities in two cases:
\begin{align*}
    m^{2 - \frac{\alpha(\omega-2)}{1-\alpha}} \cdot n^{\frac{2(\omega - 2)}{1 - \alpha}} \leq &~ \begin{cases}
    (n^2)^{2 - \frac{\alpha(\omega-2)}{1-\alpha}} \cdot n^{\frac{2(\omega - 2)}{1 - \alpha}} = n^{2 \omega} &\text{when~}m \leq n^2, \\
    m^{2 - \frac{\alpha(\omega-2)}{1-\alpha}} \cdot (m^{1/2})^{\frac{2(\omega - 2)}{1 - \alpha}} = m^{ \omega} &\text{when~}m > n^2.
    \end{cases}
\end{align*}

Adding the costs of these two parts completes the proof.
\end{proof}

This section is organized as follows:
\begin{itemize}
    \item Section~\ref{sec:initialization_time} provides the analysis of initialization cost (Lemma~\ref{lem:initialization_time}).
    \item Section~\ref{sec:cost_per_iter} studies the cost per iteration of our algorithm (Lemma~\ref{lem:cost_per_iter}).
    \item Section~\ref{sec:property_of_low_rank_update} present the property of low rank update (Theorem~\ref{thm:general_amortize_guarantee}).
    \item Section~\ref{sec:amortized_analysis} use the tools of Section~\ref{sec:property_of_low_rank_update} to bound the amortized cost of our algorithm (Lemma~\ref{lem:amortize_hessian}).
\end{itemize}

\begin{table}[]
\small
    \centering
    \begin{tabular}{|l|l|l|l|}
    \hline
        {\bf Sec.}& {\bf Statement} & {\bf Time} & {\bf Comment}  \\ \hline
        \ref{sec:initialization_time} & Lem~\ref{lem:initialization_time} & $\Tmat(m,n^2,n^2) + m^{\omega}$ & Initialization \\ \hline
        \ref{sec:cost_per_iter} & Lem~\ref{lem:cost_per_iter} & $\Tmat(m, n^2, n r_t) + \Tmat(m, m, n r_t) + (n r_t)^{\omega} + m^2 + m n^2$ & Cost per iteration \\ \hline
        \ref{sec:amortized_analysis} & Lem~\ref{lem:amortize_hessian} & $m^2 + n^4 + n^{2 \omega - 1/2} + m^{2 - \frac{\alpha(\omega-2)}{1-\alpha}} \cdot n^{\frac{2(\omega - 2)}{1 - \alpha}}$ & Amortized cost \\ \hline
    
    \end{tabular}
    \caption{Summary of Section~\ref{sec:time}. 
    }
    \label{tab:summary_of_time}
\end{table}

\subsection{Initialization cost}\label{sec:initialization_time}
The goal of this section is to prove Lemma~\ref{lem:initialization_time}.
\begin{lemma}[Initialization cost]\label{lem:initialization_time}
Initialization of Algorithm~\ref{alg:ss_ipm_woodbury} (Line~\ref{line:mathsf_A} to \ref{line:G_init}) takes time
\begin{align*}
O(\Tmat(m,n^2,n^2) + m^{\omega}).
\end{align*}
\end{lemma}
\begin{proof}
We compute the cost of each line during initialization.
\begin{itemize}
    \item {\bf Line~\ref{line:mathsf_A} of Algorithm~\ref{alg:ss_ipm_woodbury}.} Constructing $\mathsf{A} \in \R^{m \times n^2}$ by stacking vectors takes $O(m n^2)$ time.
    \item {\bf Line~\ref{line:S_init} of Algorithm~\ref{alg:ss_ipm_woodbury}.} Computing $S = \sum_{i \in [m]} y_i \cdot A_i - C$ takes $O(m n^2)$ time.
    \item {\bf Line~\ref{line:G_init} of Algorithm~\ref{alg:ss_ipm_woodbury}.} This step computes $G = (\mathsf{A} \cdot (\wt{S}^{-1} \otimes \wt{S}^{-1}) \cdot \mathsf{A}^{\top})^{-1}$. We first compute $\wt{S}^{-1} \otimes \wt{S}^{-1} \in \R^{n^2 \times n^2}$, which takes $O(n^4)$ time. We then compute $\mathsf{A} \cdot (\wt{S}^{-1} \otimes \wt{S}^{-1}) \cdot \mathsf{A}^{\top}$, which takes $O(\Tmat(m, n^2,n^2) + \Tmat(m,n^2,m))$ time. Finally computing the inverse takes $O(m^{\omega})$ time. Since $n^4 \leq \Tmat(m,n^2,n^2)$ and $\Tmat(m,n^2,m) \leq \Tmat(m,n^2,n^2) + m^{\omega}$, this step takes $O(\Tmat(m,n^2,n^2) + m^{\omega})$ time in total, 
\end{itemize}
Combining the cost of these three steps, and since $mn^2 \leq \Tmat(m,n^2,n^2)$, we have the total cost as presented in the lemma statement.
\end{proof}

\subsection{Cost per iteration}\label{sec:cost_per_iter}
The goal of this section is to prove Lemma~\ref{lem:cost_per_iter}.
\begin{lemma}[Cost per iteration]\label{lem:cost_per_iter}
For $t \in [T]$, let $r_t$ be the rank of the update in the $t$-th iteration of Algorithm~\ref{alg:ss_ipm_woodbury}. The $t$-th iteration takes time
\begin{align*}
O(\Tmat(m, n^2, n r_t) + \Tmat(m, m, n r_t) + (n r_t)^{\omega} + m^2 + m n^2).
\end{align*}
\end{lemma}
\begin{proof}

We compute the cost of each line of Algorithm~\ref{alg:ss_ipm_woodbury} in the $t$-th iteration.

{\bf Line~\ref{line:gradient_sdp_woodbury} of Algorithm~\ref{alg:ss_ipm_woodbury}: gradient computation, $O(m n^2)$ time.}

This step computes $g_{\eta^{\new}}(y)_j \leftarrow  \eta^{\new} \cdot b_j -   \tr [ S^{-1} \cdot A_j ]$.

Computing one $\tr [ S^{-1} \cdot A_j ]$ takes $n^2$ time, and computing all traces for $j \in [m]$ takes $m n^2$ time.

{\bf Line~\ref{line:delta_y_sdp_woodbury} of Algorithm~\ref{alg:ss_ipm_woodbury}: $\delta_y$ computation, $O(m^2)$ time.} 

This step computes $\delta_y \leftarrow - (\wt{H}^{\new})^{-1} \cdot g_{\eta^{\new}}(y)$.

Computing the matrix vector multiplication of $\wt{H}(y)^{-1} \in \R^{m \times m}$ with $g_{\eta^{\new}}(y) \in \R^m$ takes $O(m^2)$ time.

{\bf Line~\ref{line:S_new_sdp_woodbury} of Algorithm~\ref{alg:ss_ipm_woodbury}: $S^{\new}$ computation, $O(m n^2)$ time.}

This step computes $S^{\new} \leftarrow \sum_{i \in [m]} (y^{\new})_i A_i - C$.

Brute-forcely adding all $(y^{\new})_i A_i \in \R^{n \times n}$ takes $mn^2$ time.

{\bf Line~\ref{line:low_rank_slack_update}-\ref{line:slack_update} of Algorithm~\ref{alg:ss_ipm_woodbury}: $\wt{S}^{\new}$ computation, $O(n^{\omega})$ time.}

Line~\ref{line:low_rank_slack_update} makes a call to procedure \textsc{LowRankSlackUpdate}, and this takes $O(\omega)$ time by Lemma~\ref{lem:cost_low_rank_slack}.

Line~\ref{line:slack_update} computes $\wt{S}^{\new} \leftarrow \wt{S} + V_1 V_2^{\top}$, which takes $O(\Tmat(n, r_t, n)) \leq O(n^{\omega})$ time.

{\bf Line~\ref{line:V_3}-\ref{line:hessian_sdp_woodbury} of Algorithm~\ref{alg:ss_ipm_woodbury}: Hessian inverse computation, $O(\Tmat(m, n^2, n r_t) + \Tmat(m,m,n r_t) + (n r_t)^{\omega})$ time.} 

Computing $V_3, V_4 \in \R^{n \times r_t}$ on Line~\ref{line:V_3} and \ref{line:V_4} takes $O(n^{\omega})$ time.

Computing $Y_1, Y_2 \in \R^{n^2 \times (2 n r_t + r_t^2)}$ on Line~\ref{line:Y_1} and \ref{line:Y_2} takes the same time as the output size: $O(n^3 \cdot r_t)$.

We compute $G^{\new} \leftarrow G - G \cdot \mathsf{A} Y_1 \cdot (I + Y_2^{\top} \mathsf{A}^{\top} \mathsf{A} Y_1)^{-1} \cdot Y_2^{\top} \mathsf{A}^{\top} \cdot G$ on Line~\ref{line:hessian_sdp_woodbury} in the following order:
\begin{enumerate}
    \item Compute $\mathsf{A} Y_1, \mathsf{A} Y_2 \in \R^{m \times (2n r_t + r_t^2)}$ in $O(\Tmat(m, n^2, n r_t))$ time since $\mathsf{A} \in \R^{m \times n^2}$ and $Y_1,Y_2 \in \R^{n^2 \times (2n r_t + r_t^2)}$.
    \item Compute $I + (Y_2^{\top} \mathsf{A}^{\top}) \cdot (\mathsf{A} Y_1) \in \R^{(2n r_t + r_t^2)\times (2n r_t + r_t^2)}$ in $O(\Tmat(n r_t, m, n r_t))$ time. Then compute the inverse $(I + Y_2^{\top} \mathsf{A}^{\top} \cdot \mathsf{A} Y_1)^{-1} \in \R^{(2n r_t + r_t^2)\times (2n r_t + r_t^2)}$ in $O((n r_t)^{\omega})$ time.
    \item Compute $G \cdot \mathsf{A} Y_1 \in \R^{m \times (2n r_t + r_t^2)}$ in $O(\Tmat(m, m, n r_t))$ time since $G \in \R^{m \times m}$ and $\mathsf{A} Y_1 \in \R^{m \times (2n r_t + r_t^2)}$.
    \item Finally compute $G \mathsf{A} Y_1 \cdot (I + Y_2^{\top} \mathsf{A}^{\top} \cdot \mathsf{A} Y_1)^{-1} \cdot Y_2^{\top} \mathsf{A}^{\top} G \in \R^{m \times m}$ in $O(\Tmat(m, n r_t, m))$ time.
\end{enumerate}
Thus in total computing $G^{\new}$ takes $O(\Tmat(m, n^2, n r_t) + \Tmat(m, m, n r_t) + (n r_t)^{\omega})$ time.

{\bf Combined.} Combining the time of the four steps on Line~\ref{line:gradient_sdp_woodbury}, \ref{line:hessian_sdp_woodbury}, \ref{line:delta_y_sdp_woodbury}, and \ref{line:S_new_sdp_woodbury} of Algorithm~\ref{alg:ss_ipm_woodbury}, it is easy to see that the $t$-th iteration takes time
\begin{align*}
    &~ \text{Time~per~iteration} \\
    = & ~ \text{Line~\ref{line:gradient_sdp_woodbury}} + \text{Line~\ref{line:delta_y_sdp_woodbury}} + \text{Line~\ref{line:S_new_sdp_woodbury}} + \text{Line~\ref{line:low_rank_slack_update}-\ref{line:slack_update}} + \text{Line~\ref{line:V_3}-\ref{line:hessian_sdp_woodbury}} \\
    = & ~ \underbrace{ O(mn^2) }_{ \text{Line~\ref{line:gradient_sdp_woodbury}} } + \underbrace{O(m^2)}_{\text{Line~\ref{line:delta_y_sdp_woodbury}}} + \underbrace{ O(mn^2) }_{ \text{Line~\ref{line:S_new_sdp_woodbury}} } + \underbrace{ O( n^{\omega} ) }_{\text{Line~\ref{line:low_rank_slack_update}-\ref{line:slack_update}}} + \underbrace{ \Tmat(m, n^2, n r_t) + \Tmat(m, m, n r_t) + (n r_t)^{\omega}}_{ \text{Line~\ref{line:V_3}-\ref{line:hessian_sdp_woodbury}} }\\
    = & ~ O\left(\Tmat(m, n^2, n r_t) + \Tmat(m, m, n r_t) + (n r_t)^{\omega} + m^2 + m n^2\right).
\end{align*}
Thus, we complete the proof.

\end{proof}

\begin{lemma}[Cost of \textsc{LowRankSlackUpdate} (Algorithm~\ref{alg:low_rank_slack_update})]\label{lem:cost_low_rank_slack}
A call to procedure \textsc{LowRankSlackUpdate} (Algorithm~\ref{alg:low_rank_slack_update}) takes $O(n^{\omega})$ time.
\end{lemma}
\begin{proof}
In procedure \textsc{LowRankSlackUpdate}, the most time-consuming step is to compute the spectral decomposition of $Z^{\midd} \in \R^{n \times n}$, and this takes $O(n^{\omega})$ time.
\end{proof}

\subsection{Property of low rank update}\label{sec:property_of_low_rank_update}

The goal of this section is to prove Theorem~\ref{thm:general_amortize_guarantee}.

We first make the following definitions.
\begin{definition}[Potential function]\label{def:potential}
Let $g \in \R_+^n$ be a non-increasing vector. For any matrix $Z \in \R^{n \times n}$, let $|\lambda(Z)|_{[i]}$ to denote the $i$-th largest absolute eigenvalue of $Z$. We define a potential function $\Phi_g : \R^{n \times n} \to \R_+$,
\begin{align*}
    \Phi_g(Z) := \sum_{i=1}^n g_i \cdot |\lambda(Z)|_{[i]}.
\end{align*}
\end{definition}

In the $t$-th iteration of Algorithm~\ref{alg:ss_ipm_woodbury}, let $S, \wt{S} \in \R^{n \times n}$ be the slack matrix and the approximate slack matrix in the beginning of the iteration, and let $S^{\new}, \wt{S}^{\new}$ be the updated matrices. We define the following matrices to capture their differences. 
\begin{definition}[Difference matrices]\label{def:diff_matrix}
Define the following matrices $Z, Z^{\midd}, Z^{\new} \in \R^{n \times n}$:
\begin{align*}
    Z := &~ S^{-1/2} \cdot \wt{S} \cdot S^{-1/2} - I, \\
    Z^{\midd} := &~ (S^{\new})^{-1/2} \cdot \wt{S} \cdot (S^{\new})^{-1/2} - I, \\
    Z^{\new} := &~ (S^{\new})^{-1/2} \cdot \wt{S}^{\new} \cdot (S^{\new})^{-1/2} - I.
\end{align*}
\end{definition}

\begin{assumption}[Closeness of $S^{\new}$ and $\wt{S}$ from $S$]\label{ass:close_S}
We make the following two assumptions about $S, \wt{S}, S^{\new} \in \R^{n \times n}$:
\begin{align*}
    1. &~ \|S^{-1/2} \cdot S^{\new} \cdot S^{-1/2} - I\|_F \leq 0.02, \\
    2. &~ \|S^{-1/2} \cdot \wt{S} \cdot S^{-1/2} - I\|_2 \leq 0.01.
\end{align*}
\end{assumption}

Next we present the main theorem of this section.
\begin{theorem}[General amortized guarantee]\label{thm:general_amortize_guarantee}
Let $T$ denote the total number of iterations in Algorithm~\ref{alg:ss_ipm_woodbury}. Let $r_t$ denote the rank of the update matrices $V_1, V_2 \in \R^{n \times r_t}$ generated by Algorithm~\ref{alg:low_rank_slack_update} in the $t$-th iteration.
The ranks $r_t$'s satisfy the following condition: for any vector $g \in \R_+^n$ which is non-increasing, we have
\begin{align*}
    \sum_{t=1}^T r_t \cdot g_{r_t} \leq O(T \cdot \|g\|_2 \cdot \log n).
\end{align*}
\end{theorem}
%Note that by choosing $g_i = 1/\sqrt{i}, \forall i \in [n]$ the above lemma recovers Theorem~6.1 of \cite{jklps20}. 
\begin{proof}
Part~1 of Assumption~\ref{ass:close_S} is proved in Theorem~\ref{thm:correctness_ss_ipm_woodbury}, and Part~2 of Assumption~\ref{ass:close_S} is proved in Lemma~\ref{lem:approx_slack}. Thus we can use Lemma~\ref{lem:S_move}.

Combining Lemma~\ref{lem:S_move} and Lemma~\ref{lem:wt_S_move}, we have
\begin{align*}
    \Phi_g(Z^{\new}) - \Phi_g(Z) 
    = & ~ (\Phi_g(Z^{\midd}) - \Phi_g(Z)) - (\Phi_g(Z^{\midd}) - \Phi_g(Z^{\new})) \\
    \leq & ~ \|g\|_2 - \frac{\epsilon_S}{10 \log n} \cdot r_t \cdot g_{r_t}.
\end{align*}
With an abuse of notation, we denote the matrix $Z$ in the $t$-th iteration as $Z^{(t)}$. Since in the beginning $\Phi_g(Z^{(0)}) = 0$ and $\Phi_g(Z^{(T)}) \geq 0$, we have
\begin{align*}
    0 \leq & ~ \Phi_g(Z^{(T)}) - \Phi_g(Z^{(0)}) \\
    \leq & ~ \sum_{t=1}^T (\Phi_g(Z^{(t)}) - \Phi_g(Z^{(t-1)})) \\
    \leq & ~ T\cdot \|g\|_2 - \frac{\epsilon_S}{10 \log n} \cdot \sum_{t=1}^T r_t \cdot g_{r_t}.
\end{align*}
This completes the proof.
\end{proof}

\subsubsection{\texorpdfstring{$S$}{} move}

\begin{lemma}[Eigenvalue change]
\label{lem:part_of_lem_62_jkl}
Let matrices $Z, Z^{\midd} \in \R^{n \times n}$ be defined as in Definition~\ref{def:diff_matrix}. Under Assumption~\ref{ass:close_S}, we have
\begin{align*}
    \sum_{i=1}^n (\lambda(Z)_{[i]} - \lambda(Z^{\midd})_{[i]})^2 \leq 10^{-3},
\end{align*}
where $\lambda(Z)_{[i]}$ denotes the $i$-th largest eigenvalue of $Z$.
\end{lemma}

\begin{proof}
We notice
\begin{align*}
    \|S^{1/2}(S^{\new})^{-1}S^{1/2}- I\|_F^2= \| \nu^{-1} - {\bf 1}_n \|_2^2 % \sum_{i=1}^n (\nu_i^{-1}-1)^2,
\end{align*}
where $\{\nu_i\}_{i\in[n]}$ are the eigenvalues of $S^{-1/2}S^{\new}S^{-1/2}$. By Assumption~\ref{ass:close_S}, we have
\begin{align*}
    \max_{i\in [n]} |v_i - 1| \leq 0.02,
\end{align*}
which implies that $ \min_{i \in [n]}v_i\geq 0.98$. % for every $i\in \{1,2,\cdots, n \}$. 

  Assumption~\ref{ass:close_S} also gives 
\begin{align*}
    %\sum_{i=1}^n (v_i - 1)^2 
   \| \nu - {\bf 1}_n \|_2^2 \leq 0.0004 .
\end{align*}
Then, we have
\begin{align}\label{eq:sum_vi_minus_1_square}
    %\sum_{i=1}^n (\nu_i^{-1}-1)^2 
    \| \nu^{-1} - {\bf 1}_n \|_2^2 \leq 5\times 10^{-4}.
\end{align}

Define $F := (S^\new)^{-1/2}S^{1/2}$ and $F = U D V^\top$ be its SVD decomposition. Let $Z' := V^\top Z V$. Notice that
\begin{align*}
    Z^\midd %= &~ (S^\new)^{-1/2}\wt{S} (S^\new)^{-1/2} - I \\
    =  F Z F^\top + F F^\top - I.
\end{align*}
Combining with Eq.~\eqref{eq:sum_vi_minus_1_square},
\begin{align*}
    \sum_{i=1}^n (\lambda(Z^\new)_{[i]} - \lambda( F Z F^\top)_{[i]})^2 \leq &~ \|Z^\midd - F Z F^\top \|_F^2 \\
    = &~ \|FF^\top - I\|_F^2\\
    = &~ \| \nu^{-1} - {\bf 1}_n \|_2^2 \\% \sum_{i=1}^n (\nu_i^{-1}-1)^2\\
    \leq &~ 5\times 10^{-4}.
\end{align*}
Since $\|D^2 - I \|_F^2 = \|FF^\top - I \|_F^2 \leq 5 \cdot 10^{-4}$ and $\|Z'\|_2 = \|Z\|_2 \leq 0.01$, Fact~\ref{fac:sigma_z_sigma_z_norm} gives
\begin{align*}
    \sum_{i=1}^n (\lambda(Z)_{[i]} - \lambda(F Z F^\top)_{[i]})^2 = &~ \sum_{i=1}^n (\lambda(Z)_{[i]} - \lambda(D Z' D)_{[i]})^2\\ \leq &~ \| D Z' D - Z'\|_F^2\\
    \leq &~ 10^{-5}.
\end{align*}
Combining the above inequalities, we have
\begin{align*}
    \sum_{i=1}^n (\lambda(Z)_{[i]} - \lambda(Z^{\midd})_{[i]})^2 \leq 10^{-3}.
\end{align*}
\end{proof}

The following lemma upper bounds the increase in potential when $S$ changes to $S^{\new}$. 
\begin{lemma}[$S$ move]\label{lem:S_move}
Consider the $t$-th iteration. Let matrices $Z,Z^{\midd} \in \R^{n \times n}$ be defined as in Definition~\ref{def:diff_matrix}. Let $g \in \R_+^n$ be a non-increasing vector, and let $\Phi_g: \R^{n \times n} \to \R_+$ be defined as in Definition~\ref{def:potential}.

Under Assumption~\ref{ass:close_S}, we have
\begin{align*}
    \Phi_g(Z^{\midd}) - \Phi_g(Z) \leq \|g\|_2.
\end{align*}
\end{lemma}
\begin{proof}
Let $\pi: [n] \to [n]$ be a sorting permutation such that $|\lambda(Z^{\midd})_{\pi(1)}| \geq |\lambda(Z^{\midd})_{\pi(2)}| \geq \cdots \geq |\lambda(Z^{\midd})_{\pi(n)}|$, i.e., $|\lambda(Z^{\midd})_{\pi(i)}| = |\lambda(Z^{\midd})|_{[i]}$. Then we have
\begin{align*}
    \Phi_g(Z^{\midd}) = &~ \sum_{i=1}^n g_i \cdot |\lambda(Z^{\midd})_{\pi(i)}| \\
    \leq &~ \sum_{i=1}^n g_i \cdot |\lambda(Z)_{\pi(i)}| + \sum_{i=1}^n g_i \cdot |\lambda(Z^{\midd})_{\pi(i)} - \lambda(Z)_{\pi(i)}| \\
    \leq &~ \sum_{i\in[n]} g_i \cdot |\lambda(Z)|_{[i]} + \sum_{i\in[n]} g_i \cdot |\lambda(Z^{\midd})_{\pi(i)} - \lambda(Z)_{\pi(i)}| \\
    \leq &~ \sum_{i\in[n]} g_i \cdot |\lambda(Z)|_{[i]} + \Big( \sum_{i\in [n]} g_i^2 \Big)^{1/2} \cdot \Big( \sum_{i\in[n]} |\lambda(Z^{\midd})_{\pi(i)} - \lambda(Z)_{\pi(i)}|^2 \Big)^{1/2} \\
    \leq &~ \Phi_g(Z) + \|g\|_2
\end{align*}
where the first step follows from the definition of $\Phi_g$ (Definition~\ref{def:potential}) and $\pi$, the second step follows from triangle inequality, the third step follows from $g \in \R_+^n$ is non-increasing and $|\lambda(Z)|_{[i]}$ is the $i$-th largest absolute eigenvalue of $Z$, the fourth step follows from Cauchy-Schwarz inequality, the last step follows from definition of $\Phi_g$ and Lemma~\ref{lem:part_of_lem_62_jkl}.

Thus we have $\Phi_g(Z^{\midd}) - \Phi_g(Z) \leq \|g\|_2$.
\end{proof}

\subsubsection{\texorpdfstring{$\wt{S}$}{S~} move}
The following lemma lower bounds the decrease in potential when $\wt{S}$ is updated to $\wt{S}^{\new}$. 
\begin{lemma}[$\wt{S}$ move]\label{lem:wt_S_move}
Consider the $t$-th iteration. Let matrices $Z^{\midd}, Z^{\new} \in \R^{n \times n}$ be defined as in Definition~\ref{def:diff_matrix}. Let $g \in \R_+^n$ be a non-increasing vector, and let $\Phi_g: \R^{n \times n} \to \R_+$ be defined as in Definition~\ref{def:potential}.

Let $r_t$ denote the rank of the update matrices $V_1, V_2 \in \R^{n \times r_t}$ generated by Algorithm~\ref{alg:low_rank_slack_update} in the $t$-th iteration. We have
\begin{align*}
    \Phi_g(Z^{\midd}) - \Phi_g(Z^{\new}) \geq \frac{\epsilon_S}{10 \log n} \cdot r_t \cdot g_{r_t}.
\end{align*}
\end{lemma}
\begin{proof}
Note that $r_t = 2 r$, where $r$ is the variable of Algorithm~\ref{alg:low_rank_slack_update}. We define $\lambda \in \R^n$ and $U \in \R^{n \times n}$ (Line~\ref{line:Z_mid_decomposition}), $\pi: [n] \to [n]$ (Line~\ref{line:pi}) and $\lambda^{\new} \in \R^n$ (Line~\ref{line:eigenvalue_update}) to be the same as Algorithm~\ref{alg:low_rank_slack_update}. We extend the definition and let $\lambda_{\pi(i)} = 0$ for $i > n$. 
We have

\begin{align*}
    Z^{\new} = &~ (S^{\new})^{-1/2} \cdot \wt{S}^{\new} \cdot (S^{\new})^{-1/2} - I \\
    = &~ (S^{\new})^{-1/2} \cdot \Big( \wt{S} + (S^{\new})^{1/2} \cdot U \cdot \diag (\lambda^{\new} - \lambda) \cdot U^\top \cdot (S^{\new})^{1/2} \Big) \cdot (S^{\new})^{-1/2} - I \\
    %= &~ (S^{\new})^{-1/2} \cdot (\wt{S} + V_1 V_2^{\top}) \cdot (S^{\new})^{-1/2} - I \\
    = &~ Z^{\midd} + U \cdot \diag (\lambda^{\new} - \lambda) \cdot U^\top \\
    = &~ U \cdot \diag(\lambda^{\new}) \cdot U^{\top}
\end{align*}
where the first step follows from definition of $Z^{\new}$ (Definition~\ref{def:diff_matrix}), the second step follows from closed-form formula of $\wt{S}^{\new}$ (Lemma~\ref{lem:close_form_slack}), the third step follows from definition of $Z^{\midd}$ (Definition~\ref{def:diff_matrix}), the fourth step follows from $Z^{\midd} = U \cdot \diag(\lambda) \cdot U^{\top}$ (Line~\ref{line:Z_mid_decomposition} of Algorithm~\ref{alg:low_rank_slack_update}).

Thus from the definition of $\Phi_g$ (Definition~\ref{def:potential}), we have
\begin{align}\label{eq:phi_Z_mid_phi_Z_new}
    \Phi_g(Z^{\midd}) = \sum_{i=1}^n g_i \cdot |\lambda_{\pi(i)}|,~~~~~
    \Phi_g(Z^{\new}) = \sum_{i=1}^{n} g_i \cdot |\lambda_{\pi(2r+i)}|.
\end{align}

We consider two different cases of the outcome of the while-loop on Line~\ref{line:while_slack_update} of Algorithm~\ref{alg:low_rank_slack_update}.

{\bf Case 1. No $i \leq n/2$ satisfies both $|\lambda_{\pi(2 i)}| \leq \epsilon_S$ and $|\lambda_{\pi(2 i)}| \leq (1 - 1/\log n) |\lambda_{\pi(i)}|$.}

In this case, the while-loop exits with $r=n/2$, and hence $r_t = 2r = n$. Thus using Eq.~\eqref{eq:phi_Z_mid_phi_Z_new} we have $\Phi_g(Z^{\new}) = 0$. We consider two sub-cases.
\begin{itemize}
    \item {\bf Case 1(a). For all $i \in [n]$, $|\lambda_{\pi(i)}| > \epsilon_S$.}
    
    In this case we have 
    \begin{align*}
        \Phi_g(Z^{\midd}) - \Phi_g(Z^{\new}) 
        = & ~ \sum_{i=1}^n g_i \cdot |\lambda_{\pi(i)}| - 0 \\
        \geq & ~ n \cdot g_n \cdot \epsilon_S \\
        = & ~ \epsilon_S \cdot r_t \cdot g_{r_t}.
    \end{align*}
    \item {\bf Case 1(b). There exists a minimum $i^* \leq n/2$ such that $|\lambda_{\pi(2 i)}| \leq \epsilon_S$ for all $i \geq i^*$.}
    
    The condition of Case 1 and Case 1(b) means that for all $i \geq i^*$, we must have $|\lambda_{\pi(2 i)}| > (1 - 1/\log n) |\lambda_{\pi(i)}|$. And since $i^*$ is the minimum index such that $|\lambda_{\pi(2 i^*)}| \leq \epsilon_S$, we have $|\lambda_{\pi(i^*)}| > \epsilon_S$. Thus we have
    \begin{align*}
        |\lambda_{\pi(n)}| > (1 - 1/\log n) |\lambda_{\pi(n/2)}| > \cdots > (1 - 1/\log n)^{\log(n/i^*)} |\lambda_{\pi(i^*)}| \geq \frac{1}{e} \cdot |\lambda_{\pi(i^*)}| \geq \frac{\epsilon_S}{e}.
    \end{align*}
    So we have 
    \begin{align*}
    \Phi_g(Z^{\midd}) - \Phi_g(Z^{\new}) 
    = & ~ \sum_{i=1}^n g_i \cdot |\lambda_{\pi(i)}| - 0 \\
    \geq & ~ n \cdot g_n \cdot |\lambda_{\pi(n)}| \\
    = & ~ \frac{\epsilon_S}{e} \cdot r_t \cdot g_{r_t}.
    \end{align*}
\end{itemize}

{\bf Case 2. There exists a minimum $r \leq n/2$ such that both $|\lambda_{\pi(2 r)}| \leq \epsilon_S$ and $|\lambda_{\pi(2 r)}| \leq (1 - 1/\log n) |\lambda_{\pi(r)}|$ are satisfied.}

This $r$ is the outcome of the while-loop in Algorithm~\ref{alg:low_rank_slack_update}. Next we prove $|\lambda_{\pi(r)}| \geq \frac{\epsilon_S}{e}$.

Let $i^* \leq n/2$ be the minimum index such that $|\lambda_{\pi(2 i)}| \leq \epsilon_S$ for all $i \geq i^*$. Note that this implies $|\lambda_{\pi(i^*)}| > \epsilon_S$. We have $r \geq i^*$ since $|\lambda_{\pi(2r)}| \leq \epsilon_S$ and $i^*$ is the minimum such index. Hence $|\lambda_{\pi(2 i)}| > (1 - 1/\log n) |\lambda_{\pi(i)}|$ for all $i \in [i^*, r)$.

If $r/2 \leq i^*$, we have $|\lambda_{\pi(r)}| = |\lambda_{\pi(2(r/2))}| \geq \epsilon_S$. If $r/2 \geq i^*$, we have
\begin{align*}
    |\lambda_{\pi(r)}| > (1 - 1/\log n) |\lambda_{\pi(r/2)}| > \cdots > (1 - 1/\log n)^{\log(r/i^*)} |\lambda_{\pi(i^*)}| \geq \frac{1}{e} \cdot |\lambda_{\pi(i^*)}| \geq \frac{\epsilon_S}{e}.
\end{align*} 
Thus we have $|\lambda_{\pi(r)}| \geq \frac{\epsilon_S}{e}$ in either cases.

We have
\begin{align}\label{eq:lambda_pi_i_bound}
    1.~ \forall i \leq r,~ |\lambda_{\pi(i)}| \geq |\lambda_{\pi(r)}| \geq \frac{\epsilon_S}{e}, ~~~~~~
    2.~ \forall i \geq 2 r,~ |\lambda_{\pi(i)}| \leq |\lambda_{\pi(2 r)}| \leq (1 - 1 / \log n) |\lambda_{\pi(r)}|.
\end{align}
Thus we can bound the potential decrease as follows:
\begin{align*}
    \Phi_g(Z^{\midd}) - \Phi_g(Z^{\new})
    \geq &~ \sum_{i=1}^{r} g_i \cdot (|\lambda_{\pi(i)}| - |\lambda_{\pi(i+2r)}|) \\
    \geq &~ \sum_{i=1}^{r} g_i \cdot (|\lambda_{\pi(r)}| - |\lambda_{\pi(2r)}|) \\
    \geq &~ \sum_{i=1}^{r} g_i \cdot \frac{1}{\log n} |\lambda_{\pi(r)}| \\
    \geq &~ r_t \cdot g_{r_t} \cdot \frac{\epsilon_S}{2 e \log n},
\end{align*}
where the first step follows from Eq.~\eqref{eq:phi_Z_mid_phi_Z_new}, the second step follows from $|\lambda_{\pi(i)}|$ is decreasing, the third follows from Part~2 of Eq.~\eqref{eq:lambda_pi_i_bound}, the fourth step follows from Part~1 of Eq.~\eqref{eq:lambda_pi_i_bound}, $r_t = 2r$, and that $g$ is non-increasing.
\end{proof}

\subsection{Amortized analysis}\label{sec:amortized_analysis}
The goal of this section is to prove Lemma~\ref{lem:amortize_hessian} using Lemma~\ref{lem:general_amortize_tool}. 

\begin{lemma}[Amortization of Hessian computation]\label{lem:amortize_hessian}
Let $T$ denote the total number of iterations in Algorithm~\ref{alg:ss_ipm_woodbury}. For $t \in [T]$, the cost of Hessian computation in the $t$-th iteration can be amortized as follows:
\begin{align*}
    &~ \sum_{t=1}^T (\Tmat(m, n^2, n r_t) + \Tmat(m, m, n r_t) + (n r_t)^{\omega}) \\
    \leq &~ T \cdot O^\ast(m^2 + n^4 + n^{2 \omega - 1/2} + m^{2 - \frac{\alpha(\omega-2)}{1-\alpha}} \cdot n^{\frac{2(\omega - 2)}{1 - \alpha}-1/2} ).
\end{align*}
\end{lemma}
\begin{proof}
We have $(n r_t)^{\omega} \leq \Tmat(n^2, n^2, n r_t)$, and $\Tmat(m, n^2, n r_t) \leq \Tmat(m, m, n r_t) + \Tmat(n^2, n^2, n r_t)$, thus we can bound the cost of hessian computation as
\begin{align}\label{eq:hessian_cost_relax}
    \Tmat(m, n^2, n r_t) + \Tmat(m, m, n r_t) + (n r_t)^{\omega} \leq \Tmat(m, m, n r_t) + \Tmat(n^2, n^2, n r_t).
\end{align}

From Theorem~\ref{thm:general_amortize_guarantee} we know that for any vector $g \in \R_+^n$ which is non-increasing, we have $\sum_{t=1}^T r_t \cdot g_{r_t} \leq O(T \cdot \|g\|_2 \cdot \log n)$, so we can use Lemma~\ref{lem:general_amortize_tool}:
\begin{align*}
    \sum_{t=1}^T \Tmat(m, m, n r_t) \leq &~ O^\ast(m^2 + m^{2 - \frac{\alpha(\omega-2)}{1-\alpha}} \cdot n^{\frac{2(\omega - 2)}{1 - \alpha}-1/2}) \text{~when~}m \geq n^2,\\
    \sum_{t=1}^T \Tmat(n^2, n^2, n r_t) \leq &~ O^\ast(n^4 + n^{4 - \frac{2\alpha(\omega-2)}{1-\alpha}} \cdot n^{\frac{2(\omega - 2)}{1 - \alpha}-1/2} ) = O^\ast(n^4 + n^{2 \omega - 1/2} ).
\end{align*}
Combining these two inequalities and Eq.~\eqref{eq:hessian_cost_relax} completes the proof.
\end{proof}

\begin{lemma}[Helpful lemma for amortization of Hessian computation]\label{lem:general_amortize_tool}
Let $T$ denote the total number of iterations. Let $r_t \in [n]$ be the rank for the $t$-th iteration for $t \in [T]$. Assume $r_t$ satisfies the following condition:
for any vector $g \in \R_+^n$ which is non-increasing, we have 
\begin{align*}
\sum_{t=1}^T r_t \cdot g_{r_t} \leq O(T \cdot \|g\|_2).
\end{align*}

If the cost in the $t$-th iteration is $O(\Tmat(d, d, n r_t))$ where $d =\Omega( n^2)$ is an integer, then the amortized cost per iteration is 
\begin{align*}
O^\ast\big( d^2 + d^{2 - \frac{\alpha(\omega-2)}{1-\alpha}} \cdot n^{\frac{2(\omega - 2)}{1 - \alpha}-1/2}\big).
\end{align*}
\end{lemma}
\begin{proof}
Let $\omega$ be the matrix multiplication exponent, let $\alpha$ be the dual matrix multiplication exponent, and let $\beta = \omega(2)$ (see Section~\ref{sec:fast_matrix_multiplication} for more details). Since $d \geq n^2 \geq n r_t$, we have
\begin{align}\label{eq:amortization_eq1}
    \Tmat(d, d, n r_t) \leq & ~ d^2 + (n r_t)^{\frac{\omega-2}{1-\alpha}} \cdot d^{2 - \frac{\alpha(\omega-2)}{1-\alpha}} \notag \\
    = & ~ d^2 + d^{2 - \frac{\alpha(\omega-2)}{1-\alpha}} \cdot n^{\frac{\omega-2}{1-\alpha}} \cdot r_t^{\frac{\omega-2}{1-\alpha}},
\end{align}
where we use Fact~\ref{fact:approx_Tmat_r_leq_n} in the first step.

Define $g \in \R_+^n$ such that $\forall r \in [n]$, $g_r = r^{\frac{\omega-2}{1-\alpha}-1}$. We observe that $g$ is a non-increasing vector because $\frac{\omega-2}{1-\alpha} - 1 \leq 0$ (Fact~\ref{fact:bound_omega_alpha}). Then using the condition in the lemma statement, we have
\begin{align}\label{eq:amortization_eq2}
    \sum_{t=1}^T r_t^{\frac{\omega-2}{1-\alpha}} 
    = & ~ \sum_{t=1}^T r_t \cdot g_{r_t} \notag \\
    \leq & ~ T \cdot \|g\|_2  \notag \\
    \leq & ~ T \cdot \Big(\int_{x=1}^n x^{\frac{2(\omega-2)}{1-\alpha}-2} \mathrm{d}x \Big)^{1/2} \notag \\
    = & ~ T \cdot O(n^{\frac{(\omega-2)}{1-\alpha}-1/2}),
\end{align}
where the first step follows from the definition that $g_r = r^{\frac{\omega-2}{1-\alpha}-1}$, $\forall r \in [n]$, the second step follows from the assumption $\sum_{t=1}^T r_t \cdot g_{r_t} \leq T \cdot \|g\|_2$ in the lemma statement,  
the third step follows from upper bounding the $\ell_2$ norm $\|g\|_2^2 = \sum_{r=1}^n g_r^2 = \sum_{r=1}^n r^{\frac{2(\omega-2)}{1-\alpha}-2} \leq \int_{x=1}^n x^{\frac{2(\omega-2)}{1-\alpha}-2}$ %$g_r = r^{\frac{\omega-2}{1-\alpha}-1}$ 
and the last step follows from computing the integral $\int_{x=1}^n x^{\frac{2(\omega-2)}{1-\alpha}-2} = c \cdot x^{\frac{2(\omega-2)}{1-\alpha}-1}\big|_{1}^n = O(n^{\frac{2(\omega-2)}{1-\alpha}-1})$ where $c:= 1/(\frac{2(\omega-2)}{1-\alpha}-1)$. 

Thus we have
\begin{align*}
    \sum_{t=1}^T \Tmat(d, d, n r_t)
    \leq & ~ \sum_{t=1}^T \Big( d^2 + d^{2 - \frac{\alpha(\omega-2)}{1-\alpha}} \cdot n^{\frac{\omega-2}{1-\alpha}} \cdot r_t^{\frac{\omega-2}{1-\alpha}} \Big) \\
    = &~ T \cdot d^2 + d^{2 - \frac{\alpha(\omega-2)}{1-\alpha}} \cdot n^{\frac{\omega-2}{1-\alpha}} \cdot \sum_{t=1}^T r_t^{\frac{\omega-2}{1-\alpha}} \\ 
    \leq &~ T \cdot d^2 + d^{2 - \frac{\alpha(\omega-2)}{1-\alpha}} \cdot n^{\frac{\omega-2}{1-\alpha}} \cdot T \cdot O(n^{\frac{(\omega-2)}{1-\alpha}-1/2}) \\
    = &~ T \cdot O(d^2 + d^{2 - \frac{\alpha(\omega-2)}{1-\alpha}} \cdot n^{\frac{2(\omega - 2)}{1 - \alpha}-1/2}).
\end{align*}
where the first step follows from Eq.~\eqref{eq:amortization_eq1}, the second step follows from moving summation inside, the third step follows from Eq.~\eqref{eq:amortization_eq2}, and the last step follows from adding the terms together.

Thus, we complete the proof.
\end{proof}

\section{The Robust Interior Point Method Framework For SDP}
\label{sec:correct_previous}

%\subsection{Approximate Hessian}

%\subsection{Correctness of approximate barrier method}\label{sec:correct_central_path}

One of the contributions in this work is a more robust interior point method framework that allows errors in the Hessian matrices, the gradient vectors, and the Newton steps. We will first introduce the necessary backgrounds and definitions in Section~\ref{sec:robust_ipm_def}. Then we will perform the one step error analysis based on Newton decrements in Section~\ref{sec:our_roubst_ns}. We also list the corresponding error analysis in previous framework for comparison. In Section~\ref{sec:robust_ipm_eta_move}-\ref{sec:robust_ipm_integral}, we prove several supporting Lemmata that are used in the proof of Section~\ref{sec:our_roubst_ns}. In Section~\ref{sec:robust_ipm_duality}, we include several classical results that bound the duality gap by the Newton decrements and provide the proofs. Finally, we state the main result in Section~\ref{sec:robust_ipm_main}.

\subsection{Definitions}\label{sec:robust_ipm_def}
%{A New and Robust Invariance of Newton Step}

\iffalse
The following lemma is standard in interior point method (see e.g. Section 2.4 of \cite{r01}). In recent work \cite{jklps20}, they used the following Lemma as a black-box, see Lemma 5.5 in \cite{jklps20}.

\begin{lemma}[Invariance of Newton step]\label{lem:invariant_newton}
Given any parameters $1 \leq \alpha_H \leq 1.03$ and $0 < \epsilon_N  \leq 1/10$. Suppose that there is 
\begin{itemize}
    \item Condition 0. a feasible solution $y \in \R^m$ satisfies $\g_\eta(y)^\top (\nabla^2 \phi(y))^{-1} \g_\eta(y) \leq \epsilon_N^2$,
    \item Condition 1. a positive definite matrix $\wt{H} \in \mathbb{S}^{n \times n}_{>0}$ satisfies $\alpha_H^{-1} \nabla^2 \phi(y) \preceq \wt{H} \preceq \alpha_H  \nabla^2 \phi(y)$.
\end{itemize}
Then $\eta^{\new} = \eta (1 + \frac{\epsilon_N}{20 \sqrt{n}} ) $ and $y^{\new} = y - \wt{H}^{-1} \g(y,{\eta^{\new} })$ satisfy 
\begin{align*}
\g_{\eta^{\new}} (y^{\new})^\top (\nabla^2 \phi(y^{\new}))^{-1} \g_{\eta^{\new}}(y^{\new}) \leq \epsilon_N^2.
\end{align*}
\end{lemma}
\fi

%\newpage
We start with some definitions.

\begin{definition}\label{def:S_H}
Let $C \in \R^{n \times n}$. Let $\mathsf{A} \in \R^{m \times n^2}$ denote the matrix where the $i$-th row matrix $A_i \in \R^{n \times n}$. Consider a barrier function $\phi: \R^m \mapsto \R$ defined on the dual space $\{y \in \R^m: \sum_{i=1}^m y_i A_i - C \succeq 0 \}$.

We define function $S: \R^m \rightarrow \R^{n \times n}$ such that
\begin{align*}
    S(y) = \sum_{i=1}^m y_i \cdot A_i - C.
\end{align*}
We define function $\nabla^2 \phi : \R^m \rightarrow  \R^{m \times m}$ that maps the dual variable to the Hessian matrix of barrier function $\phi$. Notice $\nabla^2 \phi = \nabla^2 f_\eta$ for the regularized objective $f_\eta(y)$ in Eq.~\eqref{eq:augmented_obj}, since $f_\eta$ adds $\phi$ with a linear function of $y$. In particular, for logarithmic barrier:
\begin{align*}
    \nabla^2 \phi (y)= \mathsf{A} \cdot ({S}(y)^{-1} \otimes {S}(y)^{-1} ) \cdot \mathsf{A}^{\top}.
\end{align*}
We abuse the notation of $\nabla^2 \phi$ and also write $\nabla^2 \phi: \R^{n \times n} \rightarrow \R^{m \times m}$ as a function of the slack matrix. In particular, for logarithmic barrier:
\begin{align*}
    \nabla^2 \phi(S) = \mathsf{A} \cdot ({S}^{-1} \otimes {S}^{-1} ) \cdot \mathsf{A}^{\top}.
\end{align*}
We define $\g: \R^m \times \R \rightarrow \R^m$ that maps the dual variable $y$ and the learning rate $\eta$ to the gradient of the regularized objective $f_\eta(y)$ in Eq.~\eqref{eq:augmented_obj}, i.e. $\g(y,\eta) = \eta \cdot b - \nabla \phi(y)$. In particular, for logarithmic barrier:
\begin{align*}
    \g(y,\eta) = \eta \cdot b - \mathsf{A} \cdot \vect(S(y)^{-1}).
\end{align*}
We abuse the notation of $\g$ and also write $\g(S,\eta)$ as a function of slack variable $S$ and learning rate $\eta$. For example, in logarithmic barrier $\g(S,\eta)$ is given by
\begin{align*}
    \g(S,\eta) = \eta \cdot b - \mathsf{A} \cdot \vect( S^{-1} ).
\end{align*}
We define function $n : \R^m \times \R \rightarrow \R^m$ as the Newton step taken at $y$ with learning rate $\eta$:
\begin{align*}
n(y,\eta) = (\nabla^2 \phi(y))^{-1} \cdot \g(y,\eta).
\end{align*}
\end{definition}

\begin{definition}[Local norm]
We will frequently use local inner product defined w.r.t. function $\phi$ as follow: for $u,v \in \R^n$, 
$\langle u, v\rangle_x := \langle u, \nabla^2 \phi(x) v \rangle$.
It induces a local norm defined by $\|u\|_x := \langle u,u\rangle_x =  \langle u, \nabla^2 \phi(x) u \rangle$.
It also induces an operator norm for matrices defined by $\|A\|_x := \sup_{z} \frac{\| A z\|_x}{\|z\|_x}$.

In the Hilbert space equipped with local norm $\langle \cdot,  \cdot \rangle_x$, the gradient at $z$ is denoted by $\nabla \phi_x(z)$ and satisfies  $\nabla \phi_x(z)= (\nabla^2 \phi(x))^{-1} \nabla \phi(z)$. The Hessian at $z$ is denoted by $\nabla^2 \phi_x(z)$ and satisfies  $\nabla^2 \phi_x(z) = (\nabla^2 \phi(x))^{-1}\nabla^2 \phi(z)$. Similarly, $\g_x(z,\eta) = (\nabla^2 \phi(x))^{-1} \g(z,\eta)$.% $\nabla f_x(z) = (\nabla^2 f(x))^{-1} \nabla f(z)$.
\end{definition}

\begin{definition}[Self-concordant functional and barrier]\label{def:self_concordant_complexity}
Given a function with domain $D_f$. Let $B_x(y,r)$ denote the open ball of radius $r$ centered at $y$, where radius is measured with respect to the local norm $\|\cdot \|_x = \|\cdot \|_{\nabla^2 f(x)}$ defined w.r.t. $f$. 
The functional $f$ is called self-concordant if for all $x \in D_f$ we have $B_x(x,1) \subset D_f$, and if whenever $y \in B_x(x,1)$ the following holds for all $v \neq 0$
\begin{align*}
    1-\|y-x\|_x \leq \frac{\|v\|_y}{\|v\|_x} \leq \frac{1}{1-\|y-x\|_x}.
\end{align*}

A self-concordant functional $f$ is called a self-concordant barrier or barrier functional if 
\begin{align*}
    \theta_f := \sup_{x \in D_f} \|\nabla f_x(x)\|_x^2 < \infty.
\end{align*}
$\theta_f$ is referred to as the complexity of self-concordant barrier $f$.
\end{definition}

\begin{remark}
The complexity of logarithmic barrier $\phi_{\log}(y) = - \log \det ( S (y) ) $ is $n$ (\cite{nn94}). The complexity of Hybrid barrier $225(n/m)^{1/2} \cdot \left(\phi_{\mathrm{vol}}(y) + \phi_{\log}(y) \cdot (m-1)/(n-1)\right)$ is $\sqrt{mn}$ (\cite{a00}). For any barrier function $\phi$ with complexity $\theta$, since $f_\eta$ adds $\phi$ with a linear function, $f_\eta$ is also a self-concordant function with complexity $\theta$.
\end{remark}

We will use the following properties of self-concordance functions, from \cite{r01}.
\begin{theorem}[Self-concordant function, \cite{r01}]\label{thm:r01}
Given a self-concordant function $f: D_f\rightarrow \R$. For any two points $a , b \in D_f$, if $\|a - b\|_{\nabla^2 f(a)}\leq 1/4$, then we have
\begin{align}\label{eq:sc_1}
    \|(\nabla^2 f(a))^{-1} \nabla^2 f(b)\|_{\nabla^2 f(a)},~\|(\nabla^2 f(b))^{-1} \nabla^2 f(a) \|_{\nabla^2 f(a)} \leq \frac{1}{(1-\|b-a\|_{\nabla^2 f(a)})^2}.
\end{align}
Further,
\begin{align}\label{eq:sc_2}
    \|I-(\nabla^2 f(a))^{-1} \nabla^2 f(b)\|_{\nabla^2 f(a)},\|I-(\nabla^2 f(b))^{-1} \nabla^2 f(a) \|_{\nabla^2 f(a)} \leq \frac{1}{(1-\|b-a\|_{\nabla^2 f(a)})^2} - 1.
\end{align}
\end{theorem}

\begin{theorem}[Proposition 2.2.8 in \cite{r01}]\label{thm:r01_2.2.8}
Given a self-concordant function $f: D_f\rightarrow \R$. Define $n(a) := - (\nabla^2 f(a))^{-1} \nabla f(a)$.  If $\|n(a)\|_{\nabla^2 f(a)} \leq 1/4$ then $f$ has a minimizer $z$ and
\begin{align*}
    \|z - a^\new \|_{\nabla^2 f(a)} \leq \frac{3\|n(a)\|_{\nabla^2 f(a)}}{(1-\|n(a)\|_{\nabla^2 f(x)})^3}
\end{align*}
where $a^\new = a + n(a)$ and 
\begin{align*}
    \|z - a \|_{\nabla^2 f(a)} \leq \|n(a)\|_{\nabla^2 f(a)} + \frac{3\|n(a)\|_{\nabla^2 f(a)}}{(1-\|n(a)\|_{\nabla^2 f(a)})^3}.
\end{align*}
\end{theorem}

\begin{theorem}[Theorem 2.3.4 in \cite{r01}]\label{thm:r01_2.3.4}
Assume $f$ is a self-concordance barrier with domain $D_f$. Let $\theta_f$ denote its complexity. If $x,y \in D_f$ then
\begin{align*}
    \langle \nabla f(x),y-x \rangle \leq \theta_f.
\end{align*}
\end{theorem}

\subsection{One step error analysis}\label{sec:our_roubst_ns}

Classical interior point literature controls the deviation from the central path in each step by bounding the Newton decrements, given by the potential function $\| \g (y, \eta) \|_{  (\nabla^2 \phi(y))^{-1} }$. When the exact Newton step is taken, this is achieved by the following result.
\begin{lemma}[\cite{r01} exact framework]
Given $ \epsilon_N \in (0, 10^{-2})$, $\eta > 0$, and $\eta^{\new} = \eta (1 + \frac{\epsilon_N}{20 \sqrt{\theta}} ) $. Suppose that $\phi$ is a self-concordant barrier with complexity $\theta \geq 1$ and there is 
\begin{itemize}
    \item {\bf Condition 0.} a feasible solution $y \in \R^m$ satisfies  $ \| \g(y,\eta) \|_{ (\nabla^2 \phi(y))^{-1} } \leq \epsilon_N$,
\end{itemize}
Then $y^{\new} = y - (\nabla^2 \phi(y))^{-1} \g(y,{\eta^{\new} })$ satisfies \begin{align*}
\| \g (y^{\new}, \eta^{\new} ) \|_{  (\nabla^2 \phi(y^{\new}))^{-1} } \leq \epsilon_N.
\end{align*}
\end{lemma}

Further, \cite{jklps20} relaxes the exact $H$ to a PSD approximation version. This framework allows errors in the Hessian matrices.

\begin{lemma}[\cite{jklps20} semi-robust framework]\label{lem:semi_robust_one_step}
Given parameters $ \epsilon_N \in (0, 10^{-2})$, $\eta > 0$, $ \alpha_H \in [1,1+10^{-4}]$, and $\eta^{\new} = \eta (1 + \frac{\epsilon_N}{20 \sqrt{\theta}} ) $. Suppose that $\phi$ is a self-concordant barrier with complexity $\theta \geq 1$ and there is 
\begin{itemize}
    \item {\bf Condition 0.} a feasible solution $y \in \R^m$ satisfies 
    \begin{align*} 
    \| \g(y,\eta) \|_{ (\nabla^2 \phi(y))^{-1} } \leq \epsilon_N,
    \end{align*}
    \item {\bf Condition 1.} a symmetric  matrix $\wt{H} \in \mathbb{S}^{n \times n}_{>0}$ has \begin{align*}
    \alpha_H^{-1} \nabla^2 \phi(y) \preceq \wt{H} \preceq \alpha_H  \nabla^2 \phi(y).
    \end{align*}
\end{itemize}
Then $y^{\new} = y - \wt{H}^{-1} \g(y,{\eta^{\new} })$ satisfies \begin{align*} 
\| \g (y^{\new}, \eta^{\new} ) \|_{  (\nabla^2 \phi(y^{\new}))^{-1} } \leq \epsilon_N.
\end{align*}
\end{lemma}

We propose a more general framework in the following Lemma and we believe it will be useful in the future optimization tasks for semi-definite programming. This framework allows errors in the Hessian matrices, the gradient vectors, and the Newton steps. Also notice that this framework allows even more errors in the Hessian matrices. In Lemma~\ref{lem:semi_robust_one_step}, $\alpha_H \cdot \wt{H}$ must satisfy $\alpha_H^{-2} \cdot \nabla^2 \phi \preceq \alpha_H \cdot \wt{H} \preceq \nabla^2 \phi$ while $\alpha_H^{-2}$ should be close to $1$ as $\alpha_H^{-2} \in [0.99,1]$. In Lemma~\ref{lem:invariant_newton}, $c_H$ can set to smaller constants that are close to $0$ as $c_H \in [10^{-1},1]$. This fact is important to the efficient implementation of hybrid barrier in Section~\ref{sec:hybrid_barrier}.

\begin{lemma}[Our robust Newton step]\label{lem:invariant_newton}
Given any parameters %$ \alpha_S, \alpha_H \in [1, 1+10^{-4}]$, 
$ \epsilon_g , \epsilon_{\delta} \in [0, 10^{-4}], c_H \in [10^{-1},1]$, $0 < \epsilon_N  \leq 10^{-2}$, $\eta >0$, and $\eta^{\new} = \eta (1 + \frac{\epsilon_N}{20 \sqrt{\theta}} ) $. Suppose that $\phi$ is a self-concordant barrier with complexity $\theta \geq 1$. Consider the following conditions.
\begin{itemize}
    \item {\bf Condition 0.} a feasible dual solution $y \in \R^m$ satisfies $ \| \g(y,\eta) \|_{ (\nabla^2 \phi(y))^{-1} } \leq \epsilon_N$,
    \item {\bf Condition 1.} a symmetric  matrix $\wt{H} \in \mathbb{S}^{n \times n}_{>0}$ has 
    \begin{align*}
    c_H \cdot \nabla^2 \phi(y) \preceq \wt{H} \preceq \nabla^2 \phi(y) .
    \end{align*}
    \item {\bf Condition 2.} a vector $\wt{g} \in \R^m$ satisfies 
    \begin{align*}
        \|\wt{g} - \g(y,\eta^{\new} )\|_{(\nabla^2 \phi(y))^{-1}} \leq \epsilon_g \cdot \|\g(y,\eta^{\new})\|_{(\nabla^2 \phi(y))^{-1}} .
    \end{align*}
    \item {\bf Condition 3.} a vector $\wt{\delta}(y) \in \R^m$ satisfies \begin{align*}
    \|\wt{\delta}(y) - (- \wt{H}^{-1} \wt{g})\|_{\nabla^2 \phi(y)} \leq \epsilon_{\delta} \cdot \| \wt{H}^{-1} \wt{g} \|_{\nabla^2 \phi(y)}.
    \end{align*}
\end{itemize}
Suppose {Condition 0,1,2,3} hold. 
Then $y^{\new} = y + \wt{\delta}(y)$ satisfies 
\begin{align*}
\| \g (y^{\new}, \eta^{\new} ) \|_{  (\nabla^2 \phi(y^{\new}))^{-1} } \leq \epsilon_N.
\end{align*}
Furthermore, Condition 1 can also be replaced by the following
\begin{itemize}
    \item {\bf Condition 1'.} a symmetric matrix $\wt{H} \in \mathbb{S}^{n \times n}_{>0}$ satisfies 
    \begin{align*}\alpha_H^{-1} \cdot \nabla^2 \phi(y)  \preceq \wt{H} \preceq \alpha_H \cdot \nabla^2 \phi(y),
    \end{align*}
    where $\alpha_H \in [1,1+10^{-4}]$.
\end{itemize}
\end{lemma}

\begin{remark}
Notice that the error parameters $\epsilon_N,\epsilon_g, \epsilon_\delta,\alpha_H, c_H$ do not depend on the dimension nor the number of iterations. The constants $10^{-4},10^{-2},1/10$ are chosen only for simplicity. In general, if one needs smaller $\epsilon_N$, then one should inflict smaller errors in $\epsilon_g, \epsilon_\delta, \alpha_H, c_H$.
\end{remark}
 
\begin{proof}

First consider the case when Condition 0, 1', 2, 3 hold. Condition 1 is a slightly different condition than Condition 1'. Then, we explain how to modify the proof from condition 1' to condition 1.

By triangle inequality of local norm we have
\begin{align}\label{eq:approx_newton_error_1}
    &~ \|\wt{\delta}(y) - (- n(y, \eta^\new))\|_y \notag \\
    \leq &~ { \|\wt{\delta}(y) - (- \wt{H}^{-1}\wt{g})\|_y } +  \|(\wt{H}^{-1} - (\nabla^2 \phi(y))^{-1})\wt{g}\|_y +  \|(\nabla^2 \phi(y))^{-1}\wt{g} - n(y,\eta^\new)\|_y .
\end{align}

For the second term, Condition 1 and 2 give
\begin{align*}
    \|(\wt{H}^{-1} - (\nabla^2 \phi(y))^{-1})\wt{g}\|_y^2 = &~ \wt{g}^\top (\wt{H}^{-1} - (\nabla^2 \phi({S}))^{-1})\nabla^2 \phi(S)(\wt{H}^{-1} - (\nabla^2 \phi({S}))^{-1})\wt{g}\notag \\
    = &~ \wt{g}^\top (\wt{H}^{-1}\nabla^2 \phi(S)\wt{H}^{-1} -2 \wt{H}^{-1} + (\nabla^2 \phi({S}))^{-1})\wt{g}\notag\\
    \leq &~ (\alpha_H^2 - 2\alpha_H^{-1} + 1) \cdot \wt{g}^\top  (\nabla^2 \phi({S}))^{-1}\wt{g}\notag\\
    = &~ (\alpha_H^2 - 2\alpha_H^{-1} + 1) \cdot \|\wt{g}\|_{(\nabla^2 \phi(y))^{-1}}^2\\
    \leq &~ (\alpha_H^2 - 2\alpha_H^{-1} + 1) \cdot (1+\epsilon_g)^2 \cdot \|{{\g}(y,\eta^\new)}\|_{(\nabla^2 \phi(y))^{-1}}^2\\
    \leq &~ 0.001 \cdot \|n(y,\eta^\new)\|_y^2.
\end{align*}

For the first term, Condition 1, 2, 3 give
\begin{align*}
    \|\wt{\delta}(y) - (- \wt{H}^{-1}\wt{g})\|_y \leq &~ \epsilon_\delta \cdot \|\wt{H}^{-1} \wt{g}\|_y\\
    \leq &~ \epsilon_\delta \cdot \alpha_H \cdot \|(\nabla^2 \phi({S}))^{-1} \wt{g}\|_y\\
    = &~  \epsilon_\delta \cdot \alpha_H \cdot \| \wt{g}\|_{(\nabla^2 \phi({y}))^{-1}}\\
    \leq &~ \epsilon_\delta \cdot \alpha_H \cdot  (1+\epsilon_g) \cdot \| {\g}(y,\eta^\new)\|_{(\nabla^2 \phi({y}))^{-1}} \\
    \leq & ~ 0.001 \cdot \| n(y,\eta^{\new}) \|_y.
\end{align*}

For the third term, Condition 2 gives
\begin{align*}
    \|(\nabla^2 \phi(y))^{-1}\wt{g} - n(y,\eta^\new)\|_y = &~  \|(\nabla^2 \phi(y))^{-1}\wt{g} - (\nabla^2 \phi(y))^{-1}{\g}(y,\eta^\new)\|_y \\%+ \|(\nabla^2 \phi(y))^{-1}{g}(y,\eta^\new) - n(y,\eta^\new)\|_y\\
    \leq &~  \epsilon_g \cdot \|(\nabla^2 \phi(y))^{-1}{\g}(y,\eta^\new)\|_y \\ %+ (\alpha_S-1)\|n(y,\eta^\new)\|_y \\
    \leq & ~ 0.001 \cdot \| n(y,\eta^{\new}) \|_y.
\end{align*}

Combining the above bounds, we have
\begin{align}\label{eq:approx_newton_error}
    \|\wt{\delta}(y) - (- n(y,\eta^\new))\|_y 
    \leq &~ 0.1 \cdot \|n(y,\eta^\new)\|_y.
\end{align}

Combing with Lemma~\ref{lem:eta_move} and Eq.~\eqref{eq:approx_newton_error},
\begin{align}\label{eq:diff_delta_and_n}
    \|\wt{\delta}(y) \|_y \leq 1.1 \cdot \| n(y,\eta^{\new}) \|_y \leq 2 \epsilon_N, \text{~and~} \|\wt{\delta}(y) - (- n(y,\eta^\new))\|_y  \leq 0.3 \cdot \epsilon_N
\end{align}

Using Lemma~\ref{lem:y_move}, we have
\begin{align*}
    \| n(y^\new,\eta^\new) \|_{ y^\new } 
    \leq & ~ 2 \cdot ( \|\wt{\delta}(y)\|_{y}^2 + \| \wt{\delta}(y) - n(y,{\eta^{\new} })\|_y )  \\
    \leq & ~ 2 \cdot ( 4 \cdot \epsilon_N ^2 + 0.3 \cdot \epsilon_N ) \\
    \leq & ~ \epsilon_N
\end{align*}
where the second step follows from Eq.~\eqref{eq:diff_delta_and_n} and the last step follows from choice of $\epsilon_N$.

Next, we consider the case when Condition 0, 1, 2, 3 hold. By triangle inequality of local norm we still have
\begin{align*}
    &~ \|\wt{\delta}(y) - (- n(y, \eta^\new))\|_y \notag \\
    \leq &~ { \|\wt{\delta}(y) - (- \wt{H}^{-1}\wt{g})\|_y } +  \|(\wt{H}^{-1} - (\nabla^2 \phi(y))^{-1})\wt{g}\|_y +  \|(\nabla^2 \phi(y))^{-1}\wt{g} - n(y,\eta^\new)\|_y .
\end{align*}
For the second term,
\begin{align*}
    \| (\wt{H}^{-1} - (\nabla^2 \phi({S}))^{-1})\wt{g}\|_y = &~ \| (I - \nabla^2 \phi(S) \wt{H}^{-1})\wt{g}\|_{(\nabla^2 \phi({S}))^{-1}}\\
    \leq &~ \|I - \nabla^2 \phi(S) \wt{H}^{-1}\|_{(\nabla^2 \phi({S}))^{-1}} \cdot \|\wt{g}\|_{(\nabla^2 \phi({S}))^{-1}}\\
    = &~ \max_{v \in \R^m} \frac{\langle v, ((\nabla^2 \phi({S}))^{-1} - \wt{H}^{-1}) v \rangle}{\langle v, (\nabla^2 \phi({S}))^{-1} v \rangle} \cdot \|\wt{g}\|_{(\nabla^2 \phi({S}))^{-1}}\\
    \leq &~ (1-c_H) \cdot \|\wt{g}\|_{(\nabla^2 \phi({S}))^{-1}}\\
    \leq &~ (1-c_H) \cdot (1+\epsilon_g) \cdot \|n(y,\eta^\new)\|_y
\end{align*}
where the second step comes from H\"{o}lder's inequality, the third step comes from the definition of matrix norm, the penultimate step comes from $0 \preceq (\nabla^2 \phi({S}))^{-1} - \wt{H}^{-1} \preceq (1-c_H) \cdot (\nabla^2 \phi({S}))^{-1}$, and the final step comes from Condition 2. 

For the first term,
\begin{align*}
    \|\wt{\delta}(y) - (- \wt{H}^{-1}\wt{g})\|_y \leq &~ \epsilon_\delta \cdot \|\wt{H}^{-1} \wt{g}\|_y\\
    \leq &~ \epsilon_\delta \cdot \|(\nabla^2 \phi({S}))^{-1} \wt{g}\|_y\\
    = &~  \epsilon_\delta \cdot  \| \wt{g}\|_{(\nabla^2 \phi({y}))^{-1}}\\
    \leq &~ \epsilon_\delta \cdot  (1+\epsilon_g) \cdot \| {\g}(y,\eta^\new)\|_{(\nabla^2 \phi({y}))^{-1}} \\
    \leq & ~ 0.001 \cdot \| n(y,\eta^{\new}) \|_y
\end{align*}
where we use Condition 3 in the first step, we use $\wt{H}^{-1} \preceq (\nabla^2 \phi(S))^{-1}$ in the second step, the penultimate step uses condition 2, and the final step uses the choice of $\epsilon_\delta,\epsilon_g$.

For the third term, Condition 2 gives
\begin{align*}
    \|(\nabla^2 \phi(y))^{-1}\wt{g} - n(y,\eta^\new)\|_y = &~  \|(\nabla^2 \phi(y))^{-1}\wt{g} - (\nabla^2 \phi(y))^{-1}{\g}(y,\eta^\new)\|_y \\%+ \|(\nabla^2 \phi(y))^{-1}{g}(y,\eta^\new) - n(y,\eta^\new)\|_y\\
    \leq &~  \epsilon_g \cdot \|(\nabla^2 \phi(y))^{-1}{\g}(y,\eta^\new)\|_y \\ %+ (\alpha_S-1)\|n(y,\eta^\new)\|_y \\
    \leq & ~ 0.001 \cdot \| n(y,\eta^{\new}) \|_y.
\end{align*}

Thus we can replace the Eq.~\eqref{eq:approx_newton_error} by
\begin{align}\label{eq:approx_newton_error_2}
    \|\wt{\delta}(y) - (- n(y,\eta^\new))\|_y 
    \leq &~ (0.01 + 1.004 \cdot (1-c_H)) \cdot \|n(y,\eta^\new)\|_y \notag \\
    \leq &~ 0.92 \cdot \|n(y,\eta^\new)\|_y.
\end{align}

Combining with Lemma~\ref{lem:eta_move} and Eq.~\eqref{eq:approx_newton_error_2}, we have
\begin{align*}
    \|\wt{\delta}(y) \|_y \leq 2 \epsilon_N \leq 0.02.
\end{align*}

Using the above inequality and Lemma~\ref{lem:y_move}, we have
\begin{align*}
    \| n(y^\new,\eta^\new) \|_{ y^\new } 
    \leq & ~ 0.98^{-2} \cdot ( \|\wt{\delta}(y)\|_{y}^2 + \| \wt{\delta}(y) - n(y,{\eta^{\new} })\|_y )  \\
    \leq & ~ 0.98^{-2} \cdot ( 4\epsilon_N ^2 + 0.92 \epsilon_N ) \\
    \leq & ~ \epsilon_N.
\end{align*}

Thus, we complete the proof.
\end{proof}

\subsection{\texorpdfstring{$\eta$}{~} move}\label{sec:robust_ipm_eta_move}

\begin{lemma}[$\eta$ move]\label{lem:eta_move}
Let $\epsilon_N \in (0,10^{-2})$, $\| n(y,\eta) \|_y \leq \epsilon_N$ and $\eta^{\new} = \eta (1 + \frac{\epsilon_N}{20 \sqrt{\theta}} ) $. Suppose $\phi$ is a self-concordant barrier with complexity $\theta \geq 1$. 
We have
\begin{align*}
    \|n(y,\eta^\new)\|_y\leq (1+ {\epsilon_N}/{20}) \| n(y,\eta) \|_y + \epsilon_N/20 \leq  1.06 \epsilon_N.
\end{align*}
\end{lemma}
\begin{proof}

Denote the Newton step by $n (y,\eta) := (\nabla^2 \phi(y))^{-1}g (y,\eta)$, thus
\begin{align*}
n (y,\eta) = (\nabla^2 \phi(y))^{-1}(\eta b - \mathsf{g}(y) )= \eta b_y - (\nabla^2 \phi(y))^{-1} \mathsf{g}(y) = \eta b_y - \mathsf{g}_y(y)
\end{align*}
(here $\mathsf{g}(y)$ is the gradient of barrier function).

We also have
\begin{align*}
    n (y,\eta^{\new}) = \eta^{\new} b_y - \mathsf{g}_y(y)
\end{align*}
Combining the above two equations, we have
\begin{align*}
    ( n (y,\eta^{\new}) + \mathsf{g}_y(y) ) \eta =  ( n (y,\eta) + \mathsf{g}_y(y) ) \cdot \eta^{\new}
\end{align*}
which implies that
\begin{align*}
    n (y,\eta^{\new}) = & ~ \frac{\eta^{\new}}{\eta}( n (y,\eta) + \mathsf{g}_y(y) ) -  \mathsf{g}_y(y) \\
    = & ~ \frac{\eta^\new}{\eta} n(y,\eta)  + (\frac{\eta^\new}{\eta}-1)\mathsf{g}_y(y)  
\end{align*}

Since the complexity value of barrier functional $\phi$ is $\theta$, 
\begin{align*}%\label{eq:newton_step_sc}
    \|n(y,\eta^\new)\|_y = &~ \Big\|  \frac{\eta^\new}{\eta} n(y,\eta)  + (\frac{\eta^\new}{\eta}-1)\mathsf{g}_y(y) \Big\|_y \notag \\
    \leq &~ \frac{\eta^\new}{\eta} \| n(y,\eta) \|_y + \big|\frac{\eta^\new}{\eta}-1\big|\sqrt{\theta} \notag \\
    \leq & ~ \frac{\eta^\new}{\eta} \epsilon_N + \big|\frac{\eta^\new}{\eta}-1\big|\sqrt{\theta} \notag \\
    \leq & ~ (1+\frac{\epsilon_N}{20\sqrt{\theta}}) \cdot\epsilon_N + \frac{\epsilon_N}{20} \notag \\
    \leq & ~ 1.06 \cdot \epsilon_N,
\end{align*}
where the second step follows from triangle inequality and $\| \mathsf{g}_y(y)\|_y \leq \sqrt{\theta}$,  we use  $\| n(y,\eta) \|_y \leq \epsilon_N$ in the third step (see Lemma statement), we use  definition of $\eta^{\new}$ in the forth step (see Lemma statement), and we use both $\theta \geq 1$ and $\epsilon_N \in (0,10^{-2})$ in the last step.
\end{proof}

\subsection{\texorpdfstring{$y$}{~} move}\label{sec:robust_ipm_y_move}

\begin{lemma}[$y$ move]\label{lem:y_move}
Let $y^{\new} =y + \delta(y) $ and $\| \delta(y) \|_y \leq 1/4$. Suppose $\phi$ is a self-concordant barrier with complexity $\theta \geq 1$. We have
\begin{align*}
    \| n(y^\new,\eta^\new) \|_{ y^\new } 
    \leq \left( \frac{\|{\delta}(y)\|_{y}} { 1-\| {\delta}(y) \|_{ y } } \right)^2 + \frac{  \|\delta(y) - (- n(y,{\eta^{\new} }))\|_y}{1-\| {\delta}(y)\|_{ y }}
\end{align*}
Further, we have
\begin{align*}
    \| n(y^\new,\eta^\new) \|_{ y^\new } \leq 2 \cdot ( \| \delta(y) \|_y^2 + \|\delta(y) - (-n(y,{\eta^{\new} }))\|_y ).
\end{align*}
\end{lemma}
\begin{proof}

We compute the improvement by approximate Newton step. First notice that
\begin{align}\label{eq:y_move_1}
    \| n(y^\new,\eta^\new ) \|_{ y^\new }^2 = &~ \|{\nabla^2 \phi(y)}_y(y^\new)^{-1} \g_y(y^{\new},\eta^\new )\|^2_{y^\new}\notag \\
    = & ~ \left\langle {\nabla^2 \phi}(y^{\new}) {\nabla^2 \phi}_y(y^\new)^{-1} \g_y(y^\new,\eta^\new)  ,{\nabla^2 \phi}_y(y^\new)^{-1}\g_y(y^\new,\eta^\new)\right\rangle \notag\\
     = & ~ \left\langle  \nabla^2 \phi(y) \g_y(y^\new,\eta^\new)  ,{\nabla^2 \phi}_y(y^\new)^{-1}\g_y(y^\new,\eta^\new)\right\rangle \notag\\
    = &~ \left\langle \g_y(y^\new,\eta^\new),{\nabla^2 \phi}_y(y^\new)^{-1}\g_y(y^\new,\eta^\new)\right\rangle_y \notag\\
    \leq &~ \| ( {\nabla^2 \phi}_y(y^\new) )^{-1}\|_y^2 \cdot \|  \g_y(y^\new, \eta^{\new})\|_y^2.
\end{align}
where we use definition of operator norm in the final step.

By Eq.~\eqref{eq:sc_1} (in Theorem~\ref{thm:r01}), 
\begin{align}\label{eq:y_move_2}
    \| (   {\nabla^2 \phi}_y(y^\new) )^{-1}\|_y \leq \frac{1}{(1-\|y^\new - y\|_y)^2} = \frac{1}{(1-\|{\delta}(y)\|_y)^2}.
\end{align}

By Lemma~\ref{lem:integral_local_norm}, we have
\begin{align*}
    \g_y(y^\new,\eta^\new) = & ~ g _y(y,\eta^\new) + \int_0^1 {\nabla^2 \phi}_y(y+t(y^\new - y)) (y^\new - y) \d t \\
    = & ~ n(y,{\eta^{\new} }) + \int_0^1 {\nabla^2 \phi}_y(y+t(y^\new - y)) (y^\new - y) \d t \\
    = & ~ ( n(y,{\eta^{\new} }) + (y^{\new}-y) ) + \int_0^1 ( {\nabla^2 \phi}_y(y+t(y^\new - y)) -I ) (y^\new - y) \d t 
\end{align*}
We can upper bound the first term under local norm as follow:
\begin{align}\label{eq:y_move_3}
    \|  n(y,{\eta^{\new} }) + (y^{\new}-y) \|_y = \| -n(y,{\eta^{\new} }) - \delta(y) \|_y.
\end{align}
We can upper bound the second term under local norm as follow:
\begin{align}\label{eq:y_move_4}
     &~ \Big\| \int_0^1 \left({\nabla^2 \phi}_y(y+t(y^\new - y)) - I\right)(y^\new - y) \d t\Big\|_y \notag\\
    \leq &~ \|y^\new - y\|_y \int_0^1\|{\nabla^2 \phi}_y(y+t(y^\new - y)) - I\|_y \d t 
    \notag\\
    \leq &~ \|y^\new - y\|_y \int_0^1 \left( \frac{1}{(1-t \|y^\new - y\|_y)^2} - 1 \right) \d t  \notag\\
    = &~ \frac{\|y^\new - y\|^2_y}{1-\|y^\new - y\|_y} \notag\\
    = &~ \frac{\|{\delta}(y)\|^2_y}{1-\|{\delta}(y)\|_y} 
\end{align}
where the first step comes from triangle inequality of local norm, the second step comes from Eq.~\eqref{eq:sc_2} (property of self-concordance function, Theorem~\ref{thm:r01}), we use simple integration in the third step comes, finally we use $y^{\new} = y + {\delta}(y)$ in the final step.

Plugging Eq.~\eqref{eq:y_move_2}, Eq.~\eqref{eq:y_move_3}, and Eq.~\eqref{eq:y_move_4} into Eq.~\eqref{eq:y_move_1}, we have
\begin{align*}
    \| n(y^\new,\eta^\new) \|_{ y^\new } 
    \leq & ~ \frac{1}{1-\|\delta(y)\|} \cdot \| \g_y( y^{\new} , \eta^{\new} ) \|_y \\
    \leq & ~ \frac{1}{1-\|\delta(y)\|} \cdot \left( \frac{\|{\delta}(y)\|_{y}^2} { 1-\| {\delta}(y) \|_{ y } } + \|\delta(y) - (- n(y,{\eta^{\new} }))\|_y  \right) \\
    \leq & ~ \left( \frac{\|{\delta}(y)\|_{y}} { 1-\| {\delta}(y) \|_{ y } } \right)^2 + \frac{  \|\delta(y) - (- n(y,{\eta^{\new} }))\|_y}{1-\| {\delta}(y)\|_{ y }}.
\end{align*}
This completes the proof.
\end{proof}

\subsection{Integral under local norm}\label{sec:robust_ipm_integral}

\begin{lemma}\label{lem:integral_local_norm}
Let $H = g'$ and $x,y\in D_f$. It holds that
\begin{align*} 
    (\nabla^2 \phi(x))^{-1} (\g(y)-\g(x)) = \int_0^1 (\nabla^2 \phi(x))^{-1} {\nabla^2 \phi}(x+ t (y-x)) (y - x) \d t 
\end{align*}
\end{lemma}
\begin{proof}

It is sufficient to show
\begin{align*}
    \g(y)-\g(x) = \int_0^1 {\nabla^2 \phi}(x+ t (y-x)) (y - x) \d t
\end{align*}
By definition of the integral, it is sufficient to prove that for all $w$
\begin{align*}
    \langle \g(y) - \g(x), w \rangle = \int_0^1 \langle {\nabla^2 \phi}(x+t(y-x)) (y-x) , w \rangle \d t .
\end{align*}
Fix arbitrary $w$ and consider the functional 
\begin{align*}
    \psi(t):= \langle g( x + t(y-x) ), w \rangle.
\end{align*}
The basic Calculus gives
\begin{align*}
    \psi(1) - \psi(0) = \int_0^1 \psi'(t) \d t
\end{align*}
which by definition of $\psi$ is equivalent to
\begin{align*}
    \langle \g(y) - \g(x) , w \rangle = \int_0^1 \langle {\nabla^2 \phi}(x+t(y-x)) (y-x) , w \rangle \d t.
\end{align*}
\end{proof}

\subsection{Approximate dual optimality}\label{sec:robust_ipm_duality}

We make use of the following lemma that bounds the duality gap by $\eta$ and the Newton decrement.

\begin{lemma}[Approximate optimality]\label{lem:approximate_optimality}
Suppose $0 < \epsilon_N \leq 10^{-2}$. Let $\eta \geq 1$ denote a parameter. Let $y \in \mathbb{R}^m$ be dual feasible solution. % and parameter $\eta$ satisfy the following  %on Newton step size: 
Assume 
\begin{align*}
\g(y, \eta)^\top (\nabla^2 \phi(y))^{-1} \g(y,\eta) \leq \epsilon_N^2.
\end{align*}
Assume that $y^*$ is an optimal solution to the Eq.~\eqref{eq:sdp_dual}. Suppose $\phi$ is a self-concordant barrier with complexity $\theta \geq 1$. Then we have
\begin{align*}
\langle b , y  \rangle \leq \langle b , y^*  \rangle + \frac{\theta}{\eta} \cdot (1 + 2 \epsilon_N) .  
\end{align*}
\end{lemma}

\begin{proof}
Let $y(\eta)$ denote the optimal solution to the following optimization problem:
\begin{align*}
    \min_{y \in \R^m }  \eta \cdot \langle b , y \rangle + \phi(y)
\end{align*}
where $\phi(y) = - \log \det S(y)$ is the log barrier. Then due to optimality condition, $\eta b + \mathsf{g}(y(\eta)) = 0$. Therefore
\begin{align}\label{eq:by(eta)-byast}
    \langle b , y(\eta) \rangle - \langle b , y^\ast \rangle =  \frac{1}{\eta}\langle \mathsf{g}(y(\eta)), y^\ast - y(\eta) ) \rangle    \leq  \frac{\theta}{\eta}
\end{align}
where the last step comes from Theorem~\ref{thm:r01_2.3.4}. 

Furthermore, we have
\begin{align}
    \langle b , y  \rangle - \langle b , y(\eta) \rangle = &~ \frac{1}{\eta}\langle \mathsf{g}(y(\eta)),  y(\eta) - y \rangle \notag\\
    \leq &~ \frac{1}{\eta}\|\mathsf{g}(y(\eta))\|_{(\nabla^2 \phi(y(\eta)) )^{-1}} \cdot \|y - y(\eta)\|_{\nabla^2 \phi(y(\eta))} \notag\\
    \leq &~ \frac{\theta}{\eta} \cdot \|y - y(\eta)\|_{\nabla^2 \phi(y(\eta))} \label{eq:by-by(eta)}
\end{align}
where the last step used the complexity value of barrier is $\theta$.
For $\|y - y(\eta)\|_{\nabla^2 \phi(y)}$, we have
\begin{align}
    \|y - y(\eta)\|^2_{\nabla^2 \phi(y(\eta))} \leq &~ \|y - y(\eta)\|^2_{\nabla^2 \phi(y)} \cdot \underset{v}{\sup}\frac{\|v\|^2_{\nabla^2 \phi(y(\eta))}}{\|v\|^2_{\nabla^2 \phi(y)}} \notag \\
    = &~ \|y - y(\eta)\|^2_{\nabla^2 \phi(y)} \cdot \|(\nabla^2 \phi(y))^{-1}{\nabla^2 \phi(y(\eta))}\|_{\nabla^2 \phi(y)} \notag \\
    \leq &~ \|y - y(\eta)\|^2_{\nabla^2 \phi(y)} \cdot  (1- \|y - y(\eta)\|^2_{\nabla^2 \phi(y)} )^{-2} \label{eq:local_norm_y-y(eta)}
\end{align}
where the second step comes from the definition of operator norm and the third step comes from Theorem~\ref{thm:r01}.
By Theorem~\ref{thm:r01_2.2.8}, 
\begin{align*}
    \|y - y(\eta)\|_{\nabla^2 \phi(y)} \leq \|n(y,\eta)\|_{\nabla^2 \phi(y)} + \frac{2\|n(y,\eta)\|_{\nabla^2 \phi(y)}^2}{\left(1- \|n(y,\eta)\|_{\nabla^2 \phi(y)}\right)^3} \leq 1.2 \cdot \epsilon_N,
\end{align*}
thus back to Eq.~\eqref{eq:local_norm_y-y(eta)}, we have $\|y - y(\eta)\|_{\nabla^2 \phi(y(\eta))} \leq 2 \epsilon_N$. Therefore in Eq.~\eqref{eq:by-by(eta)}, we obtain
\begin{align}\label{eq:by-by(eta)_2}
    \langle b , y  \rangle - \langle b , y(\eta) \rangle 
    \leq  \frac{\theta}{\eta} \cdot 2 \epsilon_N.
\end{align}
Combining Eq.~\eqref{eq:by-by(eta)_2} and Eq.~\eqref{eq:by(eta)-byast}, we complete the proof.
\end{proof}

\begin{theorem}[Robust barrier method]\label{thm:approx_central_path_dual}
Consider a semidefinite program in Eq.~\eqref{eq:sdp_dual}. 
%Assume that any feasible solution $X \in \mathbb{S}^{n \times n}_{\geq 0}$ satisfies $\| X \|_{2} \leq R$. 
Suppose in each iteration, the $\wt{S}, \wt{H},\wt{g},\wt{\delta}$ are computed in Line~\ref{lin:approx_S}, Line~\ref{lin:approx_H}, Line~\ref{lin:approx_g}, Line~\ref{lin:approx_delta} of Algorithm~\ref{alg:robust_ipm} such that Condition 1 or 1' $\&$ Condition 2 $\&$ Condition 3 in Lemma~\ref{lem:invariant_newton} hold. 
Suppose $\phi$ is a self-concordant barrier with complexity $\theta \geq 1$. Assume $y^*$ is an optimal solution to the dual formulation Eq.~\eqref{eq:sdp_dual}. Then given a feasible initial solution that satisfies the invariant $\g(y, \eta)^\top (\nabla^2 \phi(y))^{-1} \g(y,\eta) \leq \epsilon_N^2$, for any error parameter $0 < \epsilon \leq 0.01$ and Newton step size $\epsilon_N$ satisfying  $\sqrt{\epsilon} < \epsilon_N \leq 0.01$, Algorithm~\ref{alg:robust_ipm} outputs,  in $T = 40 \epsilon_N^{-1} \sqrt{\theta} \log(\theta /\epsilon)$ iterations,  a vector $y \in \mathbb{R}^{m}$ s.t. 
\begin{align}\label{eq:dual_optimality_gap}
 b^\top  y  \leq  b^\top y^* + \epsilon^2 .  
\end{align}
Further, for logarithmic barrier $\phi_{\mathrm{log}}$, in each iteration of Algorithm~\ref{alg:robust_ipm}, the following invariant holds: 
\begin{align}\label{eq:promise_sdp}
\| S^{-1/2} S^{\new} S^{-1/2} - I \|_F \leq 1.03 \cdot \epsilon_N.
\end{align}
\end{theorem}

\begin{proof}

Since the invariant $\g(y, \eta)^\top (\nabla^2 \phi(y))^{-1} \g(y,\eta) \leq \epsilon_N^2$ holds at initialization, by Lemma~\ref{lem:invariant_newton} it then holds at any iteration. 
After $T = 40 \epsilon_N^{-1} \sqrt{\theta} \log(\theta /\epsilon)$ iterations, the step size becomes $\eta = (1 + \frac{\epsilon_N}{20 \sqrt{\theta}})^T / (\theta +2) \geq 2\theta/\epsilon^2$.  
By Lemma~\ref{lem:approximate_optimality}, we have
\begin{align*}
    \langle b , y  \rangle \leq \langle b , y^*  \rangle + \frac{\theta}{\eta} \cdot (1 + 2 \epsilon_N) \leq  \langle b , y^*  \rangle  + \epsilon^2.
\end{align*}
This completes the proof of Eq.~\eqref{eq:dual_optimality_gap}.

Finally we prove Eq.~\eqref{eq:promise_sdp} for the log-barrier $\phi_{\mathrm{log}}$. %Let $\Tilde{\delta}_{y,i}$ denote the $i$-th coordinate of vector $\Tilde{\delta}_{y}$. 
It gives
\begin{align*}
    %\| S^{-1/2} S^{\new} S^{-1/2} - I \|^2_F
    \LHS ~\text{in ~Eq.~\eqref{eq:promise_sdp}}
    = &~ \tr\left[ \left(S^{-1/2} (S^\new - S) S^{-1/2}\right)^2\right]\\
    = &~ \tr\left[ S^{-1} (\sum_{i\in[m]} \Tilde{\delta}_{y,i} A_i) S^{-1} (\sum_{i\in[m]} \Tilde{\delta}_{y,i} A_i)\right]\\
    = &~ \sum_{i\in [m]} \sum_{j\in[m]} \Tilde{\delta}_{y,i} \Tilde{\delta}_{y,j} \tr[S^{-1}A_i S^{-1}A_j]\\
    = &~ \Tilde{\delta}_{y}^\top \nabla^2 \phi(y) \Tilde{\delta}_{y}\\
    = &~ \|\Tilde{\delta}_{y}\|^2_{\nabla^2 \phi(y)}
\end{align*}
%\Ruizhe{The above equations need to be reproved for vol barrier.}
where the second step comes from $S^\new - S = \sum_{i\in [m]} \Tilde{\delta}_{y,i} A_i$. It suffices to bound $\|\Tilde{\delta}_{y}\|_{\nabla^2 \phi(y)}$. In fact we have
\begin{align*}
    \|\Tilde{\delta}_{y}\|_{\nabla^2 \phi(y)} \leq &~ \|n(y,\eta^\new)\|_{\nabla^2 \phi(y)} + \|\Tilde{\delta}_{y} - (- n(y,\eta^\new))\|_{\nabla^2 \phi(y)}\\
    \leq &~ 1.01 \|n(y,\eta^\new)\|_{\nabla^2 \phi(y)} \\ 
    \leq &~ 1.01 \cdot \left( (1+\epsilon_N/20)\|n(y,\eta)\|_{\nabla^2 \phi(y)} + \epsilon/20 \right)\\
    \leq &~ 1.03 \cdot \epsilon_N
\end{align*}
where we use triangle inequality in the first step, we use Eq.~\eqref{eq:approx_newton_error} in the second step, we use Lemma~\ref{lem:eta_move} in the  third step  and the last step comes from choice of $\epsilon_N$. This completes the proof of Eq.~\eqref{eq:promise_sdp}.
\end{proof}

\subsection{Our main result}\label{sec:robust_ipm_main}

\begin{algorithm}[!ht]
\caption{Our robust barrier method for SDP.}
\label{alg:robust_ipm}
\begin{algorithmic}[1]
\Procedure{\textsc{SolveSDP}}{$m, n, C, \{ A_i \}_{i=1}^m$, $\mathsf{A} \in \R^{m \times n^2}$, $b \in \R^m$} 
\State \Comment{Initialization} 
\State $\eta \leftarrow \frac{1}{\theta +2}$, ~~$T \leftarrow \frac{40}{\epsilon_N} \sqrt{\theta} \log (\frac{\theta }{\epsilon})$
\State Find initial feasible dual $y \in \R^m$ according to Lemma~\ref{lem:initialization}\Comment{Condition  0 in Lemma~\ref{lem:invariant_newton}}
\For {$t = 1 \to T$} {\bf do} \Comment {Iterations of approximate barrier method}
	\State $\eta^{\new} \leftarrow \eta \cdot (1 + \frac{\epsilon_N}{20 \sqrt{\theta}})$
	\State $\wt{S} \leftarrow \textsc{ApproxSlack} ()$ %\Comment{Condition  1 in Lemma~\ref{lem:invariant_newton}}
	\label{lin:approx_S}
'

\State $\wt{H} \leftarrow \textsc{ApproxHessian}()$ \Comment{Condition  1 in Lemma~\ref{lem:invariant_newton}}\label{lin:approx_H}
	\State $\wt{g} \leftarrow \textsc{ApproxGradient}( )$ \Comment{Condition  2 in Lemma~\ref{lem:invariant_newton}}\label{lin:approx_g}
	\State $\wt{\delta}(y) \leftarrow \textsc{ApproxDelta}()$\Comment{Condition  3 in Lemma~\ref{lem:invariant_newton}}\label{lin:approx_delta} 
	\State $y^{\new} \leftarrow y + \delta(y)$ 
	\State $y \leftarrow y^{\new}$ 
	\Comment{We update variables}
\EndFor
\State \Return an approximate solution to the original problem  \Comment{Lemma~\ref{lem:initialization}}
\EndProcedure
\end{algorithmic}
\end{algorithm}

\begin{theorem}[Robust central path]\label{thm:approx_central_path}
Consider an SDP instance defined in Definition~\ref{def:sdp_primal} with no redundant constraints. Assume that the feasible region is bounded, i.e., $\| X \|_{2} \leq R$. Suppose in each iteration, the $\wt{S}, \wt{H},\wt{g},\wt{\delta}$ are computed in Line~\ref{lin:approx_S}, Line~\ref{lin:approx_H}, Line~\ref{lin:approx_g}, Line~\ref{lin:approx_delta} of Algorithm~\ref{alg:robust_ipm} that satisfy Condition 1 or 1' $\&$ Condition 2 $\&$ Condition 3 in Lemma~\ref{lem:invariant_newton}. Suppose $\phi$ is a self-concordant barrier with complexity $\theta \geq 1$. Assume $X^*$ is an optimal solution to the semidefinite program in Definition~\ref{def:sdp_primal}. Then for any error parameter $0 < \epsilon \leq 0.01$ and Newton step size $\epsilon_N$ satisfying  $\sqrt{\epsilon} < \epsilon_N \leq 0.01$, Algorithm~\ref{alg:robust_ipm} outputs,  in $T = 40 \epsilon_N^{-1} \sqrt{\theta } \log(\theta /\epsilon)$ iterations,  a positive semidefinite matrix $X \in \mathbb{R}^{n \times n}_{\geq 0}$ s.t.
\begin{align}\label{eq:approx_optimality_1}
\begin{array}{l}
\langle C, X \rangle \geq \langle C, X^* \rangle - \epsilon \cdot \| C \|_2 \cdot R,~~ \text{and}\\ \sum_{i = 1 }^m \left| \langle A_i, \widehat{X} \rangle - b_i \right| \leq 4 n \epsilon \cdot \Big( R  \sum_{i = 1}^m \| A_i \|_1 + \| b \|_1 \Big),
\end{array}
\end{align}
%where  and %$\| A_i \|_1$ is the Schatten $1$-norm of matrix $A_i$. 
Furthermore, for logarithmic barrier $\phi_{\mathrm{log}}$, in each iteration of Algorithm~\ref{alg:robust_ipm}, the following invariant holds: 
\begin{align}\label{eq:promise_sdp_2}
\| S^{-1/2} S^{\new} S^{-1/2} - I \|_F \leq 1.03 \cdot \epsilon_N.
\end{align}

\end{theorem}

\begin{proof}
%In the beginning of Algorithm~\ref{alg:robust_ipm}, 
First, we use Lemma~\ref{lem:initialization} to rewrite the semidefinite programming and obtain an initial feasible solution near the dual central path with $\eta = 1/(\theta +2)$.
Thus, the induction hypothesis 
\begin{align*}
\g(y, \eta)^\top (\nabla^2 \phi(y))^{-1} \g(y,\eta) \leq \epsilon_N^2
\end{align*}
holds at the initial of  algorithm.

Say $y$ is the modified semidefinite programming's dual solution.
Theorem~\ref{thm:approx_central_path_dual} shows that $y$ has duality gap  $\leq \epsilon^2$. 

Finally, we use Lemma~\ref{lem:initialization}, to get an approximate solution to the original semidefinite programming satisfying Eq.~\eqref{eq:approx_optimality_1}.
\end{proof}

\section{Hybrid Barrier-Based SDP Solver}\label{sec:hybrid_barrier}

The hybrid barrier \cite{nn89,a00} is another useful barrier function to solve SDP and converges within a smaller number of iteration when $m \leq n$. However, it is hard to be implemented efficiently due to the complex form of Hessian matrices. In this section, we give an efficient algorithm for solving SDP using the hybrid barrier in \cite{a00} that improves the naive implementation in all parameter regimes\footnote{We also improves our straightforward algorithm in most parameter regimes. See Remark~\ref{rem:our_a00} for the different implementations of \cite{a00}.}. 

In Section~\ref{sec:hybrid_prelim}, we review some basic facts on the hybrid barrier for SDP. In Section~\ref{sec:hybrid_alg}, we give the formal version of the algorithm and time complexity result. In Section~\ref{sec:hybrid_approxQ}, we show how to low-rank approximate the change of $Q(S)$. In Section~\ref{sec:hybrid_S_slow}, we prove that the slack variable $S$ changes slowly in each iteration. Section~\ref{sec:hybrid_amortize} contains the amortized analysis for our hybrid barrier SDP solver. Combining them together, we prove the main theorem (Theorem~\ref{thm:vol_formal}) in Section~\ref{sec:hybrid_result}. 

\subsection{Basic facts on the hybrid barrier}\label{sec:hybrid_prelim}
The barrier function is defined as follows:
\begin{align*}
    \phi(y):=225\sqrt{\frac{n}{m}} \cdot \left(\phi_{\mathrm{vol}}(y) + \frac{m-1}{n-1}\cdot \phi_{\log}(y)\right),
\end{align*}
with
\begin{align*}
    \phi_{\mathrm{vol}}(y):=&~ 0.5 \log \det (H(y)),\\
    \phi_{\log}(y):= &~ -\log \det(S(y)),
\end{align*}
Note $H(y)$ and $S(y)$ are defined in Definition~\ref{def:S_H}.

According to our robust IPM framework, for every iteration, we have to calculate/estimate the gradient and Hessian of $\phi_{\mathrm{vol}}(y)$, whose closed-forms are computed in \cite{a00}. More specifically,
\begin{align*}
    (\nabla \phi_{\mathrm{vol}}(y))_i = -\tr[H(S)^{-1} \cdot \A  (S^{-1}A_i S^{-1}\otimes S^{-1}) \mathsf{A}^\top]~~~\forall i\in [m],
\end{align*}
and
\begin{align*}
    \nabla^2\phi_{\mathrm{vol}}(y) = 2Q(S) + R(S) - 2T(S),
\end{align*}
where for any $i,j\in [m]$,
\begin{align*}
    Q(S)_{i,j} = &~ \tr[H(S)^{-1} \A (S^{-1}A_iS^{-1}A_jS^{-1} \otimes_S S^{-1}) \A^\top],\\
    R(S)_{i,j} = &~ \tr[H(S)^{-1} \A (S^{-1}A_iS^{-1} \otimes_S S^{-1}A_jS^{-1}) \A^\top],\\
    T(S)_{i,j}=&~ \tr[H(S)^{-1} \A (S^{-1}A_iS^{-1} \otimes_S S^{-1}) \A^\top H(S)^{-1}\A (S^{-1}A_jS^{-1} \otimes_S S^{-1}) \A^\top ].
\end{align*}

The following fact in \cite{a00} shows that $Q(S)$ is a good PSD approximation of the Hessian $\nabla^2\phi_{\mathrm{vol}}(y)$.

\begin{fact}[\cite{a00}]\label{fac:vol_Q_Sigma}
For $S\succ 0$,
\begin{align*}
    \frac{1}{n}H(S)\preceq Q(S)\preceq \nabla^2 \phi_{\mathrm{vol}}(y) \preceq 3Q(S).
\end{align*}
\end{fact}

We also need the following lower bound on the quadratic form of hybrid barrier's Hessian.
\begin{fact}[\cite{a00}]\label{fac:hybrid_lb}
Let $S\succ 0$. For any $\xi \in \R^m$, we have
\begin{align*}
    \xi^\top \left(Q(S) + \frac{m-1}{n-1}\cdot H(S)\right) \xi \geq ~ \frac{2\sqrt{m}}{1+\sqrt{n}}\cdot \left\|S^{-1/2}\left(\sum_{i=1}^m \xi_i A_i\right) S^{-1/2}\right\|_2^2.
\end{align*}
\end{fact}

\subsection{Efficient implementation via robust SDP framework}\label{sec:hybrid_alg}

\begin{algorithm}[!ht]
\caption{Hybrid Barrier SDP solver.}
\label{alg:sdp_vol}
\begin{algorithmic}[1]
\item {\bf members}\\
%\hspace{4mm} $\{ A_i \}_{i=1}^m$\\
%\hspace{4mm} $\{ A_i \}_{i=1}^m$\\
\hspace{4mm} $S, \wt{S} \in \R^{n \times n}$ \hspace{\fill} $\triangleright$ Slack variables\\
\hspace{4mm} $y \in \R^m$ \hspace{\fill} $\triangleright$ Dual variable \\
\hspace{4mm} $H, Q \in \R^{m \times m}$ \hspace{\fill} $\triangleright$ Parts of the Hessian matrices \\
\hspace{4mm} $\eta \in \R$ \hspace{\fill} $\triangleright$ Learning rate \\
\hspace{4mm} $\A \in \R^{m \times n^2}$ \hspace{\fill} $\triangleright$ Batched constraint matrix \\
\hspace{4mm} $G \in \R^{m \times m}$ \hspace{\fill} $\triangleright$ The inverse of the Hessian matrix \\
{\bf end members}\\

\item {\bf procedure} {\textsc{HybridBarrier}}{$(m, n, C, \{ A_i \}_{i=1}^m$, $b \in \R^m)$} \\
\hspace{4mm} $\textsc{Initialize}(m, n, C, \{ A_i \}_{i=1}^m$, $b \in \R^m)$\\
\hspace{4mm} {\bf for} {$t = 1 \to T$} {\bf do} \hspace{\fill} $\triangleright$ {Iterations of approximate barrier method}\\
	\hspace{8mm} $\eta^{\new} \leftarrow \eta \cdot \left(1 + \frac{\epsilon_N}{20 (mn)^{1/4}}\right)$\\
	\hspace{8mm} $g_{\eta^{\new}}(y) \leftarrow \textsc{HybridGradient}(\eta^\new,b,C,\{ A_i \}_{i=1}^m)$\\
	\hspace{8mm} $\delta_y \leftarrow - G \cdot g_{\eta^{\new}}(y)$ \hspace{\fill} $\triangleright$ {Update on $y \in \R^m$}\label{lin:vol_delta_y}\\
	%\hspace{8mm} Let $\wt{\delta}_y$ denote the changes generated by Quantum machine \Zhao{I add this}\\
	\hspace{8mm} $y^{\new} \leftarrow y + \delta_y$\\
	\hspace{8mm} $S^{\new} \leftarrow \sum_{i \in [m]} (y^{\new})_i \cdot A_i - C$ \label{lin:vol_S_new}\\
	\hspace{8mm} $V_1, V_2 \leftarrow \textsc{LowRankSlackUpdate}(S^{\new}, \wt{S})$ \hspace{\fill} $\triangleright$ {$V_1,V_2 \in \R^{n \times r_t}$. Algorithm~\ref{alg:low_rank_slack_update}.} \label{lin:vol_low_rank_update} \\
	\hspace{8mm} $\wt{S}^{\new} \leftarrow \wt{S} + V_1 V_2^{\top}$ \hspace{\fill} $\triangleright$ {Approximate slack computation}\label{lin:vol_slack_update} \\
	\hspace{8mm} $V_3 \leftarrow - \wt{S}^{-1} V_1 (I + V_2^{\top} \wt{S}^{-1} V_1)^{-1}$ \Comment{$V_3 \in \R^{n \times r_t}$} \label{lin:vol_V_3}\\
	\hspace{8mm} $V_4 \leftarrow \wt{S}^{-1} V_2$ \hspace{\fill} $\triangleright$ {$V_4 \in \R^{n \times r_t}$} \label{lin:vol_V_4}\\
	\hspace{8mm} $y \leftarrow y^{\new}$,  $S \leftarrow S^{\new}$, $\wt{S} \leftarrow \wt{S}^{\new}$, $\eta \leftarrow \eta^\new$ \hspace{\fill} $\triangleright$ {Update variables} \label{lin:vol_update_variables}\\
	\hspace{8mm} $(\wt{Q},\wt{H},G^{\new}) \leftarrow \textsc{HybridHessian}(V_3,V_4)$ \Comment{$G^{\new} \in \R^{m \times m}$} \label{lin:vol_hessian}\\
	\hspace{8mm} \hspace{\fill} $\triangleright$ {Hessian inverse computation using Woodbury identity}\\
	\hspace{8mm} $Q \leftarrow \wt{Q}$, $H \leftarrow \wt{H}$, $G \leftarrow G^{\new}$ \hspace{\fill} $\triangleright$ {Update matrices}\\
\hspace{4mm} {\bf end for}\\
\hspace{4mm} \Return an approximate solution to the original problem \Comment{Lemma~\ref{lem:initialization}}\\
{\bf end procedure}\\

\item {\bf procedure} {\textsc{Intialize}}{$(m, n, C, \{ A_i \}_{i=1}^m$, $\mathsf{A} \in \R^{m \times n^2}$, $b \in \R^m)$} \\
%\hspace{4mm} \hspace{\fill} $\triangleright$ {Initialization} \\
\hspace{4mm} Construct $\mathsf{A} \in \R^{m \times n^2}$ by stacking $m$ vectors $\vect[A_1], \vect[A_2], \cdots, \vect[A_m] \in \R^{n^2}$ \\
\hspace{4mm} $\eta \leftarrow \frac{1}{(mn)^{1/2}+2}$, ~~$T \leftarrow \frac{40}{\epsilon_N} (mn)^{1/4} \log (\frac{mn}{\epsilon})$\\
\hspace{4mm} Find initial feasible dual vector $y \in \R^m$ according to Lemma~\ref{lem:initialization}\\
\hspace{4mm} $S \leftarrow \sum_{i \in [m]} y_i \cdot A_i - C$, ~~ $\wt{S} \leftarrow S$ \hspace{\fill} $\triangleright$ {$S, \wt{S} \in \R^{n \times n}$}\\ %\label{line:S_init}\\
\hspace{4mm} $H(S) \leftarrow \A (S^{-1} \otimes S^{-1})\A^\top$\\
\hspace{4mm} {\bf for} {$i = 1, \cdots, m$} {\bf do} \\
\hspace{8mm} $Q(S)_{i,j} \leftarrow \tr[H^{-1} \A (S^{-1}A_iS^{-1}A_jS^{-1} \otimes_S S^{-1}) \A^\top]$\\
\hspace{8mm} $R(S)_{i,j} \leftarrow \tr[H^{-1} \A (S^{-1}A_iS^{-1} \otimes_S S^{-1}A_jS^{-1}) \A^\top]$\\
\hspace{8mm} $O(S)_{i,j} \leftarrow \tr[H^{-1} \A (S^{-1}A_iS^{-1} \otimes_S S^{-1}) \A^\top H^{-1}\A (S^{-1}A_jS^{-1} \otimes_S S^{-1}) \A^\top ]$\\
\hspace{4mm} {\bf end for} \\
\hspace{4mm} $G \leftarrow \left(225\cdot \sqrt{\frac{n}{m}}\cdot \left(2Q(S) + R(S) -2 O(S) + \frac{m-1}{n-1}\cdot H(S)\right)\right)^{-1}$ \hspace{\fill} $\triangleright$ {$G \in \R^{m \times m}$}\\ %\label{line:G_init}\\
%\hspace{4mm} \hspace{\fill} $\triangleright$ {Maintain $G = \wt{H}^{-1}$ where $\wt{H} := \mathsf{A} \cdot (\wt{S}^{-1} \otimes \wt{S}^{-1}) \cdot \mathsf{A}^{\top}$}\\
%\hspace{4mm} {\bf Return} $(S,\wt{S},y,H,Q,R,T,\eta,\A,G)$\\
{\bf end procedure}
\end{algorithmic}
\end{algorithm}

\begin{algorithm}[!ht]
\caption{Hybrid Barrier SDP solver, continued.}
\label{alg:sdp_vol_cont}
\begin{algorithmic}[1]

\item {\bf procedure} {\textsc{HybridGradient}}{$(m, n, C, \{ A_i \}_{i=1}^m$, $b \in \R^m)$} \\
\hspace{4mm} {\bf for} {$i = 1, \cdots, m$} {\bf do} \\
		\hspace{8mm} $\nabla \phi_{\mathrm{log}}(y)_i \leftarrow -\tr [S^{-1} \cdot A_j]$ \hspace{\fill} $\triangleright$ {Gradient of $\phi_{\mathrm{log}}$} \label{lin:log_grad}\\
		\hspace{8mm} $\nabla \phi_{\mathrm{vol}}(y)_i \leftarrow -\tr[H(S)^{-1} \cdot \A  (S^{-1}A_i S^{-1}\otimes S^{-1}) \mathsf{A}^\top]$ \hspace{\fill} $\triangleright$ {Gradient of $\phi_{\mathrm{vol}}$}\label{lin:vol_grad}\\
\hspace{4mm} {\bf end for}\\
\hspace{4mm} $g_{\eta^{\new}}(y) \leftarrow \eta^\new b - 225\sqrt{\frac{n}{m}}\cdot \left(\nabla \phi_{\mathrm{vol}}(y) + \frac{m-1}{n-1}\cdot \nabla \phi_{\mathrm{log}}(y)\right)$\\
\hspace{4mm} {\bf return} $g_{\eta^{\new}}(y)$\\
{\bf end procedure}\\

\item {\bf procedure} {\textsc{HybridHessian}}{$(V_3,V_4)$} \\
%\hspace{4mm} \hspace{\fill} $\triangleright$ {Initialization} \\
\hspace{4mm} {\bf for} {$i, j = 1, \cdots, m$} {\bf do} \\
	%\hspace{8mm} {\bf for} {$i = 1, \cdots, m$} {\bf do} \\
	\hspace{8mm} $\wt{H}_{i,j} \leftarrow H_{i,j} +  \tr[{S}^{-1}A_i V_3V_4^\top A_j] + \tr[{S}^{-1}A_j V_3V_4^\top A_i] + \tr[V_3V_4^\top A_i V_3V_4^\top A_j]$ \label{lin:vol_H}\\
	\hspace{8mm} $\wt{Q}_{i,j} \leftarrow Q_{i,j} + \tr[H^{-1} \A (S^{-1}A_i V_3V_4^\top A_jS^{-1} \otimes_S S^{-1}) \A^\top]$ \label{lin:vol_Q}\\
	%\hspace{8mm} $\wt{R}_{i,j} \leftarrow R_{i,j} + \tr[H^{-1} \A (S^{-1}A_i^{1/2}S^{1/2} \otimes S^{-1}A_j^{1/2}S^{1/2}) (S^{-1} \otimes_S V_3V_4^\top) (S^{1/2}A_i^{1/2} S^{-1} \otimes S^{1/2}A_j^{1/2}S^{-1}) \A^\top]$\label{lin:vol_R}\\
	%\hspace{8mm} $\wt{T}_{i,j} \leftarrow \tr[H(S)^{-1} \A (S^{-1}A_iS^{-1} \otimes_S S^{-1}) \A^\top H(S)^{-1}\A (S^{-1}A_jS^{-1} \otimes_S S^{-1}) \A^\top ]$\label{lin:vol_T}\\
\hspace{4mm} {\bf end for}\\
\hspace{4mm} $G^{\new} \leftarrow \left(225 \cdot 1.001 \cdot \sqrt{\frac{n}{m}}\cdot \left(3\cdot \wt{Q} + \frac{m-1}{n-1}\cdot \wt{H}\right)\right)^{-1} \in \R^{m \times m}$ \hspace{\fill} $\triangleright$ {$G \in \R^{m \times m}$}\\
%\hspace{4mm} $G^{\new} \leftarrow \left(225\cdot \sqrt{\frac{n}{m}}\cdot \left(\wt{Q} + \wt{R} + \wt{T} + \frac{m-1}{n-1}\cdot \wt{H}\right)\right)^{-1} \in \R^{m \times m}$ \hspace{\fill} $\triangleright$ {$G \in \R^{m \times m}$}\\
%\label{line:G_init}\\
%\hspace{4mm} \hspace{\fill} $\triangleright$ {Maintain $G = \wt{H}^{-1}$ where $\wt{H} := \mathsf{A} \cdot (\wt{S}^{-1} \otimes \wt{S}^{-1}) \cdot \mathsf{A}^{\top}$}\\
%\hspace{4mm} {\bf Return} $(S,\wt{S},y,H,Q,R,T,\eta,\A,G)$\\
\hspace{4mm} {\bf return} $(\wt{Q},\wt{H},G^{\new})$ \\%\left(225\cdot \sqrt{\frac{n}{m}}\cdot \left(\wt{Q} + \wt{R} + \wt{T} + \frac{m-1}{n-1}\cdot \wt{H}\right)\right)^{-1} \in \R^{m \times m}$\\
{\bf end procedure}
\end{algorithmic}
\end{algorithm}

\begin{lemma}\label{lem:vol_time}
In Algorithm~\ref{alg:sdp_vol}, the (amortized) cost/time for every iteration is
\begin{align*}
    {O}^\ast \left(  \left(\frac{n}{m}\right)^{\frac{1}{4}} \cdot (n^2m + m^\omega n^{1/4}) + m^2 n^\omega + m^{4} + m^2 \cdot n^{\omega - \frac{1}{2}}\cdot \left(\frac{n}{m}\right)^{\frac{1}{4}} \right).
\end{align*}
\end{lemma}

\begin{proof}
Consider iteration $t$. 
In Line~\ref{lin:vol_S_new}-\ref{lin:vol_update_variables} of Algorithm~\ref{alg:sdp_vol}, we update $S^\new,\wt{S}^\new \in \R^{n \times n}$ in $O(mn^2 + n^\omega)$ time.

In Line~\ref{lin:log_grad}-\ref{lin:vol_grad} in Algorithm~\ref{alg:sdp_vol_cont}, we can first compute $S^{-1}A_iS^{-1}A_j \in \R^{n \times n}$ and $S^{-1}A_i \in \R^{n \times n}$ for all $i,j \in [m]$, in $O(m^2 n^\omega)$ time. Notice
\begin{align*}
    (\nabla \phi_{\mathrm{vol}}(y))_i = &~ -\tr[H(S)^{-1} \cdot \A  (S^{-1}A_i S^{-1}\otimes S^{-1}) \mathsf{A}^\top]\\
    = &~ \sum_{k=1}^m \sum_{l=1}^m -H(S)^{-1}_{k,l}\cdot\tr [A_k S^{-1} A_l S^{-1} A_i S^{-1}].
\end{align*}
Then it takes $O(m^3n^2)$ to find $\tr [A_k S^{-1} A_l S^{-1} A_i S^{-1}]$ for all $i,k,l \in [m]$ and subsequently $O(m^3)$-time to compute $\nabla \phi_{\mathrm{vol}}(y)$. Hence, the cost of Line~\ref{lin:log_grad}-\ref{lin:vol_grad} in Algorithm~\ref{alg:sdp_vol_cont} is $O(m^2 n^\omega + m^3 n^2 + m^3)$. 

For Line~\ref{lin:vol_H}-\ref{lin:vol_Q} of Algorithm~\ref{alg:sdp_vol_cont}, we first compute
\begin{align*}
    A_i,S^{-1}A_jV_3,~A_iV_3,,~V_4^\top A_iS^{-1}A_jS^{-1}, ~S^{-1}A_i V_3, ~V_4^\top A_j, ~V_4^\top A_i V_3
\end{align*}
for every $i \in \{1,\cdots,m\}$, for every $j \in \{1,\cdots,m\}$, in $O(m^2 \cdot \Tmat(n,n,r_t))$ time. 
Notice
\begin{align*}
    \wt{H}_{i,j} = &~ \tr[\wt{S}^{-1}A_i\wt{S}^{-1}A_j]\\
    = &~ H_{i,j} +  \tr[{S}^{-1}A_i V_3V_4^\top A_j] + \tr[{S}^{-1}A_j V_3V_4^\top A_i] + \tr[V_3V_4^\top A_i V_3V_4^\top A_j].
\end{align*}
Hence, it takes $\Tmat(m,nr_t,m)$ to compute $\tr[{S}^{-1}A_i V_3V_4^\top A_j], \tr[V_3V_4^\top A_i V_3V_4^\top A_j]$ for all $i,j \in [m]$ (by batching them together and using fast matrix multiplication on a $m$-by-$nr_t$ matrix and a $nr_t$-by-$m$ matrix) and subsequently $O(m^2)$-time to compute $\wt{Q}$. So the cost of Line~\ref{lin:vol_H} is $\Tmat(m,nr_t,m) + m^2$. 
For $\wt{Q}$, notice that for any $i,j\in [m]$,
\begin{align*}
    \wt{Q}_{i,j}
    = &~ \tr[H(S)^{-1} \A (S^{-1}A_i \wt{S}^{-1}A_jS^{-1} \otimes_S S^{-1}) \A^\top]\\
    = &~ Q_{i,j} + \sum_{k=1}^m \sum_{l=1}^m -H(S)^{-1}_{k,l}\cdot \frac{1}{2}\left(\tr [A_k S^{-1} A_l V_3V_4^\top A_i S^{-1} A_j S^{-1}]\right.\\
    &~\left.+ \tr[A_kS^{-1}A_i  V_3V_4^\top A_jS^{-1}A_\ell S^{-1}] \right).
\end{align*}
Then it takes $\Tmat(m^2,nr_t,m^2)$ to compute $\tr [A_k S^{-1} A_l V_3V_4^\top A_i S^{-1} A_j S^{-1}]$ for all $i,j,k,l \in [m]$ (by batching them together and using fast matrix multiplication on a $m^2$-by-$nr_t$ matrix and a $nr_t$-by-$m^2$ matrix) and subsequently $O(m^4)$-time to compute $\wt{Q}$. Hence the cost of Line~\ref{lin:vol_Q} is $\Tmat(m^2,nr_t,m^2) + m^4$.

In total, Line~\ref{lin:vol_H}-\ref{lin:vol_Q} of Algorithm~\ref{alg:sdp_vol_cont} takes 
\begin{align*}
    \Tmat(m^2,nr_t,m^2) + m^2 n^\omega + m^3n^2 + m^{4} +m^2 \cdot \Tmat(n,n,r_t).
\end{align*}

Summing up, the total cost in iteration $t$ is therefore given by
\begin{align*}
     \Tmat(m^2,nr_t,m^2) + m^2 n^\omega + m^3n^2 + m^{4} +m^2 \cdot \Tmat(n,n,r_t).
\end{align*}

Using Theorem~\ref{thm:general_amortize_guarantee_vol} with Lemma~\ref{lem:invariant_property_vol}, and Fact~\ref{fact:approx_Tmat_r_leq_n},
\begin{align*}
    \sum_{t=1}^T \Tmat(n,n,r_t) \leq &~ \sum_{t=1}^T {O}^\ast \left(n^2 + r_t^{\frac{\omega-2}{1-\alpha}}\cdot n^{2-\frac{\alpha(\omega-2)}{1-\alpha}} \right)\\
    \leq &~ {O}^\ast  \left( T \cdot n^2 + T \cdot n^{\frac{\omega-2}{1-\alpha} - 1/2} \cdot n^{2-\frac{\alpha(\omega-2)}{1-\alpha}} \cdot (n/m)^{\frac{1}{4}} \right)\\
    \leq &~ {O}^\ast  \left( T\cdot (n^2 + n^{\omega - 1/2} \cdot (n/m)^{\frac{1}{4}})\right).
\end{align*}

Using Corollary~\ref{cor:tmat_m2_nr_m2},
\begin{align*}
    \sum_{t=1}^T \Tmat(m^2,nr_t,m^2) \leq {O}^\ast \left( T \cdot  (n/m)^{1/4} \cdot  (n^2m + m^\omega n^{1/4})\right).
\end{align*}
Now, let us compute the (amortized) cost per iteration 
\begin{align*}
    &~ {O}^\ast \left(  \left(\frac{n}{m}\right)^{\frac{1}{4}} \cdot  (n^2m + m^\omega n^{1/4}) + m^2 n^\omega + m^{4} + m^2 \cdot \left(n^2 + n^{\omega - \frac{1}{2}}\cdot \left(\frac{n}{m}\right)^{\frac{1}{4}}\right) \right)\\
    = &~ {O}^\ast \left(  \left(\frac{n}{m}\right)^{\frac{1}{4}} \cdot (n^2m + m^\omega n^{1/4}) + m^2 n^\omega + m^{4} + m^2 \cdot n^{\omega - \frac{1}{2}}\cdot \left(\frac{n}{m}\right)^{\frac{1}{4}} \right).
\end{align*}
\end{proof}

\subsection{Approximation to \texorpdfstring{$Q$}{Q}}\label{sec:hybrid_approxQ}

The following lemma shows that $\wt{Q}$ in our algorithm is a good PSD approximation of $Q(S)$.
\begin{lemma}\label{lem:vol_approx_Q}
Let $\wt{Q} \in \R^{m \times m}$ be given by
\begin{align*}
    \wt{Q}_{i,j} = \tr[H(S)^{-1} \A (S^{-1}A_i\wt{S}^{-1}A_jS^{-1} \otimes_S S^{-1}) \A^\top].
\end{align*}
Then suppose $(1+\epsilon_S)^{-1} S \preceq \wt{S} \preceq (1+ \epsilon_S) S$ for $\epsilon_S \in (0,0.001)$, we have $(1+\epsilon_S)^{-3} Q \preceq \wt{Q} \preceq (1+ \epsilon_S)^3 Q$ where $Q := Q(S)$.
\end{lemma}

\begin{proof}
Fix $v \in \R^m$. We have
\begin{align*}
    v^\top \wt{Q} v = &~ \sum_{i=1}^m\sum_{j=1}^m v_iv_j\tr[H(S)^{-1} \A (S^{-1}A_i\wt{S}^{-1}A_jS^{-1} \otimes_S S^{-1}) \A^\top]\\
    = &~ \tr\left[ \A^\top H(S)^{-1} \A \cdot \left(S^{-1}(\sum_{i=1}^m v_iA_i )\wt{S}^{-1}(\sum_{j=1}^m v_jA_j)S^{-1} \otimes_S S^{-1}\right)\right]\\
    = &~ \frac{1}{2} \cdot (\tr[ \A^\top H^{-1} \A \cdot (S^{-1}B\wt{S}^{-1}BS^{-1} \otimes S^{-1})] + \tr[ \A^\top H^{-1} \A \cdot (S \otimes S^{-1}B\wt{S}^{-1}BS^{-1})])
\end{align*}
where we abbreviate $B = \sum_{i=1}^m v_iA_i$ and $H = H(S)$. 

For $\tr[ \A^\top H^{-1} \A \cdot (S^{-1}B\wt{S}^{-1}BS^{-1} \otimes S^{-1})]$, we have
\begin{align*}
    &~ \tr[ \A^\top H^{-1} \A \cdot (S^{-1}B\wt{S}^{-1}BS^{-1} \otimes S^{-1})]\\
    = &~ \tr[ (S^{-1}B\wt{S}^{-1/2} \otimes I)^\top \A^\top H^{-1} \A (S^{-1}B\wt{S}^{-1/2} \otimes I) (I \otimes S^{-1})]\\
    \leq &~ (1+\epsilon_S)\cdot \tr[ (S^{-1}B\wt{S}^{-1/2} \otimes I)^\top \A^\top H^{-1} \A (S^{-1}B\wt{S}^{-1/2} \otimes I) (I \otimes \wt{S}^{-1})]\\
    = &~ (1+\epsilon_S)\cdot \tr[ (S^{-1}B \otimes I)^\top \A^\top H^{-1} \A (S^{-1}B \otimes I) (\wt{S}^{-1} \otimes \wt{S}^{-1})]\\
    \leq &~ (1+\epsilon_S)^3 \cdot \tr[ (S^{-1}B \otimes I)^\top \A^\top H^{-1} \A (S^{-1}B \otimes I) ({S}^{-1} \otimes {S}^{-1})]\\
    = &~ (1+\epsilon_S)^3 \cdot\tr[ \A^\top H^{-1} \A \cdot (S^{-1}B {S}^{-1}BS^{-1} \otimes S^{-1})],
\end{align*}
where we use Fact~\ref{fac:kron_basic} in first step; we use $I \otimes {S}^{-1} \preceq (1+ \epsilon) \cdot (I \otimes \wt{S}^{-1})$ and Fact~\ref{fac:trace_loewner} in second step; the third step comes from Fact~\ref{fac:kron_basic}; the fourth step comes from Fact~\ref{fac:kron_spectral}; the last step comes from  Fact~\ref{fac:kron_basic}.

Similarly, for $\tr[ \A^\top H^{-1} \A \cdot (S^{-1} \otimes S^{-1} B\wt{S}^{-1}BS^{-1} )]$, we have
\begin{align*}
    &~ \tr[ \A^\top H^{-1} \A \cdot (S^{-1} \otimes S^{-1} B\wt{S}^{-1}BS^{-1})]\\
    = &~ \tr[ (S^{-1/2} \otimes S^{-1}B)^\top \A^\top H^{-1} \A (S^{-1/2} \otimes S^{-1}B) (I \otimes \wt{S}^{-1})]\\
    \leq &~ (1+\epsilon_S)\cdot \tr[ (S^{-1/2} \otimes S^{-1}B)^\top \A^\top H^{-1} \A (S^{-1/2} \otimes S^{-1}B) (I \otimes {S}^{-1})]\\
    = &~ (1+\epsilon_S) \cdot \tr[ \A^\top H^{-1} \A \cdot (S^{-1} \otimes S^{-1} B {S}^{-1}BS^{-1} )],
\end{align*}
where we use Fact~\ref{fac:kron_basic} in the first step; we use $I \otimes {S}^{-1} \preceq (1+ \epsilon) \cdot (I \otimes \wt{S}^{-1})$ and Fact~\ref{fac:trace_loewner} in second step; the third step comes from Fact~\ref{fac:kron_basic}. 

Summing up,
\begin{align*}
    v^\top \wt{Q} v \leq &~ (1+\epsilon_S)^3 \cdot\tr[ \A^\top H^{-1} \A \cdot (S^{-1}B {S}^{-1}BS^{-1} \otimes S^{-1})]\\
    &~ + (1+\epsilon_S) \cdot \tr[ \A^\top H^{-1} \A \cdot (S^{-1} \otimes S^{-1} B {S}^{-1}BS^{-1} )] \\
    \leq &~ (1+\epsilon_S)^3 \cdot \tr[ \A^\top H^{-1} \A \cdot (S^{-1}B {S}^{-1}BS^{-1} \otimes_S S^{-1})]\\
    = &~ (1+\epsilon_S)^3 \cdot v^\top {Q} v.
\end{align*}
Similarly, $v^\top \wt{Q} v \geq (1+\epsilon_S)^{-3} \cdot v^\top {Q} v$. Therefore, $(1+\epsilon_S)^{-3} Q \preceq \wt{Q} \preceq (1+ \epsilon_S)^3 Q$, since $v$ can be arbitrarily chosen.
\end{proof}

\subsection{\texorpdfstring{$S$}{S} move in hybrid barrier}\label{sec:hybrid_S_slow}

\begin{lemma}[$S$ move in hybrid barrier]\label{lem:invariant_property_vol}
Consider Algorithm~\ref{alg:sdp_vol}-\ref{alg:sdp_vol_cont}, in each iteration, the following invariant holds: 
\begin{align}\label{eq:promise_sdp_vol}
\| S^{-1/2} S^{\new} S^{-1/2} - I \|_F \leq&~ 1.03 \cdot \epsilon_N \cdot (n/m)^{1/2},\\
\|S^{-1/2}S^{\new}S^{-1/2} - I\|_2^2 \leq &~ 0.002\cdot \epsilon_N.\notag
\end{align}
\end{lemma}

\begin{proof}

We note that the robust framework and all corresponding results in Section~\ref{sec:robust_central_path} directly applies to hybrid barrier with $\theta = (mn)^{1/2}$. We have
\begin{align*}
    \| S^{-1/2} S^{\new} S^{-1/2} - I \|^2_F = &~ \tr\left[ \left(S^{-1/2} (S^\new - S) S^{-1/2}\right)^2\right]\\
    = &~ \tr\left[ S^{-1} \big(\sum_{i\in [m]} \Tilde{\delta}_{y,i} A_i\big) S^{-1} \big(\sum_{i\in[m]} \Tilde{\delta}_{y,i} A_i\big)\right]\\
    = &~ \sum_{i\in[m]} \sum_{j\in[m]} \Tilde{\delta}_{y,i} \Tilde{\delta}_{y,j} \tr[S^{-1}A_i S^{-1}A_j]\\
    = &~ \Tilde{\delta}_{y}^\top H(y) \Tilde{\delta}_{y}\\
    = &~ \|\Tilde{\delta}_{y}\|^2_{H(y)}
\end{align*}
where the second step comes from $S^\new - S = \sum_{i=1}^m \Tilde{\delta}_{y,i} A_i$ and the rest follows from algebra. It suffices to bound $\|\Tilde{\delta}_{y}\|_{H(y)}$. 

Since $\phi(y)=\sqrt{\frac{n}{m}} \cdot \left(\phi_{\mathrm{vol}}(y) + \frac{m-1}{n-1}\cdot \phi_{\log}(y)\right)$, we have
\begin{align*}
    \nabla^2\phi(y) = (n/m)^{1/2}\nabla^2\phi_{\mathrm{vol}}(y)+(m/n)^{1/2}H(y).
\end{align*}
By Fact~\ref{fac:vol_Q_Sigma}, we have
\begin{align*}
    Q(S)\succeq \frac{1}{n}H(y).
\end{align*}
And we also have
\begin{align*}
     \nabla^2\phi(y) \succeq &~  (n/m)^{1/2} Q(S)+(m/n)^{1/2}H(y)\\
     \succeq &~ \left((n/m)^{1/2} \cdot n^{-1} + (m/n)^{1/2} \right)H(y)\\
     \succeq &~ O((m/n)^{1/2})H(y).
\end{align*}
Thus,
\begin{align*}
    \|\Tilde{\delta}_{y}\|^2_{H(y)} \leq (n/m)^{1/2}\cdot \|\Tilde{\delta}_{y}\|^2_{\nabla^2\phi(y)}.
\end{align*}
For $\|\Tilde{\delta}_{y}\|^2_{\nabla^2\phi(y)}$, we have
\begin{align*}
    \|\Tilde{\delta}_{y}\|^2_{\nabla^2\phi(y)} \leq &~ \|n(y,\eta^\new)\|_{\nabla^2\phi(y)} + \|\Tilde{\delta}_{y} - (- n(y,\eta^\new))\|_{\nabla^2\phi(y)}\\
    \leq &~ 1.01 \|n(y,\eta^\new)\|_{\nabla^2\phi(y)} \\ 
    \leq &~ 1.01 \cdot \left( (1+\epsilon_N/20)\|n(y,\eta)\|_{\nabla^2\phi(y)} + \epsilon/20 \right)\\
    \leq &~ 1.03 \cdot \epsilon_N
\end{align*}
where the first step comes from triangle inequality, the second step comes from Eq.~\eqref{eq:approx_newton_error}, the third step comes from Lemma~\ref{lem:eta_move} and the last step comes from choice of $\epsilon_N$. 
Hence,
\begin{align*}
    \|\Tilde{\delta}_{y}\|^2_{H(y)} \leq 1.03 \cdot \epsilon_N \cdot (n/m)^{1/2},
\end{align*}
that is,
\begin{align*}
    \| S^{-1/2} S^{\new} S^{-1/2} - I \|_F \leq 1.03 \cdot \epsilon_N \cdot (n/m)^{1/2}.
\end{align*}

Moreover, by Fact~\ref{fac:hybrid_lb}, we have
\begin{align*}
    \|\Tilde{\delta}_{y}\|^2_{\nabla^2\phi(y)} \geq &~  225(n/m)^{1/2}\cdot \frac{2\sqrt{m}}{1+\sqrt{n}}\cdot \left\|S^{-1/2}\left(\sum_{i=1}^m \Tilde{\delta}_{y,i}A_i\right) S^{-1/2}\right\|_2^2\\
    = &~ O(1)\cdot \|S^{-1/2}S^{\new}S^{-1/2} - I\|_2^2.
\end{align*}
Hence, we have 
\begin{align*}
    \|S^{-1/2}S^{\new}S^{-1/2} - I\|_2^2 \leq 0.002\cdot \epsilon_N.
\end{align*}
This completes the proof of Eq.~\eqref{eq:promise_sdp_vol}.
\end{proof}

\subsection{Property of low rank update for the hybrid barrier}\label{sec:hybrid_amortize}

\begin{assumption}[Closeness of $S^{\new}$ and $\wt{S}$ from $S$]\label{ass:close_S_vol}
We make the following two assumptions about $S, \wt{S}, S^{\new} \in \R^{n \times n}$:
\begin{align*}
    1. &~ \|S^{-1/2} \cdot S^{\new} \cdot S^{-1/2} - I\|_F \leq 0.02(n/m)^{1/2}, \\
    2. &~ \|S^{-1/2} \cdot S^{\new} \cdot S^{-1/2} - I\|_2 \leq 0.005,\\
    3. &~ \|S^{-1/2} \cdot \wt{S} \cdot S^{-1/2} - I\|_2 \leq 0.01.
\end{align*}
\end{assumption}

\begin{theorem}[General amortized guarantee for the hybrid barrier]\label{thm:general_amortize_guarantee_vol}
Let $T$ denote the total number of iterations in Algorithm~\ref{alg:sdp_vol}-\ref{alg:sdp_vol_cont}. Let $r_t$ denote the rank of the update matrices $V_1, V_2 \in \R^{n \times r_t}$ generated by Algorithm~\ref{alg:low_rank_slack_update} in the $t$-th iteration.
Suppose Assumption~\ref{ass:close_S_vol} hold. The ranks $r_t$'s satisfy the following condition: for any vector $g \in \R_+^n$ which is non-increasing, we have
\begin{align}\label{eq:amortization_vol}
    \sum_{t=1}^T r_t \cdot g_{r_t} \leq O(T \cdot (n/m)^{1/4}\cdot \|g\|_2 \cdot \log n).
\end{align}
\end{theorem}

The proof of Theorem~\ref{thm:general_amortize_guarantee_vol} relies on the following three lemmas:

\begin{lemma}[Variant of Lemma~\ref{lem:part_of_lem_62_jkl}]
\label{lem:part_of_lem_62_jkl_vol}
Let matrices $Z, Z^{\midd} \in \R^{n \times n}$ be defined as in Definition~\ref{def:diff_matrix}. Under Assumption~\ref{ass:close_S_vol}, we have
\begin{align*}
    \sum_{i=1}^n (\lambda(Z)_{[i]} - \lambda(Z^{\midd})_{[i]})^2 \leq 10^{-3} (n/m)^{1/2},
\end{align*}
where $\lambda(Z)_{[i]}$ denotes the $i$-th largest eigenvalue of $Z$.
\end{lemma}
\begin{proof}[Proof sketch]
The proof is very similar to Lemma~\ref{lem:part_of_lem_62_jkl} except the upper bound of the following quantity:
\begin{align*}
    \|S^{1/2}(S^{\new})^{-1}S^{1/2}- I\|_F^2= \sum_{i=1}^n (\nu_i^{-1}-1)^2,
\end{align*}
where $\{\nu_i\}_{i\in[n]}$ are the eigenvalues of $S^{-1/2}S^{\new}S^{-1/2}$. Then, by Assumption~\ref{ass:close_S_vol} part 2, we have
\begin{align*}
    \max_{i\in [n]} |v_i - 1| \leq 0.005,
\end{align*}
which implies that $v_i\geq 0.995$ for all $i\in [n]$. By Assumption~\ref{ass:close_S_vol}, 
\begin{align*}
    \sum_{i=1}^n (v_i - 1)^2 \leq 0.02 (n/m)^{1/2}.
\end{align*}
Then, it follows that
\begin{align*}
    \sum_{i=1}^n (\nu_i^{-1}-1)^2 \leq 5\times 10^{-4} \cdot (n/m)^{1/2}.
\end{align*}
The remaining part does not change.
\end{proof}

\begin{lemma}[$S$ move]\label{lem:S_move_vol}
Consider the $t$-th iteration. Let matrices $Z,Z^{\midd} \in \R^{n \times n}$ be defined as in Definition~\ref{def:diff_matrix}. Let $g \in \R_+^n$ be a non-increasing vector, and let $\Phi_g: \R^{n \times n} \to \R_+$ be defined as in Definition~\ref{def:potential}.

Under Assumption~\ref{ass:close_S_vol}, we have
\begin{align*}
    \Phi_g(Z^{\midd}) - \Phi_g(Z) \leq \|g\|_2\cdot (n/m)^{1/4}.
\end{align*}
\end{lemma}
\begin{proof}[Proof sketch]
The proof is basically the same as Lemma~\ref{lem:S_move}. We can upper bound the LHS as follows:
\begin{align*}
    \Phi_g(Z^{\midd}) - \Phi_g(Z) \leq &~ \|g\|_2 \cdot  \Big( \sum_{i=1}^n |\lambda(Z^{\midd})_{\pi(i)} - \lambda(Z)_{\pi(i)}|^2 \Big)^{1/2}\\
    \leq &~ \|g\|_2 \cdot (n/m)^{1/4},
\end{align*}
where the last step follows from Lemma~\ref{lem:part_of_lem_62_jkl_vol}.
\end{proof}

\begin{lemma}[$\wt{S}$ move]\label{lem:wt_S_move_vol}
Consider the $t$-th iteration. Let matrices $Z^{\midd}, Z^{\new} \in \R^{n \times n}$ be defined as in Definition~\ref{def:diff_matrix}. Let $g \in \R_+^n$ be a non-increasing vector, and let $\Phi_g: \R^{n \times n} \to \R_+$ be defined as in Definition~\ref{def:potential}.

Let $r_t$ denote the rank of the update matrices $V_1, V_2 \in \R^{n \times r_t}$ generated by Algorithm~\ref{alg:low_rank_slack_update} in the $t$-th iteration. We have
\begin{align*}
    \Phi_g(Z^{\midd}) - \Phi_g(Z^{\new}) \geq \frac{\epsilon_S}{10 \log n} \cdot r_t \cdot g_{r_t}.
\end{align*}
\end{lemma}
The proof is exactly the same as Lemma~\ref{lem:wt_S_move}.
Now, we are ready to prove Theorem~\ref{thm:general_amortize_guarantee_vol}.
\begin{proof}[Proof of Theorem~\ref{thm:general_amortize_guarantee_vol}]
Combining Lemma~\ref{lem:S_move_vol} and Lemma~\ref{lem:wt_S_move_vol}, we have
\begin{align*}
    \Phi_g(Z^{\new}) - \Phi_g(Z) 
    = & ~ (\Phi_g(Z^{\midd}) - \Phi_g(Z)) - (\Phi_g(Z^{\midd}) - \Phi_g(Z^{\new})) \\
    \leq & ~ (n/m)^{1/4}\cdot \|g\|_2 - \frac{\epsilon_S}{10 \log n} \cdot r_t \cdot g_{r_t}.
\end{align*}
With an abuse of notation, we denote the matrix $Z$ in the $t$-th iteration as $Z^{(t)}$. Since in the beginning $\Phi_g(Z^{(0)}) = 0$ and $\Phi_g(Z^{(T)}) \geq 0$, we have
\begin{align*}
    0 \leq & ~ \Phi_g(Z^{(T)}) - \Phi_g(Z^{(0)}) \\
    \leq & ~ \sum_{t=1}^T (\Phi_g(Z^{(t)}) - \Phi_g(Z^{(t-1)})) \\
    \leq & ~ T\cdot (n/m)^{1/4}\cdot \|g\|_2 - \frac{\epsilon_S}{10 \log n} \cdot \sum_{t=1}^T r_t \cdot g_{r_t}.
\end{align*}
This completes the proof.
\end{proof}

\begin{corollary}\label{cor:tmat_m2_nr_m2}
Given a sequence $r_1,\dots,r_T \in [0,n]$ that satisfies Eq.~\eqref{eq:amortization_vol}. We have
\begin{align*}
    \sum_{t=1}^T \Tmat (m^2,n r_t,m^2) \leq {O}^\ast \left( T \cdot  (n/m)^{1/4} \cdot  (n^2m + m^\omega n^{1/4}) \right).
\end{align*}
\end{corollary}

\begin{proof}
Let $a_t = \log_n (r_t)$ and $b = \log_n (m^2)$. Then 
\begin{align*}
    \Tmat (m^2,n r_t,m^2) = \Tmat(n^b,n^{1+a_t},n^b) = {O}^\ast(n^{b \cdot \omega((1+a_t)/b)}).
\end{align*}
For each $i \in \{0,1,\dots,\log n) \}$, define 
\begin{align*}
    T_i = \{t \in [T]: 2^i \leq r_t \leq 2^{i+1} \}.
\end{align*}
Let $g_r = r^{-1/2}$, Theorem~\ref{thm:general_amortize_guarantee_vol} indicates
\begin{align*}
    \sum_{i=1}^{\log n} |T_i| \cdot 2^{i/2} \leq \sum_{t=1}^T r_t^{1/2}\leq O( (n/m)^{1/4} \cdot T \cdot \log^{1.5}n).
\end{align*}
This implies $|T_i| \leq O((n/m)^{1/4} \cdot T \cdot \log^{1.5}(n) / 2^{i/2}$. It thus follows that
\begin{align*}
    \sum_{t=1}^T \Tmat (m^2,n r_t,m^2) \leq &~ {O}^\ast \left(\sum_{t=1}^T n^{b \cdot \omega((1+a_t)/b)}\right)\\
    = &~ {O}^\ast \left(\sum_{i=1}^{\log n} \sum_{t \in T_i} n^{b \cdot \omega((1+a_t)/b)}\right)\\
    \leq &~ {O}^\ast \left( \max_{i \in \log n} \max_{t \in T_i} \frac{(n/m)^{1/4} \cdot T}{2^{i/2}} \cdot n^{b \cdot \omega((1+a_t)/b)}\right)\\
    %\leq &~ \title{O}(1) \cdot T \cdot \max_{i \in \log n} \max_{2^i \leq n^{a_t} \leq 2^{i+1}} \frac{1}{2^{i/2}} \cdot n^{b \cdot \omega((1+a_t)/b)}\\
    \leq &~ {O}^\ast \left( (n/m)^{1/4} \cdot T \cdot  \max_{ {a_t} \in [0,1]} n^{-a_t/2 + b \cdot \omega((1+a_t)/b)}\right).
\end{align*}
Since $\omega(\cdot)$ is a convex function (Fact~\ref{fact:omega_k_convex}), 
\begin{align*}
    \max_{ {a_t} \in [0,1]} -a_t/2 + b \cdot \omega((1+a_t)/b) \leq &~ \max_{a \in \{0,1\}}-a/2 + b \cdot \omega((1+a)/b)\\
    \leq &~ \max\{b+2, b\omega + 0.25\}.
\end{align*}

Combining the above inequalities, we have
\begin{align*}
    \sum_{t=1}^T \Tmat (m^2,n r_t,m^2) \leq &~ {O}^\ast \left((n/m)^{1/4} \cdot  T \cdot  n^{\max\{b+2, b\omega + 0.25\}}\right)\\
    \leq &~ {O}^\ast \left( T \cdot  (n/m)^{1/4} \cdot  (n^2m + m^\omega n^{1/4}) \right).
\end{align*}
This completes the proof.
\end{proof}

\subsection{Our result}\label{sec:hybrid_result}

\begin{theorem}[Main result for Algorithm~\ref{alg:sdp_vol} - \ref{alg:sdp_vol_cont}]\label{thm:vol_formal}
Given symmetric matrices $C, A_1,\cdots,A_m \in \R^{n \times n}$, and a vector $b \in \R^m$. Define matrix $\mathsf{A} \in \R^{m \times n^2}$ by stacking the $m$ vectors $\vect[A_1], \cdots, \vect[A_m] \in \R^{n^2}$ as rows. Consider the following SDP instance:
\begin{align*}
     \max_{X \in \R^{n \times n}} & ~ \langle C,X \rangle \\
    \mathrm{~s.t.~} &~ \langle A_i, X \rangle = b_i, ~ \forall i \in [m], \\
    &~ X \succeq 0,
\end{align*}
Let $X^*$ be an optimal solution of the SDP instance. 
There is a SDP algorithm (Algorithm~\ref{alg:sdp_vol}-\ref{alg:sdp_vol_cont}) that runs 
in time
\begin{align*}
    {O}^\ast \left( (mn)^{1/4}\cdot \left( m^{2} n^{\omega} +  m^4 \right) \cdot \log (1/\epsilon) \right) %m^{3}n^{2}
\end{align*}
and outputs a PSD matrix $X \in \mathbb{R}^{n \times n}$ s.t.
\begin{align*}
\langle C, X \rangle \geq \langle C, X^* \rangle - \epsilon \cdot \| C \|_{2} \cdot R \quad \text{and} \quad \sum_{i = 1}^m\left| \langle A_i, X \rangle - b_i \right| \leq 4 n \epsilon \cdot \Big( R  \sum_{i = 1}^m \| A_i \|_1 + \| b \|_1 \Big) ,
\end{align*}
 %, and $\| A_i \|_1$ is the Schatten $1$-norm of matrix $A_i$.
\end{theorem}

\begin{remark}
We improve the running time of \cite{a00} 
\begin{align*} 
O^\ast( (mn)^{1/4} \cdot (m^3 n^\omega + m^4n^2+m^{\omega+2}) \cdot \log(1/\epsilon)))
\end{align*}
for all parameters regime.

In particular, if $m = n$, the total cost of Algorithm~\ref{alg:sdp_vol}-\ref{alg:sdp_vol_cont} can be upper bounded by $n^{\omega+2.5}$. This improves the $n^{6.5}$ total cost \cite{a00}.

If $m = n^2$, this cost can be upper bounded by $n^{8.75}$. This improves the $n^{10.75}$ total cost of \cite{a00}.
\end{remark}

\begin{proof}
We first compute the running time. From Lemma~\ref{lem:vol_time}, the amortized cost per iteration is upper bounded by 
\begin{align*}
    {O}^\ast \left(  \left(\frac{n}{m}\right)^{\frac{1}{4}} \cdot (n^2m + m^\omega n^{1/4}) + m^2 n^\omega + m^{4} + m^2 \cdot n^{\omega - \frac{1}{2}}\cdot \left(\frac{n}{m}\right)^{\frac{1}{4}} \right).
\end{align*}
Since $T = O\left((mn)^{1/4}\cdot \log (mn/\epsilon)\right)$, the Algorithm~\ref{alg:sdp_vol}-\ref{alg:sdp_vol_cont} run 
in time
\begin{align*}
    {O}^\ast \left( \left(n^{\frac{5}{2}}m + m^{2+\frac{1}{4}} n^{\omega+\frac{1}{4}} + m^{4+\frac{1}{4}}n^{\frac{1}{4}} +m^\omega n^{\frac{3}{4}} \right) \cdot \log (1/\epsilon) \right).
\end{align*}
Under current matrix multiplication exponent $\omega \approx 2.373$ and $0 \leq m \leq n^2$, this simplifies to 
\begin{align*}
    & ~ {O}^\ast \left( \left( m^{2+\frac{1}{4}} n^{\omega+\frac{1}{4}} + m^{4+\frac{1}{4}}n^{\frac{1}{4}} \right) \cdot \log (1/\epsilon) \right) \\
    = & ~ {O}^\ast \left( (mn)^{1/4} \cdot  \left( m^{2 } n^{\omega } + m^{4 }  \right) \cdot \log (1/\epsilon) \right) 
\end{align*}

Now we prove the correctness of Algorithm~\ref{alg:sdp_vol}-\ref{alg:sdp_vol_cont}. We invoke the robust framework. From Fact~\ref{fac:vol_Q_Sigma} and Lemma~\ref{lem:vol_approx_Q}, $\epsilon_g = \epsilon_\delta = 0$, and $c_H = 1/3 \cdot 1.0001^{-3}$. Therefore directly applying Theorem~\ref{thm:approx_central_path} for $\theta = (mn)^{1/4}$ completes the proof.
\end{proof}

\newpage
\onecolumn
\appendix

\section*{Appendix}
\section{Initialization}\label{sec:init}
Let us state an initialization result which is very standard in literature, see Section 10 in \cite{lsw15}, Appendix A in \cite{cls19}, and Section 9 in \cite{jklps20}. It can be easily proved using the property of special matrix/Kronecker product (Fact~\ref{fac:transform_norm_to_matrix_norm}).
\begin{lemma}
\label{lem:initialization}
For an SDP instance defined in Definition~\ref{def:sdp_primal} ($m$ $n \times n$  constraint matrices, let $X^*$ be any optimal solution to SDP.), assume it has two properties : 
\begin{enumerate}
\item Bounded diameter: for any feasible solution $X \in \R^{n\times n}_{\succeq 0}$, it has $\|X\|_2 \leq R$.
\item Lipschitz objective: the objective matrix $C\in \R^{n\times n}$ has bounded spectral norm, i.e., $\|C\|_2 \leq L$. 
\end{enumerate}
Given $\epsilon \in (0, 1/2]$, we can construct the following {\em modified} SDP instance in dimension $n+2$ with $m+1$ constraints:
\begin{align*}
&\max_{\overline{X} \succeq 0} ~  \langle \overline{C}, \overline{X} \rangle \\
\mathrm{s.t.} ~ & \langle \overline{A}_i, \overline{X} \rangle = \overline{b}_i, ~\forall i \in [m + 1], 
\end{align*}
where 
\begin{align*}
\overline{A}_i =
\left[ 
\begin{matrix}  
A_i & 0_n & 0_{n} \\
0_n^\top & 0 & 0 \\
0_{n}^\top & 0 & \frac{b_i}{R} - \tr[A_i] 
\end{matrix} 
\right]~~\forall i \in [m],~\text{and}~~
\overline{A}_{m+1} &=
\left[ 
\begin{matrix}  
I_n & 0_n & 0_{n} \\
0_n^\top & 1 & 0 \\
0_{n}^\top & 0 & 0
\end{matrix} 
\right].
\end{align*}
\begin{align*}
\overline{b} = 
\left[
\begin{matrix}
\frac{1}{R} b\\
n + 1
\end{matrix}
\right] ~,~
\overline{C} = \left [
\begin{matrix}  
\frac{\epsilon}{L} \cdot C & 0_n & 0_{n} \\
0_n^\top & 0 & 0 \\
0_{n}^\top & 0 & -1 
\end{matrix} \right].
\end{align*}
Moreover, it has three properties: 
\begin{enumerate}
\item $(\overline{X}_0, \overline{y}_0, \overline{S}_0)$ are feasible primal and dual solutions of the modified instance, where
\begin{align}\label{eq:initial_modified_SDP}
\overline{X}_0 = I_{n+2} ~,~ \overline{y}_0 = \left[ \begin{matrix} 0_m \\ 1 \end{matrix} \right] ~,~
\overline{S}_0=
\left[ \begin{matrix}
I_n - C \cdot \frac{\epsilon }{L} & 0_n & 0 \\
0_n^\top & 1 & 0 \\
0_n^\top & 0 & 1 
\end{matrix}\right] .
\end{align}

\item For any feasible primal and dual solutions $(\overline{X}, \overline{y}, \overline{S})$ with duality gap at most $\epsilon^2$, the matrix $\widehat{X} = R \cdot \overline{X}_{[n] \times [n]}$, where $\overline{X}_{[n] \times [n]}$ is the top-left $n$-by-$n$ block submatrix of $\overline{X}$. The matrix $\widehat{X}$ has three properties % is an approximate solution to the original SDP instance in the following sense:
\begin{align*}
\langle C, \widehat{X} \rangle &\geq \langle C, X^* \rangle - L R \cdot \epsilon ,\\
\widehat{X} &\succeq 0, \\
\sum_{i \in [m] } |\langle A_i, \widehat{X} \rangle - b_i | &\leq 4 n \epsilon \cdot \Big(R \sum_{i \in [m] } \|A_i\|_1 + \|b\|_1 \Big), 
\end{align*}

\item If we take $\epsilon \leq \epsilon_N^2$ in Eq.~\eqref{eq:initial_modified_SDP},
then the initial dual solution satisfies the induction invariant: $$g(\overline{y}_0,\eta)H(\overline{y}_0)^{-1}g(\overline{y}_0,\eta) \leq \epsilon_N^2.$$
\end{enumerate}
\end{lemma}

\section{From Dual to Primal}
We state a lemma about transforming a nearly optimal dual solution to a primal solution for SDP, which is very standard in literature and follows directly from Section 10 in \cite{lsw15} and Fact~\ref{fac:transform_norm_to_matrix_norm}.
\begin{lemma}\label{lem:dual_2_primal}
Given parameter $\eta_{\mathrm{final}} = \frac{1}{n+2}  (1+\frac{\epsilon_N}{20\sqrt{n}})^T$ where $T = 40 \epsilon_N^{-1} \sqrt{n} \log(n/\epsilon)$, dual variable $y \in \R^m$ and slack variable $S\in \R^{n \times n}$ such that 
\begin{align*}
    g(y,\eta_{\mathrm{final}})H(y)^{-1}g(y,\eta_{\mathrm{final}}) \leq &~ \epsilon_N^2\\
    \sum_{i=1}^m y_i A_i - C = &~ S \\
    b^\top y \leq & ~ b^\top y^* + \frac{n}{\eta_{\mathrm{final}}} (1+2\epsilon_N) \\
    S \succ & ~ 0 
\end{align*}
with $\epsilon_N \leq 1/10$.
Then there is an algorithm that finds a primal variable $X \in \R^{n \times n}$ in $O(n^{\omega + o(1)})$ time such that
\begin{align}\label{eq:dual_bad_1}
\langle C , X \rangle \geq & ~ \langle C , X^* \rangle - LR \cdot \epsilon \notag\\
X \succeq & ~ 0  \\
\sum_{i \in [m]}| \langle A_i, X \rangle - b_i | \leq & ~ 4 n \epsilon \cdot (R \sum_{i\in [m]} \| A_i \|_1 + \| b \|_1) \notag %\\
\end{align}
\end{lemma}

\section{Our Straightforward Implementation of the Hybrid Barrier SDP Solver}\label{sec:straightforward_vol}
\begin{theorem}[Our straight forward implementation of the hybrid algorithm \cite{a00}]\label{thm:our_a00}
The original hybrid barrier algorithm \cite{a00} use
\begin{align*}
    O^\ast(m^2 n^{\omega} + m^4 n^2 + m^{\omega+2})
\end{align*}
cost per iteration. 
\end{theorem}
\begin{remark}\label{rem:our_a00}
A naive implementation of hybrid barrier (e.g. Table 1.1 of \cite{jklps20}\footnote{The bound claimed in \cite{jklps20} is $O^\ast(m^3 n^\omega + m^4n^2)$, since they want to consider the special parameter regime where $n\leq m \leq n^2$. In that regime, $m^{\omega+2}$ is dominated by the first two terms.}) takes time
\begin{align*}
   O^\ast(m^3 n^\omega + m^4n^2+m^{\omega+2}).
\end{align*}
Our implementation (Theorem~\ref{thm:our_a00}) improves it to $m^2 n^{\omega} + m^4 n^2 + m^{\omega+2}$ by reusing the computations in the Hessian matrix. 
%}
\end{remark}

\begin{proof}
In each iteration, the computation workload is comprised of the following:
\begin{itemize}
\item The slack variable $S \in \R^{n \times n}$, given by
\begin{align*}
    S = S(y) = \sum_{i=1}^m y_i A_i - C.
\end{align*}

\item The gradient of $\phi_{\mathrm{vol}}$, denote by $\nabla \phi_{\mathrm{vol}}(y) \in \R^m$, given by:
\begin{align*}
    \nabla \phi_{\mathrm{vol}}(y)_i = -\tr[H(S)^{-1} \cdot \A  (S^{-1}A_i S^{-1}\otimes S^{-1}) \mathsf{A}^\top] .
\end{align*}
 
\item The gradient of $\phi_{\mathrm{log}}$, denote by $\nabla \phi_{\mathrm{log}}(y) \in \R^m$, given by:
\begin{align*}
    \nabla \phi_{\mathrm{log}}(y)_i = -\tr [S^{-1} \cdot A_i].
\end{align*}

\item The Hessian matrix of $ \phi_{\mathrm{log}}$, denoted by $H(S) \in \R^{m \times m }$, given by:
\begin{align*}
    H(S) = \mathsf{A} \cdot ({S}^{-1} \otimes {S}^{-1} ) \cdot \mathsf{A}^{\top}
\end{align*}

\item The first component of the Hessian matrix of $\phi_{\mathrm{vol}}$, denoted by $Q(S) \in \R^{m \times m }$, (recall from \cite{a00}, $\nabla^2 \phi_{\mathrm{vol}}(y) = 2Q(S) + R(S) -2T(S) \in \R^{m \times m}$), given by:
\begin{align*}
    Q(S)_{i,j} = \tr[H(S)^{-1} \A (S^{-1}A_iS^{-1}A_jS^{-1} \otimes_S S^{-1}) \A^\top]. 
\end{align*}

\item The second component of Hessian matrix of $\phi_{\mathrm{vol}}$, denoted by $R(S) \in \R^{m \times m }$, given by:
\begin{align*}
    R(S)_{i,j} = \tr[H(S)^{-1} \A (S^{-1}A_iS^{-1} \otimes_S S^{-1}A_jS^{-1}) \A^\top].
\end{align*}

\item The third component of Hessian matrix of $\phi_{\mathrm{vol}}$, denoted by $T(S) \in \R^{m \times m }$, given by:
\begin{align*}
    T(S)_{i,j} = \tr[H(S)^{-1} \A (S^{-1}A_iS^{-1} \otimes_S S^{-1}) \A^\top H(S)^{-1}\A (S^{-1}A_jS^{-1} \otimes_S S^{-1}) \A^\top ].
\end{align*}

\item Newton direction, denoted by $\delta_y \in \R^{m}$, given by
\begin{align*}
    \delta_y = - (\nabla^2 \phi(y))^{-1} (\eta b - \nabla \phi(y)). 
\end{align*}
\end{itemize}

To find all these items, we carry out the following computations. 

\paragraph{Step 1.}
We compute $S \in \R^{n \times n}$ in $mn^2$ time.

\paragraph{Step 2.}
We compute $S^{-1}A_iS^{-1}A_j \in \R^{n \times n}$ and $S^{-1}A_i  $ for all $i \in \{1,\cdots,m\}$, $\forall j \in \{1,\cdots,m\}$. This step costs $n^\omega m^2$. Since $H(S)_{i,j} = \tr[S^{-1}A_iS^{-1}A_j]$, it costs an additional $n m^2$ time to compute $H(S)$, and $n^\omega$ time to find $H(S)^{-1}$ correspondingly. Since $\nabla \phi_{\mathrm{log}}(y)_i = -\tr [S^{-1} \cdot A_j]$, we already find $\nabla \phi_{\mathrm{log}}(y)$. In total, this step costs $n^\omega m^2$ time.

\paragraph{Step 3.}
We compute $\tr [S^{-1} A_k S^{-1} A_l S^{-1} A_i ]$ and $\tr [S^{-1} A_k S^{-1} A_l S^{-1} A_i S^{-1} A_j]$ for all $i,j,k,l \in [m]$. For the former, we only need to sum up all diagonal terms of $S^{-1} A_k S^{-1} A_l S^{-1} A_i $, and each diagonal term is computed by multiplying one column of $S^{-1} A_k$ and one column of $S^{-1} A_l S^{-1} A_i$. Therefore, the cost of finding $\tr [S^{-1} A_k S^{-1} A_l S^{-1} A_i ]$ for all $i,k,l \in [m]$ is $m^3n^2$. For the latter, we only need to sum up all diagonal terms of $S^{-1} A_k S^{-1} A_l S^{-1} A_i S^{-1} A_j$, and each diagonal term is computed by multiplying one column of $S^{-1} A_k S^{-1} A_l $ and one column of $S^{-1} A_i S^{-1} A_j$. Therefore, the cost of finding $\tr [S^{-1} A_k S^{-1} A_l S^{-1} A_i S^{-1} A_j]$ for all $i,j,k,l \in [m]$ is $m^4 n^2$.\footnote{If we batch them together and use matrix multiplication, this term will be ${\cal T}_{\mathrm{mat}}(m^2, n^2, m^2)$.} In total, this step costs $m^4 n^2$ time.

\paragraph{Step 4.} 
We compute $\nabla \phi_{\mathrm{vol}}(y)$,  $Q(S)$, $R(S)$, and $T(S)$. We note 
\begin{align*}
(\A  (S^{-1}A_i S^{-1}\otimes S^{-1}) \mathsf{A}^\top)_{k,l} = &~\tr [A_k S^{-1} A_l S^{-1} A_i S^{-1}]\\
(\A  (S^{-1}\otimes S^{-1}A_i S^{-1}) \mathsf{A}^\top)_{k,l} = &~ \tr[A_kS^{-1}A_i S^{-1} A_l S^{-1}]
\end{align*}
and 
\begin{align*}
(\A  (S^{-1}A_i S^{-1} A_j S^{-1}\otimes S^{-1}) \mathsf{A}^\top)_{k,l} = &~ \tr [A_k S^{-1} A_l S^{-1} A_i S^{-1} A_j S^{-1}]\\
(\A  ( S^{-1} \otimes S^{-1}A_i S^{-1} A_j S^{-1}) \mathsf{A}^\top)_{k,l} = &~ \tr[A_kS^{-1}A_iS^{-1}A_jS^{-1}A_\ell S^{-1}] 
\end{align*}
and
\begin{align*}
    (\A (S^{-1}A_iS^{-1} \otimes S^{-1}A_jS^{-1}) \A^\top)_{k,l} = \tr[A_kS^{-1}A_jS^{-1} A_l S^{-1}A_iS^{-1}] 
\end{align*}
Thus each coordinate of $\nabla \phi_{\mathrm{vol}}(y)$,  $Q(S)$, $R(S)$, and $T(S)$ can be computed by multiplying the terms in the second step (i.e. $\tr [S^{-1} A_k S^{-1} A_l S^{-1} A_i ]$ and $\tr [S^{-1} A_k S^{-1} A_l S^{-1} A_i S^{-1} A_j]$ for all $i,j,k,l \in [m]$) with $H(S)^{-1}$, and then taking trace. These cost $O(m^{\omega+2})$ time, in total. 
More specifically, for the gradient $\nabla \phi_{\mathrm{vol}}(y)$,
\begin{align*}
    (\nabla \phi_{\mathrm{vol}}(y))_i = &~ -\tr[H(S)^{-1} \cdot \A  (S^{-1}A_i S^{-1}\otimes S^{-1}) \mathsf{A}^\top]\\
    = &~ \sum_{k=1}^m \sum_{l=1}^m -H(S)^{-1}_{k,l}\cdot\tr [A_k S^{-1} A_l S^{-1} A_i S^{-1}].
\end{align*}
Hence, it takes $O(m^3)$-time to compute $\nabla \phi_{\mathrm{vol}}(y)$. 

For $Q(S)$,
\begin{align*}
    Q(S)_{i,j} = &~ \tr[H(S)^{-1} \A (S^{-1}A_iS^{-1}A_jS^{-1} \otimes_S S^{-1}) \A^\top]\\
    = &~ \sum_{k=1}^m \sum_{l=1}^m -H(S)^{-1}_{k,l}\cdot \frac{1}{2}\left(\tr [A_k S^{-1} A_l S^{-1} A_i S^{-1} A_j S^{-1}] + \tr[A_kS^{-1}A_iS^{-1}A_jS^{-1}A_\ell S^{-1}] \right).
\end{align*}
Hence, it takes $O(m^4)$-time to compute $Q(S)$. 

For $R(S)$,
\begin{align*}
    R(S)_{i,j} = &~ \tr[H(S)^{-1} \A (S^{-1}A_iS^{-1} \otimes_S S^{-1}A_jS^{-1}) \A^\top]\\
    = &~ \sum_{k=1}^m \sum_{l=1}^m -H(S)^{-1}_{k,l}\cdot \frac{1}{2}\left(\tr[A_kS^{-1}A_jS^{-1} A_l S^{-1}A_iS^{-1}] + \tr[A_kS^{-1}A_iS^{-1} A_l S^{-1}A_jS^{-1}]  \right).
\end{align*}
Hence, it takes $O(m^4)$-time to compute $Q(S)$. For $T(S)$,
\begin{align*}
    T(S)_{i,j} = &~ \tr[H(S)^{-1} \A (S^{-1}A_iS^{-1} \otimes_S S^{-1}) \A^\top H(S)^{-1}\A (S^{-1}A_jS^{-1} \otimes_S S^{-1}) \A^\top ]\\
    = &~ \vect[\A (S^{-1}A_jS^{-1} \otimes_S S^{-1}) \A^\top]^\top (H(S)^{-1} \otimes H(S)^{-1}) \vect[\A (S^{-1}A_iS^{-1} \otimes_S S^{-1}) \A^\top].
\end{align*}
The matrix $\A (S^{-1}A_iS^{-1} \otimes_S S^{-1}) \A^\top$ can be computed in $O(m^2)$-time, via
\begin{align*}
    (\A (S^{-1}A_iS^{-1} \otimes_S S^{-1}) \A^\top)_{k,l} = \frac{1}{2}\left(\tr [A_k S^{-1} A_l S^{-1} A_i S^{-1}] + \tr [A_k S^{-1} A_i S^{-1} A_l S^{-1}]\right).
\end{align*}
Then, we vectorize this matrix and multiply with the Kronecker product $H(S)^{-1}\otimes H(S)^{-1}$. It takes $O(m^\omega)$-time to obtain the vector
\begin{align*}
    (H(S)^{-1} \otimes H(S)^{-1}) \vect[\A (S^{-1}A_iS^{-1} \otimes_S S^{-1}) \A^\top]\in \R^{m^2}.
\end{align*}
Next, we do the inner product and get $T(S)_{i,j}$ in $O(m^2)$-time. Thus, $T(S)$ can be computed in $O(m^{\omega + 2})$-time. Therefore, this step takes $m^{\omega+2}$ time in total.

\paragraph{Step 5.}
We compute $\delta_y \in \R^m$. Since 
\begin{align*}
\nabla^2 \phi(y) = 225\sqrt{\frac{n}{m}}\cdot \left(Q(S) + R(S) + T(S) + \frac{m-1}{n-1}\cdot H(S)\right),
\end{align*}
it can be found in $m^2$ time using the terms in the second and fourth step. Notice $\eta b - \nabla \phi(y) = \eta b - 225\sqrt{\frac{n}{m}}\cdot \left(\nabla \phi_{\mathrm{vol}}(y) + \frac{m-1}{n-1}\cdot \nabla \phi_{\mathrm{log}}(y)\right)$, it can be computed in $m$ time. Then, $\delta_y$ can be computed in $m^\omega$ time. In total, this step costs $m^\omega$ time.

Summing up, the total cost per iteration is given by 
\begin{align*}
    mn^2 + n^\omega m^2 + m^4 n^2 + n^\omega m^2 + m^{\omega+2} =  m^2 n^{\omega} + m^4 n^2 + m^{\omega+2}.
\end{align*}
\end{proof}

\section{Maintain the Leverage Score Matrix of the Volumetric Barrier}\label{sec:leverage_matrix}

It is observed in the LP that the volumetric barrier is roughly like the log barrier weighted by the leverage score of the matrix $A_x:=S^{-1}A$, which is defined to be the diagonal elements of the orthogonal projection matrix $A_x^\top (A_x^\top A_x)^\top A_x$. Similar phenomenon also appears in the SDP. we first consider the orthogonal projection matrix $P\in \R^{n^2\times n^2}$ on the image of $\mathsf{A}(S^{-1/2}\otimes S^{-1/2})$:
\begin{align*}
    P(S):=(S^{-1/2}\otimes S^{-1/2}) \mathsf{A}^\top \left(\mathsf{A} (S^{-1}\otimes S^{-1})\mathsf{A}^\top\right)^{-1} \mathsf{A} (S^{-1/2}\otimes S^{-1/2})\in \R^{n^2\times n^2}.
\end{align*}
Then, we define the leverage score matrix as the block-trace of $P(S)$:
\begin{definition}[Leverage score matrix]\label{def:vol_projection_matrix}
For $S\succ 0$, define the leverage score matrix $\Sigma(S)\in \R^{n\times n}$ as:
\begin{align*}
    \Sigma(S)_{i,j} := \tr\left[(e_i\otimes I_n)^\top P (e_j\otimes I_n)\right]~~~\forall i,j\in [n].
\end{align*}
\end{definition}

In this section, we show how to efficiently compute the leverage score matrix in each iteration of the IPM via low-rank update and amortization, which may be of independent interest.

\subsection{Basic facts on the leverage score matrix}\label{sec:fact_leverage}
\begin{fact}[Gradient of the volumetric barrier]
For $S\succ 0$, we have
\begin{align*}
    \nabla \phi_{\mathrm{vol}}(y) = -\mathsf{A} (S^{-1/2}\otimes S^{-1/2}) \mathrm{vec}(\Sigma).
\end{align*}
The computation cost is $mn^2 + n^{\omega}$.
\end{fact}
\begin{proof}

The computation cost is $mn^2+\mathcal{T}_{\mathrm{kron}}(n)$, where $\mathcal{T}_{\mathrm{kron}}(n)$ is the time to compute $(A\otimes B)v$ for $A,B\in \R^{n\times n}, v\in \R^{n^2}$.

By Fact~\ref{fac:kron_vec}, we have $\mathcal{T}_{\mathrm{kron}}(n) = n^{\omega}$.
\end{proof}

\begin{fact}[``Proxy'' Hessian of the volumetric barrier]
For $S\succ 0$, we have
\begin{align*}
    Q(S)=\mathsf{A}(S\otimes (S^{-1/2}\Sigma S^{-1/2}))\mathsf{A}^\top.
\end{align*}
\end{fact}

\begin{fact}[Trace of the leverage score matrix]
It holds that
\begin{align*}
    \tr[\Sigma] = m.
\end{align*}
\end{fact}
%The $\Sigma$ in Definition~\ref{def:vol_projection_matrix} is another fundamental matrix in SDP. In linear programming, it reduces to the diagonal matrix of leverage scores. In the IPM framework, we actually only need to compute $\Sigma(S)\in \R^{n\times n}$. 
\subsection{Efficient algorithm for the leverage score matrix}
This section shows how to maintain $\Sigma$ efficiently.

\begin{lemma}\label{lem:compute_leverage}
Let $ \Sigma(S)_{i,j} :=  \tr\left[(e_i\otimes I_n)^\top P(S) (e_j\otimes I_n)\right]$. $S^{-1/2} \in \R^{n \times n}$, $S^{-1} \in \R^{n \times n}$, and $H(S)^{-1} \in \R^{n^2 \times n^2}$ are known. Then we can compute $\Sigma(S) \in \R^{n \times n}$ it in
\begin{align*}
    \min \Big\{ n^4 + {\cal T}_{M(S)}, n^\omega m^2 \Big\}
\end{align*}
\end{lemma}
\begin{proof}
We provide two different approaches for computing the matrix $\Sigma(S)$.

\paragraph{Approach 1.}
Each entry of $\Sigma(S)$ can be expressed as follows:
\begin{align*}
    \Sigma(S)_{i,j} := &~ \tr\left[(e_i\otimes I_n)^\top P(S) (e_j\otimes I_n)\right]\\
    = &~ \tr\left[(e_i\otimes I_n)^\top (S^{-1/2}\otimes S^{-1/2}) \mathsf{A}^\top H(S)^{-1} \mathsf{A} (S^{-1/2}\otimes S^{-1/2}) (e_j\otimes I_n)\right]\\
    = &~ \tr\left[\left((S^{-1/2})_i^\top\otimes S^{-1/2}\right) \mathsf{A}^\top H(S)^{-1} \mathsf{A} \left((S^{-1/2})_j\otimes S^{-1/2}\right)\right]\\
    = &~ \tr\left[\left((S^{-1/2})_i^\top\otimes S^{-1/2}\right) M(S) \left((S^{-1/2})_j\otimes S^{-1/2}\right)\right]\\
    = &~ \tr\left[\left((S^{-1/2})_j\otimes S^{-1/2}\right)\left((S^{-1/2})_i^\top\otimes S^{-1/2}\right) M(S) \right]\\
    = &~ \tr\left[\left((S^{-1/2})_j(S^{-1/2})_i^\top \otimes S^{-1}\right) M(S) \right],
\end{align*}
where $M(S):=\mathsf{A}^\top H(S)^{-1} \mathsf{A} \in \R^{n^2 \times n^2}$.
Notice that
\begin{align*}
    &\left((S^{-1/2})_j(S^{-1/2})_i^\top \otimes S^{-1}\right) M(S) \\
    = &~ \left((S^{-1/2})_j(S^{-1/2})_i^\top \otimes S^{-1}\right) \begin{bmatrix}
    \vect[M(S)_{(1,1)}] & \vect[M(S)_{(1,2)}] & \cdots & \vect[M(S)_{(n,n)}]
    \end{bmatrix}\\
    =&~ \begin{bmatrix}
    \vect[S^{-1}M_{(1,1)}(S^{-1/2})_i (S^{-1/2})_j^\top] &
    \cdots& 
    \vect[S^{-1}M_{(n,n)}(S^{-1/2})_i (S^{-1/2})_j^\top]
    \end{bmatrix}\in \R^{n^2\times n^2},
\end{align*}
where $M_{(k,\ell)}\in \R^{n\times n}:=\mathrm{mat}[M(S)_{(k,\ell)}]$ is the matrix form of the $(k,\ell)$-th column of $M(S)\in \R^{n^2\times n^2}$ for $(k,\ell)\in [n]\times [n]$.
Hence,
\begin{align*}
    \Sigma(S)_{i,j} = &~
    \sum_{k=1}^n \sum_{\ell=1}^n \vect\left[S^{-1}M_{(k,\ell)}(S^{-1/2})_i(S^{-1/2})_j^\top\right]_{(k-1)n+\ell}\\ 
    = &~ \sum_{k=1}^n \sum_{\ell=1}^n \left(S^{-1}M_{(k,\ell)}(S^{-1/2})_i(S^{-1/2})_j^\top\right)_{\ell, k}\\
    = &~\sum_{k=1}^n \sum_{\ell = 1}^n e_\ell^\top \cdot S^{-1} M_{(k,\ell)} (S^{-1/2})_i (S^{-1/2})_j^\top \cdot  e_k\\
    = &~ \sum_{k=1}^n \sum_{\ell = 1}^n (S^{-1})_\ell ^\top\cdot  M_{(k,\ell)}\cdot  (S^{-1/2})_i \cdot (S^{-1/2})_{j,k}.
\end{align*}
Then, we get that 
\begin{align}
    \Sigma(S) = &~ \sum_{i=1}^n \sum_{j=1}^n e_ie_j^\top\sum_{k=1}^n \sum_{\ell = 1}^n (S^{-1})_\ell ^\top\cdot  M_{(k,\ell)}\cdot  (S^{-1/2})_i \cdot (S^{-1/2})_{j,k}\notag\\
    = &~ \sum_{i=1}^n \sum_{k=1}^n \sum_{\ell = 1}^n \left((S^{-1})_\ell ^\top\cdot  M_{(k,\ell)}\cdot  (S^{-1/2})_i\right) e_i \cdot \left(\sum_{j=1}^n (S^{-1/2})_{j,k}\cdot  e_j^\top\right)\notag\\
    = &~ \sum_{i=1}^n \sum_{k=1}^n \sum_{\ell = 1}^n \left((S^{-1})_\ell ^\top\cdot  M_{(k,\ell)}\cdot  (S^{-1/2})_i\right) e_i \cdot  (S^{-1/2})_{k}^\top\notag\\
    = &~ \sum_{k=1}^n \sum_{\ell = 1}^n \left(\sum_{i=1}^n (S^{-1})_\ell ^\top\cdot  M_{(k,\ell)}\cdot  (S^{-1/2})_i\cdot e_i\right) \cdot  (S^{-1/2})_{k}^\top\notag\\
    = &~ \sum_{k=1}^n \sum_{\ell = 1}^n \left((S^{-1})_\ell^\top \cdot M_{(k,\ell)}\cdot  S^{-1/2}\right)^\top \cdot (S^{-1/2})_k^\top\notag\\
    = &~ \sum_{k=1}^n \sum_{\ell = 1}^n S^{-1/2} \cdot M_{(k,\ell)}^\top\cdot (S^{-1})_\ell\cdot  (S^{-1/2})_k^\top\notag\\
    = &~ S^{-1/2} \cdot \sum_{k = 1}^n \sum_{\ell = 1}^n M_{(k,\ell)}^\top \cdot (S^{-1})_\ell \cdot e_k^\top\cdot  S^{-1/2}\label{eq:Sigma_S_submatrices}.
\end{align}

Therefore, once we have $S^{-1/2} \in \R^{n \times n}$, $S^{-1} \in \R^{n \times n}$, and $M(S) \in \R^{n^2 \times n^2}$, then the $\Sigma(S)$ can be computed exactly in $O(n^4)$-time using Eq.~\eqref{eq:Sigma_S_submatrices}. More specifically, 
\paragraph{Step 1.}
For $k,\ell\in [n]$, we form $M_{(k,\ell)}$ from $M(S)$, which takes $O(n^2)$-time.
\paragraph{Step 2.}
We compute the matrix
\begin{align*}
    W_{k,\ell}:=M_{(k,\ell)}^\top \cdot (S^{-1})_\ell \cdot e_k^\top \in \R^{n\times n}
\end{align*}
in $O(n^2)$ time. And it takes $O(n^4)$-time to compute $\{W_{k,\ell}\}_{k,\ell\in [m]}$.
\paragraph{Step 3.}
We sum all the $W_{k,\ell}$ together in $O(n^4)$-time. And 
\begin{align*}
    \Sigma(S) =S^{-1/2} \cdot \left(\sum_{k = 1}^n \sum_{\ell = 1}^n W_{k,\ell}\right)\cdot  S^{-1/2}, 
\end{align*}
which can be done in $O(n^\omega)$-time.

Hence, $\Sigma(S)$ can be computed in
\begin{align*}
    O\left(n^\omega + n^4 + {\cal T}_{M(S)}\right) = O\left(n^4 + {\cal T}_{M(S)}\right)
\end{align*}
time, where ${\cal T}_{M(S)}$ is the computation cost for computing $M(S)$.

\paragraph{Approach 2.}

Another approach for computing $\Sigma(S)$ can be done in $O\left(n^\omega m^2\right)$-time, without maintaining $M(S)$.

Recall that
\begin{align*}
    \Sigma(S)_{i,j}=&~ \tr\left[\left((S^{-1/2})_i^\top\otimes S^{-1/2}\right) \mathsf{A}^\top \cdot H(S)^{-1} \cdot  \mathsf{A} \left((S^{-1/2})_j\otimes S^{-1/2}\right)\right]\\
    = &~ \tr\left[M_{i}\cdot H(S)^{-1}\cdot M_{j}^\top\right]\\
    = &~ \left\langle M_j^\top M_{i}, H(S)^{-1}\right\rangle
\end{align*}
Consider the matrix $M_i$:
\begin{align*}
    M_i:=&~\left((S^{-1/2})_i^\top\otimes S^{-1/2}\right) \mathsf{A}^\top\\
    = &~ \left((S^{-1/2})_i^\top\otimes S^{-1/2}\right) \begin{bmatrix}
    \vect[A_1] & \vect[A_2]&\cdots & \vect[A_m]
    \end{bmatrix}\\
    = &~ \begin{bmatrix}
    S^{-1/2}A_1(S^{-1/2})_i&\cdots & S^{-1/2}A_m(S^{-1/2})_i
    \end{bmatrix}\in \R^{n\times m}.
\end{align*}
\iffalse
Similarly, the third matrix can be written as
\begin{align*}
    \mathsf{A}\left((S^{-1/2})_j\otimes S^{-1/2}\right)
    = ~ \begin{bmatrix} 
    \vect[S^{-1/2}A_1(S^{-1/2})_j]^\top\\\vdots \\ \vect[S^{-1/2}A_m(S^{-1/2})_j]^\top
    \end{bmatrix}\in \R^{m\times n^2}.
\end{align*}
\fi

For $k,l\in [m]$, the $(k,l)$-entry of $M_j^\top M_i$ is
\begin{align*}
    (S^{-1/2}A_k (S^{-1/2})_j)^\top\cdot  (S^{-1/2}A_l (S^{-1/2})_i)=&~ (S^{-1/2})_j^\top A_k S^{-1/2}\cdot S^{-1/2}A_l (S^{-1/2})_i\\
    = &~ (S^{-1/2})_j^\top A_k S^{-1}A_l (S^{-1/2})_i\\
    = &~ (S^{-1/2}A_kS^{-1}A_l S^{-1/2})_{j,i}\\
    = &~ (S^{-1/2}A_l S^{-1}A_k S^{-1/2})_{i,j}.
\end{align*}
Hence, $(i,j)$-entry of $\Sigma(S)$ is
\begin{align*}
    \Sigma(S)_{i,j} = &~ \sum_{k=1}^m \sum_{l=1}^m (S^{-1/2}A_lS^{-1}A_k S^{-1/2})_{i,j}\cdot  (H(S)^{-1})_{k,l}.
\end{align*}
It implies that
\begin{align*}
    \Sigma(S) = &~ \sum_{k=1}^m \sum_{l=1}^mS^{-1/2}A_lS^{-1}A_k S^{-1/2}\cdot  (H(S)^{-1})_{k,l}.
\end{align*}

We use the following steps to compute $\Sigma(S)$:
\paragraph{Step 1.} 
We compute $S^{-1/2}A_l S^{-1}A_k S^{-1/2}$ for all $k,l\in [m]$. It can be done in $O(n^\omega m^2)$-time.
\paragraph{Step 2.}
We compute $\Sigma(S)$ by summing the $m^2$ matrices with weights $(H(S)^{-1})_{k,l}$. It can be done in $O(m^2n^2)$-time.

Hence, it takes $O(n^\omega m^2)$-time in total to compute $\Sigma(S)$, assuming $H(S)^{-1}$ is given.
\end{proof}

\subsection{Maintain intermediate matrix}

The following lemma shows how to efficiently compute $M(S)$ in each iteration:
\begin{lemma}[Compute $M(S)$]\label{lem:maintain_M}
The matrix $M(S):=\mathsf{A}^\top H(S)^{-1}\mathsf{A}\in \R^{n^2\times n^2}$ can be computed as follows:
\begin{itemize}
\item {\bf Part 1.}
In the initialization, if $H(S)^{-1} \in \R^{m \times m}$ is already known, then $M(S)$ can be computed in ${\cal T}_{\mathrm{mat}}(n^2, m, n^2)$ time.

\item {\bf Part 2.}
In each iteration, $M(\wt{S})$ can be computed in ${\cal T}_{\mathrm{mat}}(n^2, n^2, nr_t)$ time.
\end{itemize}
\end{lemma}
\begin{proof}
{\bf Part 1.}
Given $H(S)^{-1}$, it costs ${\cal T}_{\mathrm{mat}}(n^2, m, n^2)$ to compute $M(S) = \A^\top H(S)^{-1} \A \in \R^{n^2 \times n^2}$.

{\bf Part 2.}
We can maintain $M(S)$ in each iteration. Let $\wt{S} \in \R^{n \times n}$ be the approximated slack variable in the previous iteration and $\wt{S}^{\new} \in \R^{n \times n}$ be the current approximated slack variable. Let $G:=H(\wt{S})^{-1} \in \R^{m \times m}$ and $G^{\new}:=H(\wt{S}^{\new})^{-1} \in \R^{m \times m}$. By Lemma~\ref{lem:close_form_hessian_inverse}, we have
\begin{align*}
    G^{\new} = G - G \cdot \mathsf{A} Y_1 \cdot (I + Y_2^{\top} \mathsf{A}^{\top} \cdot \mathsf{A} Y_1)^{-1} \cdot Y_2^{\top} \mathsf{A}^{\top} \cdot G
\end{align*}
Then, 
\begin{align*}
    M(\wt{S}^{\new}) = &~ \mathsf{A}^\top G^{\new} \mathsf{A}\\
    = &~ \mathsf{A}^\top G \mathsf{A} - \mathsf{A}^\top G \cdot \mathsf{A} Y_1 \cdot (I + Y_2^{\top} \mathsf{A}^{\top} \cdot \mathsf{A} Y_1)^{-1} \cdot Y_2^{\top} \mathsf{A}^{\top} \cdot G\mathsf{A}\\
    = &~ M(\wt{S}) - M(\wt{S})\cdot Y_1 \cdot (I + Y_2^{\top} \mathsf{A}^{\top} \cdot \mathsf{A} Y_1)^{-1} \cdot Y_2^{\top} \cdot M(\wt{S}),
\end{align*}
where $Y_1,Y_2\in \R^{n^2 \times (2 n r_t + r_t^2)}$ and $(I + Y_2^{\top} \mathsf{A}^{\top} \cdot \mathsf{A} Y_1)^{-1} \in \R^{(2n r_t + r_t^2)\times (2n r_t + r_t^2)}$. Hence, we first compute $M(\wt{S})Y_1\in \R^{n^2\times nr_t}$ and $Y_2^\top M(\wt{S})\in \R^{nr_t\times n^2}$ in ${\cal T}_{\mathrm{mat}}(n^2, n^2, nr_t)$. $M(\wt{S}^{\new}) \in \R^{n^2 \times n^2}$ can be directly computed from $G^{\new} \in \R^{m \times m}$ in ${\cal T}_{\mathrm{mat}}(n^2, nr_t, n^2)$-time.

\end{proof}

\subsection{Amortized running time}
\begin{theorem}
There is an algorithm that compute $\Sigma(S)$ in each iteration of Algorithm~\ref{alg:sdp_vol} with amortized cost-per-iteration
\begin{align*}
    \min \Big\{ n^{2\omega-\frac{1}{2}}, n^\omega m^2 \Big\}.
\end{align*}
\end{theorem}
\begin{proof}
By Lemma~\ref{lem:compute_leverage}, we know that the cost-per-iteration to compute $Sigma(S)$ is
\begin{align*}
    \min \Big\{ n^4 + {\cal T}_{M(S)}, n^\omega m^2 \Big\},
\end{align*}
where ${\cal T}_{M(S)} = {\cal T}_{\mathrm{mat}}(n^2, nr_t, n^2)$ by Lemma~\ref{lem:maintain_M}.

Then, by Corollary~\ref{cor:tmat_m2_nr_m2}, we have
\begin{align*}
    \sum_{t=1}^T {\cal T}_{\mathrm{mat}}(n^2, nr_t, n^2) = O^*\left(T\cdot n^{2\omega-\frac{1}{2}}\right).
\end{align*}
Therefore, the amortized running time per iteration is
\begin{align*}
    \min \Big\{ n^4 + n^{2\omega-\frac{1}{2}}, n^\omega m^2 \Big\} = \min \Big\{ n^{2\omega-\frac{1}{2}}, n^\omega m^2 \Big\}.
\end{align*}
\end{proof}

\begin{remark}
The second term $n^{\omega}m^2$ in the running time represents the approach that does not use the maintenance technique. When $m=\Omega(n^{0.94})$, the first approach using low-rank update and amortization is faster.
\end{remark}

\addcontentsline{toc}{section}{References}
\bibliographystyle{alpha}
\bibliography{ref}
\end{document}